\documentclass[hidelinks,onefignum,onetabnum]{siamart220329}
\usepackage{amssymb}
\usepackage{amsmath}
\usepackage{enumitem}
\usepackage{bm}
\usepackage{color}
\usepackage{exscale}
\usepackage{relsize}
\usepackage{graphicx}
\usepackage{wrapfig}
\usepackage{multicol}
\usepackage{verbatim}
\newtheorem{rmk}{Remark}
\usepackage{appendix}

\usepackage{lipsum}
\usepackage{amsfonts}
\usepackage{graphicx}
\usepackage{epstopdf}
\usepackage{algorithmic}
\ifpdf
  \DeclareGraphicsExtensions{.eps,.pdf,.png,.jpg}
\else
  \DeclareGraphicsExtensions{.eps}
\fi

\newsiamremark{remark}{Remark}
\newsiamremark{hypothesis}{Hypothesis}
\crefname{hypothesis}{Hypothesis}{Hypotheses}
\newsiamthm{claim}{Claim}

\headers{Multi-scale method for nonlinear thermo-mechanical problems}{H. Dong, Y. F. Nie, X. F. Guan and J. Z. Cui }

\title{Higher-order multi-scale method for high-accuracy nonlinear thermo-mechanical simulation of heterogeneous shells\thanks{Submitted to the editors DATE.
\funding{This work was funded by the the National Natural Science Foundation of China (No.\hspace{1mm}12001414), Young Talent Fund of Association for Science and Technology in Shaanxi, China (No.\hspace{1mm}20220506), Young Talent Fund of Association for Science and Technology in Xi'an, China (No.\hspace{1mm}095920221338), the Fundamental Research Funds for the Central Universities (No.\hspace{1mm}KYFZ23020), the National Natural Science Foundation of China (Nos.\hspace{1mm}11971386 and 51739007).}}}

\author{Hao Dong\thanks{School of Mathematics and Statistics, Xidian University, Xi'an 710071, PR China
  (\email{donghaoxd@xidian.edu.cn}).}
\and Xiaofei Guan\thanks{School of Mathematical Sciences, Tongji University, Shanghai 200092, PR China (\email{guanxf@tongji.edu.cn}).}
\and Yufeng Nie\thanks{School of Mathematics and Statistics, Northwestern Polytechnical University, Xi'an 710129, PR China (\email{yfnie@nwpu.edu.cn}).}
}

\usepackage{amsopn}

\ifpdf
\hypersetup{
  pdftitle={Higher-order multi-scale method for nonlinear thermo-mechanical simulation},
  pdfauthor={H. Dong, X. F. Guan and Y. F. Nie, }
}
\fi

\begin{document}

\maketitle

\begin{abstract}
In the present work, we consider multi-scale computation and convergence for nonlinear time-dependent thermo-mechanical equations of inhomogeneous shells possessing temperature-dependent material properties and orthogonal periodic configurations. The first contribution is that a novel higher-order macro-micro coupled computational model is rigorously devised via multi-scale asymptotic technique and Taylor series approach for high-accuracy simulation of heterogeneous shells. Benefitting from the higher-order corrected terms, the higher-order multi-scale computational model keeps the conservation of local energy and momentum for nonlinear thermo-mechanical simulation. Moreover, a global error estimation with explicit rate of higher-order multi-scale solutions is first derived in the energy norm sense. Furthermore, an efficient space-time numerical algorithm with off-line and on-line stages is presented in detail. Adequate numerical experiments are conducted to confirm the competitive advantages of the presented multi-scale approach, exhibiting not only the exceptional numerical accuracy, but also the less computational expense for heterogeneous shells.
\end{abstract}



\section{Introduction}
In aerospace and aviation fields, various heterogeneous materials are often manufactured as multi-layer, sandwich, reinforced plate and shell structures. These heterogeneous plate and shell structures are extensively employed in spacecraft, aircraft and space structures, which always served under extreme thermo-mechanical coupled environments and must survive high thermo-mechanical coupled loads \cite{R6}. According to the existing researches \cite{R9,R10,R13,R14,R15,R16}, the material properties will change significantly with temperature variation and exhibit powerful dependency on temperature in high temperature environments, which necessarily enables us to take into account the prominent nonlinearities of material behaviors.

Specifically, on the basis of classical elastic shell and thermo-mechanical multiphysics theories \cite{R1,R2,R3,R4,R5}, we establish time-dependent governing equations for nonlinear thermo-mechanical problems of inhomogeneous shells with temperature-dependent properties and orthogonal periodic configurations as follows
\begin{equation}
\left\{
\begin{aligned}
&{{\rho ^\xi }({\bm{\alpha }},T^{\xi})\frac{\partial^2 u_i^\xi({\bm{\alpha }},t)}{\partial t^2}}-\sum\limits_{j = 1}^3 {\frac{1}{H}\frac{\partial }{{\partial {\alpha _j}}}\left( {\frac{H}{{{H_j}}}\sigma _{ij}^\xi({\bm{\alpha }},t) } \right)}\\
&\quad\quad\quad\quad\quad  -\sum\limits_{j \ne i,j = 1}^3 {\frac{1}{{{H_i}{H_j}}}\frac{{\partial {H_i}}}{{\partial {\alpha _j}}}} \sigma _{ij}^\xi({\bm{\alpha }},t)+ \sum\limits_{j \ne i,j = 1}^3 {\frac{1}{{{H_i}{H_j}}}\frac{{\partial {H_j}}}{{\partial {\alpha _i}}}} \sigma _{jj}^\xi({\bm{\alpha }},t) \\
&\quad\quad\quad\quad\quad  = {f_i}({\bm{\alpha }},t),\;\;\text{in}\;\;\Omega\times(0,\mathcal T),\;\;i=1,2,3,\\
&{{\rho ^\xi }({\bm{\alpha }},T^{\xi}){c^\xi }({\bm{\alpha }},T^{\xi})\frac{{\partial {T^\xi({\bm{\alpha }},t) }}}{{\partial t}}{\rm{ + }}}\sum\limits_{i = 1}^3 {\frac{1}{H}\frac{\partial }{{\partial {\alpha _i}}}\left( {\frac{H}{{{H_i}}}q_i^\xi({\bm{\alpha }},t) } \right)} \\
&\quad\quad\quad\quad\quad{ + \vartheta _{ij}^\xi({\bm{\alpha }},T^{\xi}) \frac{{\partial \varepsilon _{ij}^\xi({\bm{\alpha }},t) }}{{\partial t}}}= h({\bm{\alpha }},t),\;\;\text{in}\;\;\Omega\times(0,\mathcal T),\\
&\bm{u}^{\xi}({\bm{\alpha }},t)=\widehat{\bm{u}}({\bm{\alpha }},t),\;\;\text{on}\;\;\partial\Omega_{u}\times(0,\mathcal T),\\
&T^{\xi}({\bm{\alpha }},t)=\widehat{T}({\bm{\alpha }},t),\;\;\text{on}\;\;\partial\Omega_{T}\times(0,\mathcal T),\\
&\bm{u}^\xi (\bm{\alpha },0)=\widetilde{\bm{u}}^0,\;\;\frac{\partial \bm{u}^\xi (\bm{\alpha },t)}{\partial t}\big|_{t=0}=\widetilde{\bm{u}}^1 (\bm{\alpha }),\;\;T^{\xi}(\bm{\alpha },0)=\widetilde T,\;\;\text{in}\;\;\Omega.
\end{aligned} \right.
\end{equation}
In above multi-scale nonlinear initial-boundary system (1.1), the displacement field $\bm{u}^{\xi}({\bm{\alpha }},t)=(u_1^{\xi}({\bm{\alpha }},t),u_2^{\xi}({\bm{\alpha }},t),u_3^{\xi}({\bm{\alpha }},t))^T$ and temperature field $T^{\xi}({\bm{\alpha }},t)$ are the targeted and undetermined, that are embedded in the stresses $\sigma_{ij}^\xi({\bm{\alpha }},t)$, strains ${\varepsilon}_{ij}^{\xi}({\bm{\alpha }},t)$ and heat flux $q_i^\xi({\bm{\alpha }},t)$. The equation coefficient $H=H_1H_2H_3$ and $H_i$ is the Lam$\rm{\acute{e}}$ coefficient in $\alpha_i$ direction of orthogonal curvilinear coordinates \cite{R1}. $\Omega$ is a bounded convex domain in $\mathbb{R}^3$ with a boundary $\partial\Omega=\partial\Omega_{u}\cup\partial\Omega_{T}$. $\widehat{\bm{u}}(\bm{\alpha },t)$ and $\widehat{T}(\bm{\alpha },t)$ are, respectively, the prescribed displacement and temperature on the domain boundaries $\partial\Omega_{u}$ and $\partial\Omega_{T}$. $\widetilde{\bm{u}}^0$, $\widetilde{\bm{u}}^1(\bm{\alpha })$ and $\widetilde{T}$ are, respectively, the initial displacement, velocity and temperature when inhomogeneous shells are in initial stress-free state. $\xi$ denotes the characteristic length of periodic unit cell. $\displaystyle\rho^\xi ({\bm{\alpha }},T^{\xi})$, $\displaystyle c^\xi ({\bm{\alpha }},T^{\xi})$ and $\{\vartheta _{ij}^\xi({\bm{\alpha }},T^{\xi})\}$ represent the mass density, specific heat and two-way thermo-mechanical effect tensor. $f_i(\bm{\alpha },t)$ and $h(\bm{\alpha },t)$ denote the body force and internal heat source.

Furthermore, the constitutive law of multi-scale nonlinear problem (1.1) takes the form as below
\begin{equation}
\sigma _{ij}^{\xi}({\bm{\alpha }},t) = {a_{ijkl}^\xi}({\bm{\alpha }},T^{\xi})\varepsilon _{kl}^{\xi}({\bm{\alpha }},t) - {b _{ij}^\xi}({\bm{\alpha }},T^{\xi})( {{T^{\xi}}({\bm{\alpha }},t) - \widetilde T} ).
\end{equation}
The heat flux $q_i^\xi$ in multi-scale nonlinear problem (1.1) takes the expression as below
\begin{equation}
q_i^\xi({\bm{\alpha }},t) =  - \sum\limits_{j = 1}^3 {k_{ij}^\xi({\bm{\alpha }},T^{\xi}) \frac{1}{{{H_j}}}} \frac{{\partial {T^\xi({\bm{\alpha }},t) }}}{{\partial {\alpha _j}}}.
\end{equation}
In the above formulas (1.2) and (1.3), $\{
a_{ijkl}^\xi ({\bm{\alpha }},T^{\xi})\}$, $\{b _{ij}^\xi ({\bm{\alpha }},T^{\xi})\}$ and $\{k_{ij}^\xi ({\bm{\alpha }},T^{\xi})\}$ denote the tensor coefficients of elastic modulus, thermal modulus and thermal conductivity of inhomogeneous shells separately. Additionally, in the governing equations (1.1) and constitutive law (1.2), the strain field $\bm{\varepsilon}^{\xi}$ is presented with regard to the displacement field $\bm{u}^{\xi}$ as below \cite{R1}
\begin{equation}
\begin{aligned}
&\varepsilon _{11}^\xi({\bm{\alpha }},t)  =\varepsilon _{11}^\xi(\bm{u}^{\xi}({\bm{\alpha }},t))  = \frac{1}{{{H_1}}}\frac{{\partial u_1^\xi }}{{\partial {\alpha _1}}} + \frac{1}{{{H_1}{H_2}}}\frac{{\partial {H_1}}}{{\partial {\alpha _2}}}u_2^\xi  + \frac{1}{{{H_1}{H_3}}}\frac{{\partial {H_1}}}{{\partial {\alpha _3}}}u_3^\xi ,\\
&\varepsilon _{22}^\xi({\bm{\alpha }},t)  =\varepsilon _{22}^\xi(\bm{u}^{\xi}({\bm{\alpha }},t))  =  \frac{1}{{{H_2}}}\frac{{\partial u_2^\xi }}{{\partial {\alpha _2}}} + \frac{1}{{{H_2}{H_3}}}\frac{{\partial {H_2}}}{{\partial {\alpha _3}}}u_3^\xi  + \frac{1}{{{H_2}{H_1}}}\frac{{\partial {H_2}}}{{\partial {\alpha _1}}}u_1^\xi ,\\
&\varepsilon _{33}^\xi({\bm{\alpha }},t)  = \varepsilon _{33}^\xi(\bm{u}^{\xi}({\bm{\alpha }},t))  = \frac{1}{{{H_3}}}\frac{{\partial u_3^\xi }}{{\partial {\alpha _3}}} + \frac{1}{{{H_3}{H_1}}}\frac{{\partial {H_3}}}{{\partial {\alpha _1}}}u_1^\xi  + \frac{1}{{{H_3}{H_2}}}\frac{{\partial {H_3}}}{{\partial {\alpha _2}}}u_2^\xi,\\
&\varepsilon _{ij}^\xi({\bm{\alpha }},t)  = \varepsilon _{ij}^\xi(\bm{u}^{\xi}({\bm{\alpha }},t))  = \frac{1}{2}\bigg[\frac{{{H_j}}}{{{H_i}}}\frac{\partial }{{\partial {\alpha _i}}}\Big( {\frac{{u_j^\xi }}{{{H_j}}}} \Big) + \frac{{{H_i}}}{{{H_j}}}\frac{\partial }{{\partial {\alpha _j}}}\Big( {\frac{{u_i^\xi }}{{{H_i}}}} \Big)\bigg],\;\;i\neq j.
\end{aligned}
\end{equation}
Mathematical model (1.1) appears in the mechanical and thermal coupling simulation of heterogeneous shells with temperature-dependent properties, whose multi-scale property originates from the high-frequency oscillatory equation coefficients due to orthogonal periodic configurations when $0<\xi\ll1$, which lead to refined discretization for high-fidelity simulation and make the overall cost prohibitively high, especially for dynamic time-dependent problems. In addition, the Laplace transform technique for linear thermo-mechanical system in \cite{R4,R54,R55} can not work for numerically simulating and theoretically analyzing the investigated nonlinear thermo-mechanical system due to its nonlinear coefficients.

To resolve the challenging issue of numerical computation for heterogeneous materials and structures, researchers have established many kinds of multi-scale approaches which possessing the efficiency of macro-scale model and the accuracy of the micro-scale model, such as the asymptotic homogenization method (AHM) \cite{R18,R19}, multi-scale finite element method (MsFEM) \cite{R22}, heterogeneous multi-scale method (HMM) \cite{R20}, variational multi-scale method (VMS) \cite{R21}, multi-scale eigenelement method (MEM) \cite{R23}, localized orthogonal decomposition method (LOD) \cite{R24} and finite volume based asymptotic homogenization theory (FVBAHT) \cite{R25}, etc. In recent years, Cui and his research group systematically presented a class of higher-order multi-scale approach to efficiently compute the thermo-mechanical coupling behaviors of inhomogeneous materials and structures, as pointed out in the literatures \cite{R27,R28,R29,R30}. Through introducing the enriched higher-order correctors, these novel multi-scale methods are capable of precisely capturing the microscopic fluctuation information inside inhomogeneous medias than conventional homogenized approaches and lower-order multi-scale approaches. However, it should be noted that the practical applications for above-mentioned multi-scale approaches are often restricted by the hypotheses on heterogeneous periodic configurations in cartesian coordinate system \cite{R33,R34,R35,R36}, etc. In practical engineering applications, heterogeneous shells are extensively employed in plenty of engineering areas, which possesses periodic configurations in non-cartesian coordinate system. Inevitably, research attention of scientists and engineers is focused on multi-scale modeling and computation for inhomogeneous solids with non-cartesian periodic configurations. As pointed out in references \cite{R39,R40,R42}, multi-scale methods are presented to compute thermo-mechanical problems of inhomogeneous structures with microscopic periodic configurations in cylindrical coordinates. Additionally, in the literatures \cite{R43,R44,R45,R27,R46}, multi-scale approaches are developed for simulating thermal, mechanical and thermo-mechanical problems of inhomogeneous structures with microscopic periodicity in general curvilinear coordinates, respectively.

Concerning nonlinear simulations of inhomogeneous materials and structures with temperature-dependent properties, this nonlinear thermo-mechanical coupled effect has a great impact on the deformation and failure of materials and structures. Consequently, continued research is of great significance and an important task for accurately and efficiently nonlinear thermo-mechanical coupling simulation and prediction, which is the premise of vibration control and optimal design of heterogeneous shell structures \cite{R16}. In references \cite{R48,R49,R50,R52,R53}, researchers employed asymptotic homogenization method to analyze and simulate static and transient nonlinear heat conduction problems of inhomogeneous solids with temperature-dependent conductivity coefficients. In \cite{R51}, Pankov systematically studied the G-convergence and homogenized theory of nonlinear partial differential operators. Moreover, researchers in \cite{R9} developed a multi-scale micromechanical-structural approach to simulate coupling heat transfer and deformations of functionally graded materials (FGM) with temperature-dependent material properties. In \cite{R10}, a variational asymptotic method was presented by Teng et al. for simulating inhomogeneous materials comprising temperature-dependent constituents. In \cite{R13}, researchers presented a multi-scale model for the thermal shock resistance prediction of porous ceramics with temperature-dependent material performances. Skinner and Chattopadhyay in \cite{R14} established a multi-scale thermo-mechanical computation framework to simulate the temperature-dependent damage behaviors of ceramic composites. In \cite{R5}, a pre-correction multi-scale technique is developed by Zhang et al. based on coupling extended multi-scale finite element method (CEMsFEM) for nonlinear thermo-elastic simulation of inhomogeneous medias with temperature-dependent properties. To sum up, few researches have been published in which the nonlinear thermo-mechanical properties are considered in the multi-scale modeling and computation of inhomogeneous medias with temperature-dependent properties. However, extensively engineering demands strongly prompt the continued research about this challenging issue.

The target of the present study is to devise a novel higher-order multi-scale computational approach for effective nonlinear thermo-mechanical simulation of heterogeneous shells with temperature-dependent material performances, especially for large-scale heterogeneous shells that can not be numerically computed by direct numerical simulation (DNS). The way out of computational and theoretical analysis difficulties is that we innovatively devise a novel higher-order macro-micro coupled computational model via multi-scale asymptotic technique and Taylor series approach to precisely analyze nonlinear thermo-mechanical behaviors of heterogeneous shells. On the other hand, by imposing the homogeneous Dirichlet condition on auxiliary cell functions and under some simplifications, the explicit convergence rate of the higher-order multi-scale solutions is successfully derived since the higher-order multi-scale solutions will automatically fulfil the boundary condition of original governing equations. By comparison to direct numerical simulation of the nonlinear thermo-mechanical equations of inhomogeneous shells, the proposed multi-scale computational approach and corresponding space-time numerical algorithm only need to solve a homogenization problem with non-discontinuous and non-oscillatory coefficient on macro-scale and several small-scale auxiliary cell problems on micro-scale. Hence, the space-time algorithm we established can greatly economize computational overhead both in computer storage space and CPU time, so as to quickly and accurately simulate the nonlinear thermo-mechanical behaviors of inhomogeneous shells. Additionally, we should highlight that the multi-scale computational models for heterogeneous solids possessing microscopic periodicity in cartesian, cylindrical and spherical coordinates can be derived without difficulty from the multi-scale computational model for heterogeneous shells in the present work.

The remainder of this article is organized as below. Section 2 devotes to construct the higher-order macro-micro coupled computational model for nonlinear thermo-mechanical coupling simulations of heterogeneous shells with orthogonal periodic configurations using multi-scale asymptotic technique and Taylor series approach. In Section 3, the local and global error theories are established for the proposed multi-scale solutions. Section 4 provides a space-time numerical algorithm with off-line microscopic computation and on-line macroscopic computation to effectively simulate nonlinear thermo-mechanical problems of heterogeneous shells at length. Extensive numerical examples including distinct types of heterogeneous shells motivated by engineering applications are designed in Section 5 to demonstrate the validity of the presented multi-scale computational model and corresponding numerical algorithm. Eventually, research conclusions and future directions are declared in Section 6. Throughout this text for readability, Einstein summation convention is employed to simplify repetitious indices apart from those repetitious indices with summation notations.
\section{Higher-order multi-scale computational model of nonlinear shells}
\subsection{The setting of multi-scale nonlinear system}
Within the framework of the asymptotic homogenization method, let us define $\displaystyle\bm{\beta}={\bm{\alpha}}/{\xi}=({\alpha_1}/{\xi},{\alpha_2}/{\xi},{\alpha_3}/{\xi})=(\beta_1,\beta_2,\beta_3)$ as microscopic coordinates of PUC $\Theta=(0,1)^3$. Whereupon material parameters $\rho^\xi ({\bm{\alpha }},T^{\xi})$, $a_{ijkl}^\xi ({\bm{\alpha }},T^{\xi})$, $b _{ij}^\xi ({\bm{\alpha }},T^{\xi})$, $c^\xi ({\bm{\alpha }},T^{\xi})$, $k_{ij}^\xi ({\bm{\alpha }},T^{\xi})$ and $\vartheta_{ij}^\xi ({\bm{\alpha }},T^{\xi})$ in nonlinear thermo-mechanical system (1.1) can be expressed with new forms $\rho({\bm{\beta }},T^{\xi})$, ${a_{ijkl}}({\bm{\beta }},T^{\xi})$, ${b _{ij}}({\bm{\beta }},T^{\xi})$, $c({\bm{\beta }},T^{\xi})$, ${k_{ij}}({\bm{\beta }},T^{\xi})$ and ${\vartheta_{ij}}({\bm{\beta }},T^{\xi})$. These new expressions imply that $\rho({\bm{\beta }},T^{\xi})$, $a_{ijkl}({\bm{\beta }},T^{\xi})$, $b_{ij}({\bm{\beta }},T^{\xi})$, $c({\bm{\beta }},T^{\xi})$, $k_{ij}({\bm{\beta }},T^{\xi})$ and $\vartheta_{ij}({\bm{\beta }},T^{\xi})$ are 1-periodic functions in variable ${\bm{\beta }}$.

Referring to \cite{R4,R27,R28,R29,R30}, some hypotheses of equation coefficients in multi-scale nonlinear equations (1.1) are presented as below.
\begin{enumerate}
\item[(A)]
$a_{ijkl}^{\xi}({\bm{\alpha }},T^{\xi})$, $b_{ij}^{\xi}({\bm{\alpha }},T^{\xi})$, $k_{ij}^{\xi}({\bm{\alpha }},T^{\xi})$ and $\vartheta_{ij}^{\xi}({\bm{\alpha }},T^{\xi})$ are symmetrical, and there exist two positive constants $\gamma_0$ and $\gamma_1$ irrespective of $\xi$ for the following uniform elliptic conditions
\begin{displaymath}
\begin{aligned}
&a_{ijkl}^{\xi}=a_{ijlk}^{\xi}=a_{klij}^{\xi},\gamma_0\eta_{ij}\eta_{ij}\leq a_{ijkl}^\xi({\bm{\alpha }},T^{\xi})\eta_{ij}\eta_{kl} \le\gamma_1\eta_{ij}\eta_{ij},\\
&b_{ij}^{\xi}=b_{ji}^{\xi},\gamma_0|\bm{\zeta}|^2\leq b_{ij}^\xi({\bm{\alpha }},T^{\xi})\zeta_i\zeta_j  \le\gamma_1|\bm{\zeta}|^2,\\
&k_{ij}^{\xi}=k_{ji}^{\xi},\gamma_0|\bm{\zeta}|^2\leq k_{ij}^\xi({\bm{\alpha }},T^{\xi})\zeta_i\zeta_j  \le\gamma_1|\bm{\zeta}|^2,\\
&\vartheta_{ij}^{\xi}=\vartheta_{ji}^{\xi},\gamma_0|\bm{\zeta}|^2\leq \vartheta_{ij}^\xi({\bm{\alpha }},T^{\xi})\zeta_i\zeta_j  \le\gamma_1|\bm{\zeta}|^2,
\end{aligned}
\end{displaymath}
where $\{\eta_{ij}\}$ is an arbitrary symmetrical matrix in $\mathbb{R}^{3\times 3}$, $\bm{\zeta}=(\zeta_1,\zeta_2,\zeta_3)$ is an arbitrary vector with
real elements in $\mathbb{R}^3$, and ${\bm{\alpha }}$ is an arbitrary point in $\Omega$.
\item[(B)]
$\rho^\xi ({\bm{\alpha }},T^{\xi})$, $a_{ijkl}^{\xi}({\bm{\alpha }},T^{\xi})$, $b_{ij}^{\xi}({\bm{\alpha }},T^{\xi})$, $ c^\xi ({\bm{\alpha }},T^{\xi})$, $k^\xi_{ij}({\bm{\alpha }},T^{\xi})$ and $\vartheta_{ij}^{\xi}({\bm{\alpha }},T^{\xi})\in L^\infty (\Omega)$; $0<\rho^{*}\leq\rho^\xi ({\bm{\alpha }},T^{\xi})$, $0<c^{*}\leq c^\xi ({\bm{\alpha }},T^{\xi})$, where $\rho^{*}$ and $c^{*}$ are constants irrespective of $\xi$.
\item[(C)]
$f_{i}\in L^2(\Omega\times(0,\mathcal T))$, $h\in L^2(\Omega\times(0,\mathcal T))$, $\widehat{\bm{u}}({\bm{\alpha }},t)\in L^2(0,\mathcal T;H^{1}(\Omega))^{3}$, $\widehat{T}({\bm{\alpha }},t)\in L^2(0,\mathcal T;H^{1}(\Omega))$, $\bar{\sigma}_i(\bm{\alpha },t)\in L^2(\Omega\times(0,\mathcal T))$, $\bar{q}(\bm{\alpha },t)\in L^2(\Omega\times(0,\mathcal T))$, $\widetilde{\bm{u}}^1(\bm{\alpha })\in ({L^2 (\Omega)})^{3}$.
\end{enumerate}
\subsection{Higher-order multi-scale analysis for nonlinear model problem}
In order to establish higher-order multi-scale computational model, we firstly assume that the basic physical quantities $\sigma _{ij}^\xi(\bm{\alpha },t)$ and $q_i^\xi(\bm{\alpha },t)$ are expanded as the following multi-scale asymptotic forms enlightened by \cite{R27,R45,R46}.
\begin{equation}
\left\{\begin{array}{l}
\varepsilon _{ij}^\xi(\bm{\alpha },t)  = {\xi ^{ - 1}}\varepsilon _{ij}^{[ - 1]}(\bm{\alpha },\bm{\beta },t) + {\xi ^0}\varepsilon _{ij}^{[0]}(\bm{\alpha },\bm{\beta },t) + {\xi ^1}\varepsilon _{ij}^{[1]}(\bm{\alpha },\bm{\beta },t) + {\rm O}({\xi ^2}),\\
\sigma _{ij}^\xi(\bm{\alpha },t)  = {\xi ^{ - 1}}\sigma _{ij}^{[ - 1]}(\bm{\alpha },\bm{\beta },t) + {\xi ^0}\sigma _{ij}^{[0]}(\bm{\alpha },\bm{\beta },t) + {\xi ^1}\sigma _{ij}^{[1]}(\bm{\alpha },\bm{\beta },t) + {\rm O}({\xi ^2}),\\
q_i^\xi(\bm{\alpha },t)  = {\xi ^{ - 1}}q_i^{[ - 1]}(\bm{\alpha },\bm{\beta },t) + {\xi ^0}q_i^{[0]}(\bm{\alpha },\bm{\beta },t) + {\xi ^1}q_i^{[1]}(\bm{\alpha },\bm{\beta },t) + {\rm O}({\xi ^2}).
\end{array}\right.
\end{equation}

Next, the physical quantities $u_i^\xi(\bm{\alpha },t)$ and ${T^\xi }(\bm{\alpha },t)$ are assumed to be expanded as the succeeding multi-scale expansion forms.
\begin{equation}
\left\{\begin{array}{l}
u_i^\xi(\bm{\alpha },t)= u_i^{[0]}(\bm{\alpha },\bm{\beta },t) + \xi u_i^{[1]}(\bm{\alpha },\bm{\beta },t) + {\xi ^2}u_i^{[2]}(\bm{\alpha },\bm{\beta },t) + {\rm O}({\xi ^3}),\\
{T^\xi }(\bm{\alpha },t) = {T^{[0]}}(\bm{\alpha },\bm{\beta },t) + \xi {T^{[1]}}(\bm{\alpha },\bm{\beta },t) + {\xi ^2}{T^{[2]}}(\bm{\alpha },\bm{\beta },t) + {\rm O}({\xi ^3}),
\end{array}\right.
\end{equation}
where $u_i^{[0]}$ and $T^{[0]}$ are zeroth-order expansion terms, $u_i^{[1]}$ and $T^{[1]}$ are first-order asymptotic terms (lower-order asymptotic terms), $u_i^{[2]}$ and $T^{[2]}$ are second-order asymptotic terms (higher-order asymptotic terms).

Then, we should introduce the key idea to handle the temperature-dependent material parameters. By means of the Taylor's formula and multi-index notation $\displaystyle\mathbf{D}^{(\alpha_1,\cdots,\alpha_m)}f(x_1,\cdots,x_m)=\frac{\partial^{\alpha} f(x_1,\cdots,x_m)}{\partial x_1^{\alpha_1}\cdots\partial x_m^{\alpha_m}}$ with $\alpha=\alpha_1+\cdots+\alpha_m$, material parameter $a_{ijkl}^\xi ({\bm{\alpha}},{T^\xi })$ is strongly reliant on temperature $T^\varepsilon$ and can be expanded around leading term $T^{[0]}$ as below
\begin{equation}
\begin{aligned}
&{a_{ijkl}}({\bm{\beta}},{T^{\xi}})= {a_{ijkl}}({\bm{\beta}},{T^{[0]}} + \xi {T^{[1]}} + {\xi ^2}{T^{[2]}} + {\rm O}({\xi ^3}))\\
& = {a_{ijkl}}({\bm{\beta}},{T^{[0]}}) + \mathbf{D}^{(0,1)}{a_{ijkl}}({\bm{\beta}},{T^{[0]}})\big[ {\xi {T^{[1]}} + {\xi ^2}{T^{[2]}} + {\rm O}({\xi ^3})} \big]\\
&+ \frac{1}{2}\mathbf{D}^{(0,2)}{a_{ijkl}}({\bm{\beta}},{T^{[0]}}){\big[ {\xi {T^{[1]}} + {\xi ^2}{T^{[2]}} + {\rm O}({\xi ^3})} \big]^2} +  \cdots \\
& = {a_{ijkl}}({\bm{\beta}},{T^{[0]}}) + \xi {T^{[1]}}\mathbf{D}^{(0,1)}{a_{ijkl}}({\bm{\beta}},{T^{[0]}})\\
&+ {\xi ^2}\big[ {{T^{[2]}}\mathbf{D}^{(0,1)}{a_{ijkl}}({\bm{\beta}},{T^{[0]}}) + \frac{1}{2}{{( {{T^{[1]}}} )}^2}\mathbf{D}^{(0,2)}{a_{ijkl}}({\bm{\beta}},{T^{[0]}})} \big] + {\rm O}({\xi ^3})\\
&= a_{ijkl}^{[0]}({\bm{\beta}},{T^{[0]}}) + \xi a_{ijkl}^{[1]}({\bm{\alpha}},{\bm{\beta}},{T^{[0]}}) + {\xi ^2}a_{ijkl}^{[2]}({\bm{\alpha}},{\bm{\beta}},{T^{[0]}}) + {\rm O}({\xi ^3}).
\end{aligned}
\end{equation}
Employing the aforementioned expanding approach as (2.3), remaining material parameters ${\rho^\xi}({\bm{\alpha}},{T^\xi })$, $b _{ij}^\xi ({\bm{\alpha}},{T^\xi })$, $c^\xi({\bm{\alpha}},{T^\xi })$, $k_{ij}^\xi ({\bm{\alpha}},{T^\xi })$ and $\vartheta_{ij}^\xi({\bm{\alpha}},{T^\xi })$ can be continuously expanded as below
\begin{equation}
\begin{aligned}
{\rho}({\bm{\beta}},{T^{\xi}})&= {\rho^{[0]}}({\bm{\beta}},{T^{[0]}}) + \xi {\rho^{[1]}}({\bm{\alpha}},{\bm{\beta}},{T^{[0]}}) + {\xi ^2}{\rho^{[2]}}({\bm{\alpha}},{\bm{\beta}},{T^{[0]}}) + {\rm O}({\xi ^3}),\\
b _{ij}({\bm{\beta}},{T^{\xi}}) &=b _{ij}^{[0]}({\bm{\beta}},{T^{[0]}}) + \xi b _{ij}^{[1]}({\bm{\alpha}},{\bm{\beta}},{T^{[0]}}) + {\xi ^2}b _{ij}^{[2]}({\bm{\alpha}},{\bm{\beta}},{T^{[0]}}) + {\rm O}({\xi ^3}),\\
{c}({\bm{\beta}},{T^{\xi}})& = {c^{[0]}}({\bm{\beta}},{T^{[0]}}) + \xi {c^{[1]}}({\bm{\alpha}},{\bm{\beta}},{T^{[0]}}) + {\xi ^2}{c^{[2]}}({\bm{\alpha}},{\bm{\beta}},{T^{[0]}}) + {\rm O}({\xi ^3}),\\
k_{ij}({\bm{\beta}},{T^{\xi}})& = k_{ij}^{[0]}({\bm{\beta}},{T^{[0]}}) + \xi k_{ij}^{[1]}({\bm{\alpha}},{\bm{\beta}},{T^{[0]}}) + {\xi ^2}k_{ij}^{[2]}({\bm{\alpha}},{\bm{\beta}},{T^{[0]}}) + {\rm O}({\xi ^3}),\\
\vartheta_{ij}({\bm{\beta}},{T^{\xi}})&= \vartheta_{ij}^{[0]}({\bm{\beta}},{T^{[0]}}) + \xi \vartheta_{ij}^{[1]}({\bm{\alpha}},{\bm{\beta}},{T^{[0]}}) + {\xi ^2}\vartheta_{ij}^{[2]}({\bm{\alpha}},{\bm{\beta}},{T^{[0]}}) + {\rm O}({\xi ^3}).
\end{aligned}
\end{equation}

Taking into account the relations between macroscopic and microscopic coordinates, we can define the chain rule of multi-scale asymptotic analysis in orthogonal coordinates as follows
\begin{equation}
\frac{\partial \aleph^\xi(\bm{\alpha },t)}{{\partial {\alpha _i}}} = \frac{\partial \aleph(\bm{\alpha },\bm{\beta },t)}{{\partial {\alpha _i}}} + \frac{1}{\xi }\frac{\partial \aleph(\bm{\alpha },\bm{\beta },t)}{{\partial {\beta _i}}},\;\;(i=1,2,3).
\end{equation}
Then substituting (2.1)-(2.4) into two-scale nonlinear initial-boundary problem (1.1) and expanding its derivatives by chain rule (2.5), we hence have a series of equations by collecting the power-like terms of small periodic parameter $\xi$ as follows
\begin{equation}
O({\xi ^{ - 2}}):\left\{ {\begin{aligned}
&{ - \sum\limits_{j = 1}^3 {\frac{1}{H}\frac{\partial }{{\partial {\beta _j}}}\left( {\frac{H}{{{H_j}}}\sigma _{ij}^{[ - 1]}} \right)}  = 0,}\\
&{\sum\limits_{i = 1}^3 {\frac{1}{H}\frac{\partial }{{\partial {\beta _i}}}\left( {\frac{H}{{{H_i}}}q_i^{[ - 1]}} \right)}  = 0.}
\end{aligned}} \right.
\end{equation}
\begin{equation}
O({\xi ^{ - 1}}):\left\{ {\begin{aligned}
&- \sum\limits_{j = 1}^3 {\frac{1}{H}\frac{\partial }{{\partial {\alpha _j}}}\left( {\frac{H}{{{H_j}}}\sigma _{ij}^{[- 1]}} \right)}  - \sum\limits_{j = 1}^3 {\frac{1}{H}\frac{\partial }{{\partial {\beta _j}}}\left( {\frac{H}{{{H_j}}}\sigma _{ij}^{[0]}} \right)} \\
&- \sum\limits_{j \ne i,j = 1}^3 {\frac{1}{{{H_i}{H_j}}}\frac{{\partial {H_i}}}{{\partial {\alpha _j}}}} \sigma _{ij}^{[ - 1]} + \sum\limits_{j \ne i,j = 1}^3 {\frac{1}{{{H_i}{H_j}}}\frac{{\partial {H_j}}}{{\partial {\alpha _i}}}} \sigma _{jj}^{[ - 1]} = 0,\\
&{\sum\limits_{i = 1}^3 {\frac{1}{H}\frac{\partial }{{\partial {\alpha _i}}}\left( {\frac{H}{{{H_i}}}q_i^{[ - 1]}} \right)} {\rm{ + }}\sum\limits_{i = 1}^3 {\frac{1}{H}\frac{\partial }{{\partial {\beta _i}}}\left( {\frac{H}{{{H_i}}}q_i^{[0]}} \right)}{ + {\vartheta_{ij}^{[0]}}\frac{{\partial \varepsilon _{ij}^{[ - 1]}}}{{\partial t}}} = 0.}
\end{aligned}} \right.
\end{equation}
\begin{equation}
O({\xi ^0}):\left\{ {\begin{aligned}
&\rho^{[0]} \frac{{{\partial ^2}u_i^{[0]}}}{{\partial {t^2}}}-\sum\limits_{j = 1}^3 {\frac{1}{H}\frac{\partial }{{\partial {\alpha _j}}}\left( {\frac{H}{{{H_j}}}\sigma _{ij}^{[0]}} \right)}  - \sum\limits_{j = 1}^3 {\frac{1}{H}\frac{\partial }{{\partial {\beta _j}}}\left( {\frac{H}{{{H_j}}}\sigma _{ij}^{[1]}} \right)} \\
&- \sum\limits_{j \ne i,j = 1}^3 {\frac{1}{{{H_i}{H_j}}}\frac{{\partial {H_i}}}{{\partial {\alpha _j}}}} \sigma _{ij}^{[0]} + \sum\limits_{j \ne i,j = 1}^3 {\frac{1}{{{H_i}{H_j}}}\frac{{\partial {H_j}}}{{\partial {\alpha _i}}}} \sigma _{jj}^{[0]} = {f_i},\\
&{\rho^{[0]} c^{[0]}\frac{{\partial {T^{[0]}}}}{{\partial t}}{\rm{ + }}}\sum\limits_{i = 1}^3 {\frac{1}{H}\frac{\partial }{{\partial {\alpha _i}}}\left( {\frac{H}{{{H_i}}}q_i^{[0]}} \right)} {\rm{ + }}\sum\limits_{i = 1}^3 {\frac{1}{H}\frac{\partial }{{\partial {\beta _i}}}\left( {\frac{H}{{{H_i}}}q_i^{[1]}} \right)} \\
&+ {\vartheta_{ij}^{[0]}}\frac{{\partial \varepsilon _{ij}^{[0]}}}{{\partial t}}+ {\vartheta_{ij}^{[1]}}\frac{{\partial \varepsilon _{ij}^{[-1]}}}{{\partial t}}= h.
\end{aligned}} \right.
\end{equation}
In the subsequent section, each equation with order $O({\xi^{-2}})$, $O({\xi^{- 1}})$ and $O({\xi^{0}})$ is analyzed respectively so that the higher-order multi-scale analysis can be fulfilled.

For simplifying the analysis procedures of this work, new differential operators $\displaystyle\Psi_i=\frac{1}{H_i}\frac{\partial }{{\partial {\alpha _i}}}$ and $\displaystyle\widetilde{\Psi}_i=\frac{1}{H_i}\frac{\partial }{{\partial {\beta _i}}}$ in macro-scale and micro-scale are introduced. Thereupon, the new chain rule for multi-scale asymptotic analysis of heterogeneous shells can be expressed as
\begin{equation}
\Psi_i\big[\aleph^\xi(\bm{\alpha },t)\big]=\Psi_i\big[\aleph(\bm{\alpha },\bm{\beta },t)\big]+\xi^{-1}\widetilde{\Psi}_i\big[\aleph(\bm{\alpha },\bm{\beta },t)\big].
\end{equation}

Next, inspired by the fundamental results in \cite{R27,R46}, we shall obtain the specific expansion form of strain component $\varepsilon_{ij}^\xi(\bm{\alpha },t)$. Each asymptotic terms have the following detailed expressions by substituting (2.2) into (1.4).
\begin{equation}
\begin{aligned}
&\varepsilon _{ij}^{[ - 1]} = \frac{1}{2}\big[ {{{\widetilde \Psi }_i}(u _j^{[0]}) + {{\widetilde \Psi }_j}(u_i^{[0]})} \big],\\
&\varepsilon _{ij}^{[s]} =\breve\varepsilon _{ij}^{[s]} + \frac{1}{2}\big[{{{\widetilde \Psi }_i}(u_j^{[s + 1]}) + {{\widetilde \Psi }_j}(u_i^{[s + 1]})} \big],\\
&\breve\varepsilon _{11}^{[s]} = {\Psi _1}(u_1^{[s]})+\frac{{\Psi _2}(H_1)}{H_1}u_2 ^{[s]}+\frac{{\Psi _3}(H_1)}{H_1}u_3 ^{[s]},\\
&\breve\varepsilon _{22}^{[s]} = {\Psi _2}(u_2^{[s]})+\frac{{\Psi _3}(H_2)}{H_2}u_3 ^{[s]}+\frac{{\Psi _1}(H_2)}{H_2}u_1 ^{[s]},\\
&\breve\varepsilon _{33}^{[s]} = {\Psi _3}(u_3^{[s]})+\frac{{\Psi _1}(H_3)}{H_3}u_1 ^{[s]}+\frac{{\Psi _2}(H_3)}{H_3}u_2 ^{[s]},\\
&\breve\varepsilon _{ij}^{[s]} = \frac{1}{2}\big[ {{\Psi _i }(u_j^{[s]}) + {\Psi _j}(u_i^{[s]})} -\frac{{\Psi _i}(H_j)}{H_j}u_j ^{[s]}-\frac{{\Psi _j}(H_i)}{H_i}u_i^{[s]}\big],\;i\neq j,\;s=0,1.
\end{aligned}
\end{equation}
Combining (2.1)-(2.4) and (2.9)-(2.10) together, the specific expressions of multi-scale expansion terms of $\sigma _{ij}^\xi(\bm{\alpha },t)$ and $q_i^\xi(\bm{\alpha },t)$ can be derived as follows.
\begin{equation}
\begin{aligned}
&\sigma _{ij}^{[- 1]} = {a_{ijkl}^{[0]}}{\widetilde \Psi _k}(u_l^{[0]}),\\
&\sigma _{ij}^{[0]} = {a_{ijkl}^{[0]}}\big[\breve\varepsilon_{kl}^{[0]} + {\widetilde \Psi _k}(u_l^{[1]})\big]+ {a_{ijkl}^{[1]}}{\widetilde \Psi _k}(u_l^{[0]}) - {b _{ij}^{[0]}}( {{T^{[0]}} - \widetilde T)}, \\
&\sigma _{ij}^{[1]} = {a_{ijkl}^{[0]}}\big[\breve\varepsilon_{kl}^{[1]} + {\widetilde \Psi _k}(u_l^{[2]})\big]+ {a_{ijkl}^{[1]}}\big[\breve\varepsilon_{kl}^{[0]}+ {\widetilde \Psi _k}(u_l^{[1]})\big]\\
&\quad\;\;+ {a_{ijkl}^{[2]}}{\widetilde \Psi _k}(u_l^{[0]})- {b _{ij}^{[0]}}{T^{[1]}}- {b _{ij}^{[1]}}{( {{T^{[0]}} - \widetilde T)}}.
\end{aligned}
\end{equation}
\begin{equation}
\begin{aligned}
&q_i^{[ - 1]} =  - {k_{ij}^{[0]}}{\widetilde \Psi _j}{(T^{[0]})},\\
&q_i^{[0]} =  - {k_{ij}^{[0]}}\big[{\Psi _j}{(T^{[0]})}+{\widetilde \Psi _j}{(T^{[1]})}\big]- {k_{ij}^{[1]}}{\widetilde \Psi _j}{(T^{[0]})},\\
&q_i^{[1]} =  - {k_{ij}^{[0]}}\big[{\Psi _j}{(T^{[1]})} +{\widetilde \Psi _j}{(T^{[2]})}\big]- {k_{ij}^{[1]}}\big[{\Psi _j}{(T^{[0]})}+{\widetilde \Psi _j}{(T^{[1]})}\big] - {k_{ij}^{[2]}}{\widetilde \Psi _j}{(T^{[0]})}.
\end{aligned}
\end{equation}

From now on, we initiate to sequentially solve each asymptotic expansion component of $u_i^\xi$ and $T^\xi$ for deriving the detailed expressions of lower-order multi-scale (LOMS) and higher-order multi-scale (HOMS) solutions. Firstly, replacing $\sigma_{ij}^{[-1]}$ and $q_i^{[-1]}$ in $O({\xi^{-2}})$-order equations (2.6) with their detailed expressions in (2.11) and (2.12), classical asymptotic homogenization theory and the orthogonal periodicity of $u_i^{[0]}$ and $T^{[0]}$ imply that $u_i^{[0]}$ and $T^{[0]}$ are irrespective of microscopic variable $\bm{\beta}$ such that
\begin{equation}
\left\{
\begin{aligned}
&u_i^{[0]}({\bm{\alpha }},{\bm{\beta }},t) = u_i^{[0]}({\bm{\alpha }},t),\\
&{T^{[0]}}({\bm{\alpha }},{\bm{\beta }},t) = {T^{[0]}}({\bm{\alpha }},t).
\end{aligned}\right.
\end{equation}
For brevity, we denote $\bm{H}=({H_1},{H_2},{H _3})$. In the following, replacing the associated asymptotic terms $\sigma _{ij}^{[-1]}$, $q_i^{[-1]}$, $\sigma _{ij}^{[0]}$, $q_i^{[0]}$ and $\varepsilon_{ij}^{[-1]}$ in $O({\xi^{-1}})$-order equations (2.7) and employing (2.13), the effective expressions for $u_i^{[1]}$ and $T^{[1]}$ can be established as below
\begin{equation}
\left\{
\begin{aligned}
&u_i^{[1]}({\bm{\alpha }},{\bm{\beta }},t) = N_i^{mn}(\bm{H},T^{[0]},\bm{\beta })\breve\varepsilon_{mn}^{[0]} - {P_i}(\bm{H},T^{[0]},\bm{\beta })({T^{[0]}} - \widetilde T),\\
&{T^{[1]}}({\bm{\alpha }},{\bm{\beta }},t) = {M^m}(\bm{H},T^{[0]},\bm{\beta }){\Psi _m}( {{T^{[0]}}})
,\;\;m,n=1,2,3,
\end{aligned}\right.
\end{equation}
where $N_i^{mn}$, $P_i$ and ${M^m}$ are termed as first-order cell functions defined in PUC $\Theta$, that all rely upon the macroscopic geometric parameter $\bm{H}$ and macroscopic temperature field $T^{[0]}$. These macroscopic parameters play a role of variable parameters. This is an apparent discrepancy compared to classical heterogeneous materials possessing microscopic periodicity in cartesian coordinates. Moreover, these first-order auxiliary cell functions are obtained by solving the following equations equipped with periodic boundary condition
\begin{equation}
\left\{ {\begin{aligned}
&{{\widetilde \Psi }_j}\big[{a_{ijkl}^{[0]}}{{\widetilde \Psi }_k}(N_l^{mn})\big] =  - {{\widetilde \Psi }_j}({a_{ijmn}^{[0]}}),&&\;\;{\rm{  }}{\bm{\beta }} \in \Theta,\\
&N_l^{mn}(\bm{H},T^{[0]},\bm{\beta })\;\mathrm{is}\;1-\mathrm{periodic}\;\mathrm{in}\;{\bm{\beta }}{\rm{, }}&&{\int_{\Theta}}N_l^{mn}d\Theta=0.
\end{aligned}} \right.
\end{equation}
\begin{equation}
\left\{ {\begin{aligned}
&{{\widetilde \Psi }_j}\big[{a_{ijkl}^{[0]}}{{\widetilde \Psi }_k}({P_l})\big] =  - {{\widetilde \Psi }_j}({b _{ij}^{[0]}}),&&{\rm{}}\;\;{\bm{\beta }} \in \Theta,\\
&{P_l}(\bm{H},T^{[0]},\bm{\beta })\;\mathrm{is}\;1-\mathrm{periodic}\;\mathrm{in}\;{\bm{\beta }}{\rm{, }}&&{\int_{\Theta}}{P_l}d\Theta=0.
\end{aligned}} \right.
\end{equation}
\begin{equation}
\left\{ \begin{aligned}
&{{\widetilde \Psi }_i}\big[{k_{ij}^{[0]}}{{\widetilde \Psi }_j}({M^m})\big] =  - {{\widetilde \Psi }_i}({k_{im}^{[0]}}),&&\;\;{\bm{\beta }} \in \Theta,\\
&{M^m}(\bm{H},T^{[0]},\bm{\beta })\;\mathrm{is}\;1-\mathrm{periodic}\;\mathrm{in}\;{\bm{\beta }}{\rm{, }}&&{\int_{\Theta}}{M^m}d\Theta=0.
\end{aligned} \right.
\end{equation}
Afterwards, applying the homogenized operator to both sides of $O({\xi^{0}})$-order equations (2.8) with respect to $\bm{\beta}$ over PUC $\Theta$ and employing the Green's formula on (2.8), we shall derive the macroscopic homogenized equations related with (1.1) as below
\begin{equation}
\left\{ {\begin{aligned}
&\widehat \rho \frac{{{\partial ^2}u_i^{[0]}}}{{\partial {t^2}}}- \sum\limits_{j = 1}^3 {\frac{H_j}{H}{\Psi _j}\left[ {\frac{H}{{{H_j}}}\left( {{{\widehat a}_{ijkl}}\breve\varepsilon _{kl}^{[0]} - {{\widehat b }_{ij}}({T^{[0]}} - \widetilde T)} \right)} \right]}\\
&\quad\quad-\sum\limits_{j \ne i,j = 1}^3 {\frac{{\Psi _j}({H_i})}{{{H_i}}}} \left[ {{{\widehat a}_{ijkl}}\breve\varepsilon _{kl}^{[0]} - {{\widehat b }_{ij}}({T^{[0]}} - \widetilde T)} \right]\\
&\quad\quad+\sum\limits_{j \ne i,j = 1}^3 {\frac{{\Psi _i}({H_j})}{{{H_j}}}} \left[ {{{\widehat a}_{jjkl}}\breve\varepsilon _{kl}^{[0]} - {{\widehat b }_{jj}}({T^{[0]}} - \widetilde T)} \right] = {f_i},\;\;\text{in}\;\;\Omega\times(0,\mathcal T),\\
&{\widehat S\frac{{\partial {T^{[0]}}}}{{\partial t}}}{ - \sum\limits_{i = 1}^3 {\frac{H_i}{H}{\Psi _i}\left[ {\frac{H}{{{H_i}}}{{\widehat k}_{ij}}{\Psi _j}({T^{[0]}})} \right]}+ \widehat\vartheta_{ij}\frac{{\partial \breve\varepsilon _{ij}^{[0]}}}{{\partial t}} = h,\;\;\text{in}\;\;\Omega\times(0,\mathcal T),}\\
&\bm{u}^{[0]}({\bm{\alpha }},t)=\widehat{\bm{u}}({\bm{\alpha }},t),\;\;\text{on}\;\;\partial\Omega_{u}\times(0,\mathcal T),\\
&T^{[0]}({\bm{\alpha }},t)=\widehat{T}({\bm{\alpha }},t),\;\;\text{on}\;\;\partial\Omega_{T}\times(0,\mathcal T),\\
&\bm{u}^{[0]} (\bm{\alpha },0)=\widetilde{\bm{u}}^0,\;\;\frac{\partial \bm{u}^{[0]} (\bm{\alpha },t)}{\partial t}\big|_{t=0}=\widetilde{\bm{u}}^1 (\bm{\alpha }),\;\;T^{[0]}(\bm{\alpha },0)=\widetilde T,\;\;\text{in}\;\;\Omega.
\end{aligned}}\right.
\end{equation}
In aforementioned macroscopic homogenized equations (2.18), the homogenized material parameters without multi-scale property can be defined in macro-scale as below
\begin{equation}
\begin{aligned}
&\widehat \rho(T^{[0]})= \left\langle \rho^{[0]}  \right\rangle,\\
&{\widehat a_{ijkl}}(\bm{H},T^{[0]})= \left\langle {{a_{ijkl}^{[0]}} + {a_{ijmn}^{[0]}}{{\widetilde \Psi }_m}(N_n^{kl})} \right\rangle,\\
&{\widehat b _{ij}}(\bm{H},T^{[0]})= \left\langle {{b _{ij}^{[0]}} + {a_{ijkl}^{[0]}}{{\widetilde \Psi }_k}({P_l})} \right\rangle,\\
&\widehat S(\bm{H},T^{[0]}) = \left\langle {\rho^{[0]} c^{[0]} - {\vartheta_{ij}^{[0]}}{{\widetilde \Psi }_j}({P_i})} \right\rangle,\\
&{\widehat k_{ij}}(\bm{H},T^{[0]})= \big\langle {k_{ij}^{[0]} + k_{im}^{[0]}{{\widetilde \Psi }_m}{(M^j)}} \big\rangle,\\
&\widehat\vartheta_{ij}(\bm{H},T^{[0]}) = \left\langle {{\vartheta_{ij}^{[0]}} + {\vartheta_{mn}^{[0]}}{{\widetilde \Psi }_n}(N_m^{ij})} \right\rangle,
\end{aligned}
\end{equation}
where we introduced a bracket operator $\left\langle\bullet\right\rangle$ as $\left\langle\bullet\right\rangle=\displaystyle\frac{1}{|\Theta|}{\int_{\Theta}}\bullet d\Theta$. It should also
be noted that the homogenized material parameters excluding mass density $\widehat \rho$ alter with macroscopic geometric parameter $\bm{H}$ and macroscopic temperature field $T^{[0]}$. In essence, heterogeneous shells are physically homogenized as functionally graded shells in macro-scale.

To seek the essential second-order cell functions, we first replace the terms $\sigma _{ij}^{[0]}$, $q_i^{[0]}$, $\sigma _{ij}^{[1]}$, $q_i^{[1]}$, $\varepsilon_{ij}^{[-1]}$ and $\varepsilon_{ij}^{[0]}$ involving in (2.8) with  their detailed expressions in (2.10)-(2.12). Then, eliminating the terms $f_i$ and $h$ involving in (2.8) from the homogenized equations (2.18), we establish the expansion expressions for $
u_i^{[2]}({\bm{\alpha }},{\bm{\beta }},t)$ and ${T^{[2]}}({\bm{\alpha }},{\bm{\beta }},t)$ as below
\begin{equation}
\left\{
\begin{aligned}
&u_i^{[2]}({\bm{\alpha }},{\bm{\beta }},t) = N_i^{pmn}(\bm{H},T^{[0]},\bm{\beta }){\Psi _p}(\breve\varepsilon_{mn}^{[0]}) + H_i^p(\bm{H},T^{[0]},\bm{\beta }){\Psi _p}({T^{[0]}})\\
&\quad\quad\quad\quad\quad+ F_i^p(\bm{H},T^{[0]},\bm{\beta })\frac{{{\partial ^2}u_p^{[0]}}}{{\partial {t^2}}} + W_i^{mn}(\bm{H},T^{[0]},\bm{\beta })\breve\varepsilon_{mn}^{[0]} \\
&\quad\quad\quad\quad\quad+ {Q_i}(\bm{H},T^{[0]},\bm{\beta })({T^{[0]}} - \widetilde T)\\
&\quad\quad\quad\quad\quad + Z_i^p(\bm{H},T^{[0]},\bm{\beta }){\Psi _p}({T^{[0]}})({T^{[0]}} - \widetilde T)\\
&\quad\quad\quad\quad\quad -X_i^{pmn}(\bm{H},T^{[0]},\bm{\beta }){\Psi }_p( T^{[0]})\breve\varepsilon _{mn}^{[0]},\\
&{T^{[2]}}({\bm{\alpha }},{\bm{\beta }},t) ={M^{mn}}(\bm{H},T^{[0]},\bm{\beta }){\Psi _m}{\Psi _n}({T^{[0]}})+ A(\bm{H},T^{[0]},\bm{\beta })\frac{{\partial {T^{[0]}}}}{{\partial t}}\\
&\quad\quad\quad\quad\quad\;+ {G^{mn}}(\bm{H},T^{[0]},\bm{\beta })\frac{{\partial \breve\varepsilon_{mn}^{[0]}}}{{\partial t}} + {R^m}(\bm{H},T^{[0]},\bm{\beta }){\Psi _m}({T^{[0]}}) \\
&\quad\quad\quad\quad\quad\;-B^{mn}(\bm{H},T^{[0]},\bm{\beta }){{\Psi }_m}( T^{[0]}){{\Psi }_n} (T^{[0]}),\;\;p=1,2,3,
\end{aligned}\right.
\end{equation}
In (2.20), $N_i^{pmn}$, $H_i^p$, $F_i^p$, $W_i^{mn}$, $Q_i$, $Z_i^p$, $X_i^{pmn}$, $M^{mn}$, $A$, $G^{mn}$, $R^m$ and $B^{mn}$ are the second-order cell functions defined in PUC $\Theta$, which can be obtained by solving an array of equations equipping with periodic boundary condition as below
\begin{equation}
\left\{ \begin{aligned}
&{{\widetilde \Psi }_j}\big[{a_{ijkl}^{[0]}}{{\widetilde \Psi }_k}(N_l^{pmn})\big] = {{\widehat a}_{ipmn}} - {a_{ipmn}^{[0]}} - {a_{ipkl}^{[0]}}{{\widetilde \Psi }_k}(N_l^{mn})\\
&\quad\quad\quad\quad\quad\quad\quad\quad\;\; - {{\widetilde \Psi }_l}({a_{ilkp}^{[0]}}N_k^{mn}),\;\;\;\;{\bm{\beta }} \in \Theta,\\
&N_l^{pmn}(\bm{H},T^{[0]},\bm{\beta })\;\mathrm{is}\;1-\mathrm{periodic}\;\mathrm{in}\;{\bm{\beta }}{\rm{, }}\;\;\;\;{\int_{\Theta}}N_l^{pmn}d\Theta=0.
\end{aligned} \right.
\end{equation}
\begin{equation}
\left\{ \begin{aligned}
&{{\widetilde \Psi }_j}\big[{a_{ijkl}^{[0]}}{{\widetilde \Psi }_k}(H_l^p)\big] = {b_{ip}^{[0]}} - {{\widehat b}_{ip}} + {a_{ipkl}^{[0]}}{{\widetilde \Psi }_k}({P_l}) + {{\widetilde \Psi }_l}({a_{ilkp}^{[0]}}{P_k})\\
&\quad\quad\quad\quad\quad\quad\quad\;\; + {{\widetilde \Psi }_k}({b _{ik}^{[0]}}{M^p}),\;\;\;\;{\bm{\beta }} \in \Theta,\\
&H_l^p(\bm{H},T^{[0]},\bm{\beta })\;\mathrm{is}\;1-\mathrm{periodic}\;\mathrm{in}\;{\bm{\beta }}{\rm{, }}\;\;\;\;{\int_{\Theta}}H_l^pd\Theta=0.
\end{aligned} \right.
\end{equation}
\begin{equation}
\left\{ \begin{aligned}
&{{{\widetilde \Psi }_j}[{a_{ijkl}^{[0]}}{{\widetilde \Psi }_k}(F_l^p)] = {\delta _{ip}}\big( {\rho^{[0]}  - \widehat \rho } \big),\;\;\;\;{\rm{ }}{\bm{\beta }} \in \Theta,}\\
&{F_l^p}(\bm{H},T^{[0]},\bm{\beta })\;\mathrm{is}\;1-\mathrm{periodic}\;\mathrm{in}\;{\bm{\beta }}{\rm{, }}\;\;\;\;{\int_{\Theta}}{F_l^p}d\Theta=0.
\end{aligned} \right.
\end{equation}
\begin{equation}
\left\{ \begin{aligned}
&{{\widetilde \Psi }_j}\big[{a_{ijkl}^{[0]}}{{\widetilde \Psi }_k}(W_l^{mn})\big] =
- \sum\limits_{j = 1}^3 {\frac{{{H_j}}}{H}{\Psi _j}\Big[ {\frac{H}{{{H_j}}}\big( {{a_{ijmn}^{[0]}} - {{\widehat a}_{ijmn}} + {a_{ijkl}^{[0]}}{{\widetilde \Psi }_k}(N_l^{mn})} \big)} \Big]} \\
&\quad\;\;\;- \sum\limits_{j \ne i,j = 1}^3 {\frac{{\Psi _j}({H_i} )}{{{H_i}}}}\big( {{a_{ijmn}^{[0]}} - {{\widehat a}_{ijmn}} + {a_{ijkl}^{[0]}}{{\widetilde \Psi }_k}(N_l^{mn})} \big)- {\widetilde \Psi _j}({a_{ijkl}^{[0]}}{D_{klmn}})\\
&\quad\;\;\;+\sum\limits_{j \ne i,j = 1}^3 {\frac{{\Psi _i}({H_j} )}{{{H_j}}}}\big( {{a_{jjmn}^{[0]}} - {{\widehat a}_{jjmn}} + {a_{jjkl}^{[0]}}{{\widetilde \Psi }_k}(N_l^{mn})} \big),\;\;\;\;\bm{\beta } \in \Theta,\\
&W_l^{mn}(\bm{H},T^{[0]},\bm{\beta })\;\mathrm{is}\;1-\mathrm{periodic}\;\mathrm{in}\;{\bm{\beta }}{\rm{, }}\;\;\;\;{\int_{\Theta}}W_l^{mn}d\Theta=0.
\end{aligned}\right.
\end{equation}
\begin{equation}
\left\{ \begin{aligned}
&{{\widetilde \Psi }_j}\big[{a_{ijkl}^{[0]}}{{\widetilde \Psi }_k}({Q_l})\big] =
\sum\limits_{j = 1}^3 {\frac{{{H_j}}}{H}{\Psi _j}\Big[ {\frac{H}{{{H_j}}}\big( {{b_{ij}^{[0]}} - {{\widehat b}_{ij}} + {a_{ijkl}^{[0]}}{{\widetilde \Psi }_k}({P_l})} \big)} \Big]} \\
&\quad\quad\quad\;+ \sum\limits_{j \ne i,j = 1}^3 {\frac{{\Psi _j}({H_i})}{{{H_i}}}} \big( {{b_{ij}^{[0]}} - {{\widehat b }_{ij}} + {a_{ijkl}^{[0]}}{{\widetilde \Psi }_k}({P_l})} \big)+ {\widetilde \Psi _j}({a_{ijkl}^{[0]}}{E_{kl}})\\
&\quad\quad\quad\;- \sum\limits_{j \ne i,j = 1}^3 {\frac{{\Psi _i}({H_j})}{{{H_j}}}} \big( {{b _{jj}^{[0]}} - {{\widehat b }_{jj}} + {a_{jjkl}^{[0]}}{{\widetilde \Psi }_k}({P_l})} \big)
,\;\;\;\;{\bm{\beta }} \in \Theta,\\
&{Q_l}(\bm{H},T^{[0]},\bm{\beta })\;\mathrm{is}\;1-\mathrm{periodic}\;\mathrm{in}\;{\bm{\beta }}{\rm{, }}\;\;\;\;{\int_{\Theta}}{Q_l}d\Theta=0.
\end{aligned} \right.
\end{equation}
\begin{equation}
\left\{ \begin{aligned}
&{{{\widetilde \Psi }_j}[{a_{ijkl}^{[0]}}{{\widetilde \Psi }_k}(Z_l^p)] = {\widetilde \Psi }_j\big[{{M^{p}}\mathbf{D}^{(0,1)}b _{ij}^{[0]} + {M^{p}}\mathbf{D}^{(0,1)}a_{i jkl}^{[0]}{\widetilde \Psi }_l( P_k)} \big],\;\;\;\;{\rm{ }}{\bm{\beta }} \in \Theta,}\\
&{Z_l^p}(\bm{H},T^{[0]},\bm{\beta })\;\mathrm{is}\;1-\mathrm{periodic}\;\mathrm{in}\;{\bm{\beta }}{\rm{, }}\;\;\;\;{\int_{\Theta}}{Z_l^p}d\Theta=0.
\end{aligned} \right.
\end{equation}
\begin{equation}
\left\{ \begin{aligned}
&{{\widetilde \Psi }_j}\big[{a_{ijkl}^{[0]}}{{\widetilde \Psi }_k}(X_l^{pmn})\big] = {\widetilde \Psi }_j\big[ {{M^{p}}\mathbf{D}^{(0,1)}a_{ijmn}^{[0]} + {M^{p}}\mathbf{D}^{(0,1)}a_{ijkl}^{[0]}{\widetilde \Psi }_l( N^{mn}_{k}}) \big],\;\;\;\;{\bm{\beta }} \in \Theta,\\
&X_l^{pmn}(\bm{H},T^{[0]},\bm{\beta })\;\mathrm{is}\;1-\mathrm{periodic}\;\mathrm{in}\;{\bm{\beta }}{\rm{, }}\;\;\;\;{\int_{\Theta}}X_l^{pmn}d\Theta=0.
\end{aligned} \right.
\end{equation}
\begin{equation}
\left\{ \begin{aligned}
&{{\widetilde \Psi }_i}\big[{k_{ij}^{[0]}}{{\widetilde \Psi }_j}({M^{mn}})\big] = {{\widehat k}_{mn}} - {k_{mn}^{[0]}} - {k_{mj}^{[0]}}{{\widetilde \Psi }_j}({M^n}) - {{\widetilde \Psi }_j}({k_{mj}^{[0]}}{M^n}),\;\;\;\;{\bm{\beta }} \in \Theta,\\
&{M^{mn}}(\bm{H},T^{[0]},\bm{\beta })\;\mathrm{is}\;1-\mathrm{periodic}\;\mathrm{in}\;{\bm{\beta }}{\rm{, }}\;\;\;\;{\int_{\Theta}}{M^{mn}}d\Theta=0.
\end{aligned} \right.
\end{equation}
\begin{equation}
\left\{ \begin{aligned}
&{{{\widetilde \Psi }_i}[{k_{ij}^{[0]}}{{\widetilde \Psi }_j}(A)] = \rho ^{[0]}c ^{[0]}- \widehat S - {\vartheta_{ij}^{[0]}}{{\widetilde \Psi }_i}({P_j}),\;\;\;\;{\rm{}}{\bm{\beta }} \in \Theta,}\\
&{A}(\bm{H},T^{[0]},\bm{\beta })\;\mathrm{is}\;1-\mathrm{periodic}\;\mathrm{in}\;{\bm{\beta }}{\rm{, }}\;\;\;\;{\int_{\Theta}}{A}d\Theta=0.
\end{aligned} \right.
\end{equation}
\begin{equation}
\left\{ \begin{aligned}
&{{{\widetilde \Psi }_i}[{k_{ij}^{[0]}}{{\widetilde \Psi }_j}({G^{mn}})] = {\vartheta_{mn}^{[0]}} - {{\widehat \vartheta}_{mn}} + {\vartheta _{ij}^{[0]}}{{\widetilde \Psi }_i}(N_j^{mn}),\;\;\;\;{\rm{}}{\bm{\beta }} \in \Theta,}\\
&{G^{mn}}(\bm{H},T^{[0]},\bm{\beta })\;\mathrm{is}\;1-\mathrm{periodic}\;\mathrm{in}\;{\bm{\beta }}{\rm{, }}\;\;\;\;{\int_{\Theta}}{G^{mn}}d\Theta=0.
\end{aligned} \right.
\end{equation}
\begin{equation}
\left\{ \begin{aligned}
&{{\widetilde \Psi }_i}\big[{k_{ij}^{[0]}}{{\widetilde \Psi }_j}({R^m})\big] =  {\sum\limits_{i = 1}^3 {\frac{{{H_i}}}{H}{\Psi _i}\Big[ {\frac{H}{{{H_i}}}\big( {{{\widehat k}_{im}} - {k_{im}^{[0]}} - {k_{ij}^{[0]}}{{\widetilde \Psi }_j}({M^m})} \big)} \Big]}}\\
&\quad\quad\quad\quad\quad\quad\quad\;{{- {{\widetilde \Psi }_i}\big[ {{k_{ij}^{[0]}}{\Psi _j}({M^m})} \big]} } ,\;\;\;\;{\bm{\beta }} \in \Theta,\\
&{R^m}(\bm{H},T^{[0]},\bm{\beta })\;\mathrm{is}\;1-\mathrm{periodic}\;\mathrm{in}\;{\bm{\beta }}{\rm{, }}\;\;\;\;{\int_{\Theta}}{R^m}d\Theta=0.
\end{aligned}\right.
\end{equation}
\begin{equation}
\left\{ \begin{aligned}
&{{\widetilde \Psi }_i}\big[{k_{ij}^{[0]}}{{\widetilde \Psi }_j}({B^{mn}})\big] = {{\widetilde \Psi }_i}\big[{{M^{m}}\mathbf{D}^{(0,1)}{k_{in}^{[0]}}}+ {{M^{m}}\mathbf{D}^{(0,1)}{k_{ij}^{[0]}}{{\widetilde \Psi }_j}( M^{n})}\big],\;\;\;\;{\bm{\beta }} \in \Theta,\\
&{B^{mn}}(\bm{H},T^{[0]},\bm{\beta })\;\mathrm{is}\;1-\mathrm{periodic}\;\mathrm{in}\;{\bm{\beta }}{\rm{, }}\;\;\;\;{\int_{\Theta}}{B^{mn}}d\Theta=0.
\end{aligned} \right.
\end{equation}
\begin{rmk}
By the hypotheses (A)-(B) and Lax-Milgram theorem, we can demonstrate that equations (2.15)-(2.17) and (2.21)-(2.32) all possess a unique solution for any fixed macroscopic geometric parameter $\bm{H}$ and macroscopic temperature field $T^{[0]}$. In addition, the detailed expressions of undefined terms $D_{ijmn}$ in (2.24) and $E_{ij}$ in (2.25) can be found in reference \cite{R27}.
\end{rmk}

Summing up, two types of multi-scale asymptotic solutions are established for time-dependent nonlinear thermo-mechanical equations (1.1) of heterogeneous shells. The LOMS solutions for displacement $u_i^\xi$ and temperature $T^\xi$ are given by
\begin{equation}
{\begin{aligned}
u^{[1\xi]}_{i}(\bm{\alpha },t)&=u^{[0]}_{i}(\bm{\alpha },t)+\xi u^{[1]}_{i}({\bm{\alpha }},{\bm{\beta }},t)\\
&= u_i^{[0]}+\xi\big[N_i^{mn}(\bm{H},T^{[0]},{\bm{\beta }})\breve\varepsilon_{mn}^{[0]} - {P_i}(\bm{H},T^{[0]},{\bm{\beta }})( {{T^{[0]}} - \widetilde T} )\big].
\end{aligned}}
\end{equation}
\begin{equation}
{\begin{aligned}
T^{[1\xi]}(\bm{\alpha },t)&=T^{[0]}(\bm{\alpha },t)+\xi T^{[1]}({\bm{\alpha }},{\bm{\beta }},t)\\
&=T^{[0]}+\xi{M^m}(\bm{H},T^{[0]},{\bm{\beta }}){\Psi _m}( {{T^{[0]}}}).
\end{aligned}}
\end{equation}
Moreover, the HOMS solutions for displacement $u_i^\xi$ and temperature $T^\xi$ are given by
\begin{equation}
{\begin{aligned}
u^{[2\xi]}_{i}(\bm{\alpha },t)&=u^{[0]}_{i}(\bm{\alpha },t)+\xi u^{[1]}_{i}({\bm{\alpha }},{\bm{\beta }},t)+\xi^2 u^{[2]}_{i}({\bm{\alpha }},{\bm{\beta }},t)\\
&= u_i^{[0]}+\xi\big[N_i^{mn}(\bm{H},T^{[0]},{\bm{\beta }})\breve\varepsilon_{mn}^{[0]} - {P_i}(\bm{H},T^{[0]},{\bm{\beta }})( {{T^{[0]}} - \widetilde T} )\big]\\
&+ \xi^2\big[ N_i^{pmn}(\bm{H},T^{[0]},\bm{\beta }){\Psi _p}(\breve\varepsilon_{mn}^{[0]}) + H_i^p(\bm{H},T^{[0]},\bm{\beta }){\Psi _p}({T^{[0]}})\\
&+ F_i^p(\bm{H},T^{[0]},\bm{\beta })\frac{{{\partial ^2}u_p^{[0]}}}{{\partial {t^2}}} + W_i^{mn}(\bm{H},T^{[0]},\bm{\beta })\breve\varepsilon_{mn}^{[0]} \\
&+ {Q_i}(\bm{H},T^{[0]},\bm{\beta })({T^{[0]}} - \widetilde T) + Z_i^p(\bm{H},T^{[0]},\bm{\beta }){\Psi _p}({T^{[0]}})({T^{[0]}} - \widetilde T)\\
&-X_i^{pmn}(\bm{H},T^{[0]},\bm{\beta }){\Psi }_p( T^{[0]})\breve\varepsilon _{mn}^{[0]}\big].
\end{aligned}}
\end{equation}
\begin{equation}
{\begin{aligned}
T^{[2\xi]}(\bm{\alpha },t)&=T^{[0]}(\bm{\alpha },t)+\xi T^{[1]}({\bm{\alpha }},{\bm{\beta }},t)+\xi^2 T^{[2]}({\bm{\alpha }},{\bm{\beta }},t)\\
&=T^{[0]}+\xi{M^m}(\bm{H},T^{[0]},{\bm{\beta }}){\Psi _m}( {{T^{[0]}}})\\
&+ \xi^2\big[{M^{mn}}(\bm{H},T^{[0]},\bm{\beta }){\Psi _m}{\Psi _n}({T^{[0]}})+ A(\bm{H},T^{[0]},\bm{\beta })\frac{{\partial {T^{[0]}}}}{{\partial t}}\\
&+ {G^{mn}}(\bm{H},T^{[0]},\bm{\beta })\frac{{\partial \breve\varepsilon_{mn}^{[0]}}}{{\partial t}} + {R^m}(\bm{H},T^{[0]},\bm{\beta }){\Psi _m}({T^{[0]}})\\
&-B^{mn}(\bm{H},T^{[0]},\bm{\beta }){{\Psi }_m}( T^{[0]}){{\Psi }_n} (T^{[0]})\big].
\end{aligned}}
\end{equation}
\noindent As appeared in (2.35) and (2.36), we can conclude that only the HOMS solutions can characterize the impact of displacement fields on temperature field owing to the existence of corrected term $\displaystyle {G^{mn}}(\bm{H},T^{[0]},\bm{\beta })\frac{{\partial \breve\varepsilon_{mn}^{[0]}}}{{\partial t}}$ in (2.36). In addition, by substituting (2.35) and (2.36) into (1.2)-(1.4), we can gain the high-precision HOMS solutions for the strains $\varepsilon _{ij}^\xi$, stresses $\sigma _{ij}^{\xi}$ and heat fluxes $q_i^\xi$ of multi-scale nonlinear problem (1.1). Furthermore, on the basis of (2.35) and (2.36), we can derive the HOMS asymptotic solutions for different heterogeneous structures with cartesian, cylindrical and orthogonal periodic configurations.
\section{The error analyses of multi-scale approximate solutions}
\subsection{The conservation proof of local energy and momentum by local error analysis}
First, we further define residual functions for LOMS and HOMS approximate solutions as below
\begin{align}
{\begin{array}{*{20}{l}}
u^{[1\xi[}_{\Delta i}=u^\xi_{i}-u^{[1\xi]}_{i},\;\;T^{[1\xi]}_{\Delta}=T^\xi-T^{[1\xi]},\\
u^{[2\xi]}_{\Delta i}=u^\xi_{ i}-u^{[2\xi]}_{i},\;\;T^{[2\xi]}_{\Delta}=T^\xi-T^{[2\xi]}.
\end{array}}
\end{align}

To analyze the local energy and momentum conservation of multi-scale solutions, substituting the above residual functions $u^{[1\xi]}_{\Delta i}$ and $T^{[1\xi]}_{\Delta}$ into (1.1), the residual equations for the LOMS solutions can be derived as below.
\begin{equation}
\left\{
\begin{aligned}
&{\rho ^\xi }\frac{{\partial ^2}u_{\Delta i}^{[1\xi]}}{{\partial t^2}}-\sum\limits_{j = 1}^3 {\frac{1}{H}\frac{\partial }{{\partial {\alpha _j}}}\left( {\frac{H}{{{H_j}}}\sigma _{ij}^\xi }(\bm{u}_{\Delta }^{[1\xi]},T_{\Delta}^{[1\xi]}) \right)}\\
&\quad\quad\quad\quad\; -\sum\limits_{j \ne i,j = 1}^3 {\frac{1}{{{H_i}{H_j}}}\frac{{\partial {H_i}}}{{\partial {\alpha _j}}}} \sigma _{ij}^\xi(\bm{u}_{\Delta }^{[1\xi]},T_{\Delta}^{(1\xi)})\\
&\quad\quad\quad\quad\; + \sum\limits_{j \ne i,j = 1}^3 {\frac{1}{{{H_i}{H_j}}}\frac{{\partial {H_j}}}{{\partial {\alpha _i}}}} \sigma _{jj}^\xi (\bm{u}_{\Delta }^{[1\xi]},T_{\Delta}^{(1\xi)})\\
&\quad\quad\quad\quad\;  = \mathbb{F}_{0i}({\bm{\alpha }},{\bm{\beta }},t)+\xi \mathbb{F}_{1i}({\bm{\alpha }},{\bm{\beta }},t),\;\;\text{in}\;\;\Omega\times(0,\mathcal T),\\
&{{\rho ^\xi }{c^\xi }\frac{{\partial {T_{\Delta}^{[1\xi]}}}}{{\partial t}}{\rm{ + }}}\sum\limits_{i = 1}^3 {\frac{1}{H}\frac{\partial }{{\partial {\alpha _i}}}\left( {\frac{H}{{{H_i}}}q_i^\xi(T_{\Delta}^{[1\xi]}) } \right)} { + \vartheta_{ij}^\xi \frac{{\partial \varepsilon_{ij}^\xi (\bm{u}_{\Delta }^{[1\xi]})}}{{\partial t}}}\\
&\quad\quad\quad\quad\;\; = \mathbb{S}_{0}({\bm{\alpha }},{\bm{\beta }},t)+\xi \mathbb{S}_{1}({\bm{\alpha }},{\bm{\beta }},t),\;\;\text{in}\;\;\Omega\times(0,\mathcal T),
\end{aligned} \right.
\end{equation}
where the operators $\sigma _{ij}^{\xi}(\bm{u}_{\Delta }^{[1\xi]},T_{\Delta}^{[1\xi]})$ and $q_i^{\xi}(T_{\Delta}^{[1\xi]})$ are defined as $\sigma _{ij}^{\xi}(\bm{u}_{\Delta }^{[1\xi]},T_{\Delta}^{[1\xi]})={a_{ijkl}^\xi}\varepsilon _{kl}^{\xi}(\bm{u}_{\Delta }^{[1\xi]}) - {b_{ij}^\xi}T_{\Delta}^{[1\xi]}$ and $q_i^{\xi}(T_{\Delta}^{[1\xi]})= -{k_{ij}^\xi}{\Psi _j}(T_{\Delta}^{[1\xi]})$ in (3.2). The specific expressions of functions $\mathbb{F}_{0i}({\bm{\alpha }},{\bm{\beta }},t)$, $\mathbb{F}_{1i}({\bm{\alpha }},{\bm{\beta }},t)$, $\mathbb{S}_{0}({\bm{\alpha }},{\bm{\beta }},t)$ and $\mathbb{S}_{1}({\bm{\alpha }},{\bm{\beta }},t)$ are uncomplicated to obtain and not exhibited in the present study because of their lengthy forms.

Then putting the residual functions $u^{[2\xi]}_{\Delta i}$ and $T^{[2\xi]}_{\Delta}$ into (1.1), we have the following residual equations for the HOMS solutions.
\begin{equation}
\left\{
\begin{aligned}
&{\rho ^\xi }\frac{{\partial ^2}u_{\Delta i}^{[2\xi]}}{{\partial t^2}}-\sum\limits_{j = 1}^3 {\frac{1}{H}\frac{\partial }{{\partial {\alpha _j}}}\left( {\frac{H}{{{H_j}}}\sigma _{ij}^\xi (\bm{u}_{\Delta }^{[2\xi]},T_{\Delta}^{[2\xi]})} \right)}\\
&\quad\quad\quad\quad\;-\sum\limits_{j \ne i,j = 1}^3 {\frac{1}{{{H_i}{H_j}}}\frac{{\partial {H_i}}}{{\partial {\alpha _j}}}} \sigma _{ij}^\xi(\bm{u}_{\Delta }^{[2\xi]},T_{\Delta}^{(2\xi)})\\
&\quad\quad\quad\quad\;+ \sum\limits_{j \ne i,j = 1}^3 {\frac{1}{{{H_i}{H_j}}}\frac{{\partial {H_j}}}{{\partial {\alpha _i}}}} \sigma _{jj}^\xi(\bm{u}_{\Delta }^{[2\xi]},T_{\Delta}^{[2\xi]}) \\
&\quad\quad\quad\quad\; = \xi \mathbb{G}_i({\bm{\alpha }},{\bm{\beta }},t),\;\;\text{in}\;\;\Omega\times(0,\mathcal T),\\
&{{\rho ^\xi }{c^\xi }\frac{{\partial {T_{\Delta}^{[2\xi]}}}}{{\partial t}}{\rm{ + }}}\sum\limits_{i = 1}^3 {\frac{1}{H}\frac{\partial }{{\partial {\alpha _i}}}\left( {\frac{H}{{{H_i}}}q_i^\xi(T_{\Delta}^{[2\xi]}) } \right)}{ + \vartheta_{ij}^\xi \frac{{\partial \varepsilon_{ij}^\xi (\bm{u}_{\Delta }^{[2\xi]})}}{{\partial t}}}\\
&\quad\quad\quad\quad\quad= \xi \mathbb{E}({\bm{\alpha }},{\bm{\beta }},t),\;\;\text{in}\;\;\Omega\times(0,\mathcal T),
\end{aligned} \right.
\end{equation}
where the operators $\sigma _{ij}^{\xi}(\bm{u}_{\Delta }^{[2\xi]},T_{\Delta}^{[2\xi]})$ and $q_i^{\xi}(T_{\Delta}^{[2\xi]})$ are defined as $\sigma _{ij}^{\xi}(\bm{u}_{\Delta }^{[2\xi]},T_{\Delta}^{[2\xi]})={a_{ijkl}^\xi}\varepsilon _{kl}^{\xi}(\bm{u}_{\Delta }^{[2\xi]}) - {b_{ij}^\xi}T_{\Delta}^{[2\xi]}$ and $q_i^{\xi}(T_{\Delta}^{[2\xi]})=- {k_{ij}^\xi}{\Psi _j}(T_{\Delta}^{[2\xi]})$ in (3.3). And the detailed expressions of functions $\mathbb{G}_i({\bm{\alpha }},{\bm{\beta }},t)$ and $\mathbb{E}({\bm{\alpha }},{\bm{\beta }},t)$ are also uncomplicated to attain and not displayed in the present study owing to its lengthy forms.

Through the local error analysis in the pointwise sense, we can demonstrate that LOMS solutions can not maintain the conservation of local energy and momentum since the $\xi$-independent terms $\mathbb{F}_{0i}({\bm{\alpha }},{\bm{\beta }},t)$ and $\mathbb{S}_0({\bm{\alpha }},{\bm{\beta }},t)$ in (3.2) can not approach to zero along the variable $\xi$ approaches to zero. Benefitting from the higher-order corrected terms, the HOMS solutions can guarantee the local momentum conservation of mechanical equations and local energy conservation of thermal equation in the original governing equations (1.1) due to their $O(\xi)$-order pointwise errors. This is the principal motivation for this study to establish the HOMS solutions exhibiting high-accuracy computation performance for heterogeneous shells.
\subsection{The convergence proof by global error estimation}
To obtain global error estimation of higher-order multi-scale approximate solutions, three main difficulties was firstly claimed. The first challenge of global error estimation is that the nonlinear two-way thermo-mechanical coupling term $\displaystyle\vartheta _{ij}^\xi({\bm{\alpha }},T^{\xi}) \frac{{\partial \varepsilon _{ij}^\xi({\bm{\alpha }},t) }}{{\partial t}}$ makes it impossible to incorporate the elastic-dynamics equations and the transient heat conduction equation into single equation for error estimation. The second challenge of global error estimation is that lack of a prior estimate for elastic-dynamics equations in time-dependent nonlinear thermo-mechanical system with mixed boundary condition. The third challenge of global error estimation is that the Laplace transform technique for linear thermo-mechanical problems in \cite{R4,R54,R55} is invalid due to its nonlinear equation coefficients. Consideration of the global error estimation for the HOMS solutions is one novelty of the present study.

For gaining the global error estimation in the energy norm sense, some assumptions are further presented as below.
\begin{enumerate}
\item[(i)]
Assume that $\Omega$ is the union of the entire periodic cells, i.e. $\bar{\Omega}=\cup_{\mathbf{z}\in I_{\xi}}\xi(\mathbf{z}+\bar{\Theta})$, where the index set $I_{\xi}=\{\mathbf{z}=(z_1,z_2,z_3)\in Z^3,\xi(\mathbf{z}+\bar{\Theta})\subset \bar{\Omega}\}$. Besides, let $E_\mathbf{z}=\xi(\mathbf{z}+\Theta)$ be
translational unit cell and $\partial E_\mathbf{z}$ be the boundary of $E_\mathbf{z}$.
\item[(ii)]
Employ the homogeneous Dirichlet boundary condition to displace the periodic boundary condition for whole auxiliary cell functions \cite{R30,R56,R57}.
\item[(iii)]
Let $b_{ij}(\bm{\beta},{T^\xi})=b_{ii}(\bm{\beta},{T^\xi})\delta_{ij}$, $k_{ij}(\bm{\beta},{T^\xi})=k_{ii}(\bm{\beta},{T^\xi})\delta_{ij}$ and $\vartheta_{ij}(\bm{\beta},{T^\xi})=\vartheta_{ii}(\bm{\beta},{T^\xi})\delta_{ij}$, and $\delta_{ij}$ is a Kronecker symbol. Moreover, let $\Delta_1,\Delta_2,\Delta_3$ be the middle hyperplanes of PUC $\Theta=(0,1)^3$. Then suppose that material parameters $\rho(\bm{\beta},{T^\xi})$, $a_{ijkl}(\bm{\beta},{T^\xi})$, $b_{ii}(\bm{\beta},{T^\xi})$, $c(\bm{\beta},{T^\xi})$, $k_{ii}(\bm{\beta},{T^\xi})$ and $\vartheta_{ii}(\bm{\beta},{T^\xi})$ are symmetric with respect to $\Delta_1,\Delta_2,\Delta_3$ for stationary ${T^\varepsilon }\in [T_{min},T_{max}+C_*]$.
\item[(iv)]
The nonlinear term $\displaystyle\vartheta _{ij}^{\xi}=T^{\xi}\big[b _{ij}^{\xi}+\frac{\partial b_{ij}^{\xi}}{\partial T^{\xi}}(T^{\xi}-\widetilde{T})\big]$ in time-dependent nonlinear thermo-mechanical equations (1.1) is simplified as linear term $\displaystyle\vartheta _{ij}^{\xi}=\widetilde Tb_{ij}^{\xi}$.
\item[(v)]
Assume that $\displaystyle\frac{\partial \rho^{\xi}({\bm{\alpha }},T^{\xi})}{\partial t}$, $\displaystyle\frac{\partial a^{\xi}_{ijkl}({\bm{\alpha }},T^{\xi})}{\partial t}$ and $\displaystyle\frac{\partial c^{\xi}({\bm{\alpha }},T^{\xi})}{\partial t}\in L^\infty(\Omega\times(0,\mathcal T))$.
\end{enumerate}

On the basis of the above assumptions, two important lemmas are obtained, which will be utilized to conduct global error estimation of the HOMS solutions.
\begin{lemma}
Denote the derivative operators $\displaystyle\sigma_{i \Theta}(\bm{\phi})=n_j a_{ijkl}(\bm{\beta},{T^\varepsilon })\breve\varepsilon_{kl}^{[0]}(\bm{\phi})$ and $\displaystyle\sigma_{T \Theta}(\chi)=n_i k_{ij}(\bm{\beta},{T^\varepsilon })\Psi_j(\chi)$. Then under assumptions (A)-(B) and (ii)-(iii), the normal derivatives $\sigma_{i\Theta}(\mathbf{N}^{mn})$, $\sigma_{i\Theta}(\mathbf{P})$, $\sigma_{i\Theta}(\mathbf{N}^{jmn})$, $\sigma_{i\Theta}(\mathbf{H}^j)$, $\sigma_{i\Theta}(\mathbf{F}^j)$, $\sigma_{i\Theta}(\mathbf{W}^{mn})$, $\sigma_{i\Theta}(\mathbf{Q})$, $\sigma_{i\Theta}(\mathbf{Z}^j)$, $\sigma_{i\Theta}(\mathbf{X}^{jmn})$, $\sigma_{T\Theta}(M^m)$, $\sigma_{T\Theta}(M^{mn})$, $\sigma_{T\Theta}(A)$, $\sigma_{T\Theta}(G^{mn})$, $\sigma_{T\Theta}(R^m)$ and $\sigma_{T\Theta}(B^{mn})$ are proved to be continuous on the boundary of PUC $\Theta$ via employing the identical approach in Refs. \cite{R30,R56,R57}. Furthermore, $\sigma_{i\Theta}(\bm{u}^{[2\xi]})$ and $\sigma_{T\Theta}(T^{[2\xi]})$ can be proved to be continuous on the boundary of PUC $\Theta$.
\end{lemma}

\begin{lemma}
Under assumptions (i)-(ii) and (iv), the residual equations (3.3) of the HOMS solutions are further presented and attached with initial-boundary conditions as follows.
\begin{equation}
\left\{
\begin{aligned}
&{\rho ^\xi }\frac{{\partial ^2}u_{\Delta i}^{[2\xi]}}{{\partial t^2}}-\sum\limits_{j = 1}^3 {\frac{1}{H}\frac{\partial }{{\partial {\alpha _j}}}\left( {\frac{H}{{{H_j}}}\sigma _{ij}^\xi (\bm{u}_{\Delta }^{[2\xi]},T_{\Delta}^{[2\xi]})} \right)} \\
&\quad\quad\quad\quad\;-\sum\limits_{j \ne i,j = 1}^3 {\frac{1}{{{H_i}{H_j}}}\frac{{\partial {H_i}}}{{\partial {\alpha _j}}}} \sigma _{ij}^\xi(\bm{u}_{\Delta }^{[2\xi]},T_{\Delta}^{[2\xi]})\\
&\quad\quad\quad\quad\; + \sum\limits_{j \ne i,j = 1}^3 {\frac{1}{{{H_i}{H_j}}}\frac{{\partial {H_j}}}{{\partial {\alpha _i}}}} \sigma _{jj}^\xi(\bm{u}_{\Delta }^{[2\xi]},T_{\Delta}^{[2\xi]}) \\
&\quad\quad\quad\quad\; = \xi \mathbb{G}_i({\bm{\alpha }},{\bm{\beta }},t),\;\;\text{in}\;\;\Omega\times(0,\mathcal T),\\
&{{\rho ^\xi }{c^\xi }\frac{{\partial {T_{\Delta}^{[2\xi]}}}}{{\partial t}}{\rm{ + }}}\sum\limits_{i = 1}^3 {\frac{1}{H}\frac{\partial }{{\partial {\alpha _i}}}\left( {\frac{H}{{{H_i}}}q_i^\xi(T_{\Delta}^{[2\xi]}) } \right)}{ + \widetilde Tb_{ij}^{\xi}\frac{{\partial \varepsilon_{ij}^\xi (\bm{u}_{\Delta }^{[2\xi]})}}{{\partial t}}} \\
&\quad\quad\quad\quad\;\; = \xi \mathbb{E}({\bm{\alpha }},{\bm{\beta }},t),\;\;\text{in}\;\;\Omega\times(0,\mathcal T),\\
&u_{\Delta i}^{[2\xi]}({\bm{\alpha }},t)=0,\;\;\text{on}\;\;\partial\Omega_{u}\times(0,\mathcal T),\\
&T_{\Delta}^{[2\xi]}({\bm{\alpha }},t)=0,\;\;\text{on}\;\;\partial\Omega_{T}\times(0,\mathcal T),\\
&u_{\Delta i}^{[2\xi]}(\bm{\alpha},0)=\xi \widehat{\varphi}_{2i}(\bm{\alpha}),\;\;\frac{\partial u_{\Delta i}^{[2\xi]}(\bm{\alpha},t)}{\partial t}\big|_{t=0}=\xi \widehat{\varpi}_{2i}(\bm{\alpha}),\;\;\text{in}\;\;\Omega,\\
&\;T_{\Delta}^{[2\xi]}(\bm{\alpha},0)=\xi \widehat{\omega}_2(\bm{\alpha}),\;\;\text{in}\;\;\Omega.
\end{aligned} \right.
\end{equation}
\end{lemma}

Next, we give the final result of global error estimation for the HOMS solutions of the time-dependent multi-scale nonlinear equations (1.1) as the following theorem.
\begin{theorem}
Let $\bm{u}^\xi(\bm{\alpha},t)$ and $T^\xi(\bm{\alpha},t)$ be the weak solutions of multi-scale nonlinear equations (1.1), $\bm{u}^{[0]}(\bm{\alpha},t)$ and $T^{[0]}(\bm{\alpha},t)$ be the weak solutions of corresponding homogenized equations (2.18), $\bm{u}^{[2\xi]}(\bm{\alpha},t)$ and $T^{[2\xi]}(\bm{\alpha},t)$ be the HOMS solutions stated in formulas (2.35) and (2.36). Under the above hypotheses (A)-(C) and (i)-(v), the following global error estimation are obtained.
\begin{equation}
{\begin{aligned}
&\Big\| \frac{\partial\bm{u}_{\Delta }^{[2\xi]}(\bm{\alpha},t)}{\partial t}\Big\|_{L^\infty(0,\mathcal T;(L^2(\Omega))^3)}+
\big\| \bm{u}_{\Delta }^{[2\xi]}(\bm{\alpha},t)\big\|_{L^\infty(0,\mathcal T;(H^1_0(\Omega))^3)}\\
&+\big\| T_{\Delta}^{[2\xi]}(\bm{\alpha},t)\big\|_{L^\infty(0,\mathcal T;L^2(\Omega))}+
\big\| T_{\Delta}^{[2\xi]}(\bm{\alpha},t)\big\|_{L^2(0,\mathcal T;H^1_0(\Omega))}\leq C(\Omega,\mathcal T)\xi,
\end{aligned}}
\end{equation}
\end{theorem}
where $C(\Omega,\mathcal T)$ is a positive constant irrespective of $\xi$, but dependent of $\Omega$ and $\mathcal T$.\\
$\mathbf{Proof:}$ The residual equations (3.4) in lemma 3.2 are employed to accomplish the global error estimation. Since $\displaystyle\frac{\partial u_{\Delta i}^{[2\xi]}}{\partial t}$ belongs to $L^\infty(0,\mathcal T;L^2(\Omega))$, $\displaystyle\frac{\partial u_{\Delta i}^{[2\xi]}}{\partial t}$ cannot be directly utilized as test function for hyperbolic parts of residual equations (3.4). To tackle this challenge, the so-called density argument approach in \cite{R18,R56} is employed. For simplifying the proof, the details about density argument approach are omitted. We further suppose that the equation parameter $H\in[\kappa_1,\kappa_2]$, where $\kappa_1$ and $\kappa_2$ are defined as the minimum and maximum $H$ of inhomogeneous shells separately. In what follows, multiplying $\displaystyle\frac{\partial u_{\Delta i}^{[2\xi]}}{\partial t}H$ and $\displaystyle T_{\Delta}^{[2\xi]}H$ on both sides of (3.4) and integrating over $\Omega$, we can derive following equalities
\begin{equation}
\left\{ 
\begin{aligned}
&{\int_{\Omega}}\rho^{\xi}({\bm{\alpha }},T^{\xi}) \frac{{{\partial ^2}{u_{\Delta i}^{[2\xi]}}}}{{\partial {t^2}}}\frac{\partial u_{\Delta i}^{[2\xi]}}{\partial t}Hd\Omega - {\int_{\Omega}}\Big[\sum\limits_{j = 1}^3 {\frac{1}{H}\frac{\partial }{{\partial {\alpha _j}}}\Big( {\frac{H}{{{H_j}}}\sigma _{ij}^\xi (\bm{u}_{\Delta }^{(2\xi)},T_{\Delta}^{[2\xi]})} \Big)}\\
&+\sum\limits_{j \ne i,j = 1}^3 {\frac{1}{{{H_i}{H_j}}}\frac{{\partial {H_i}}}{{\partial {\alpha _j}}}} \sigma _{ij}^\xi(\bm{u}_{\Delta }^{[2\xi]},T_{\Delta}^{[2\xi]})\\
&-\sum\limits_{j \ne i,j = 1}^3 {\frac{1}{{{H_i}{H_j}}}\frac{{\partial {H_j}}}{{\partial {\alpha _i}}}} \sigma _{jj}^\xi(\bm{u}_{\Delta }^{[2\xi]},T_{\Delta}^{[2\xi]})\Big]\frac{\partial u_{\Delta i}^{[2\xi]}}{\partial t}Hd\Omega\\
& = {\int_{\Omega}}{\xi \mathbb{G}_i({\bm{\alpha }},{\bm{\beta }},t)}\frac{\partial u_{\Delta i}^{[2\xi]}}{\partial t}Hd\Omega,\\
&{\int_{\Omega}}\rho^{\xi}({\bm{\alpha }},T^{\xi}) c^{\xi}({\bm{\alpha }},T^{\xi})\frac{{\partial T_{\Delta}^{[2\xi]}}}{{\partial t}}T_{\Delta}^{[2\xi]}Hd\Omega +{\int_{\Omega}}\Big[\sum\limits_{i = 1}^3 {\frac{1}{H}\frac{\partial }{{\partial {\alpha _i}}}\Big( {\frac{H}{{{H_i}}}q_i^\xi(T_{\Delta}^{[2\xi]}) } \Big)}\\
&{ + \widetilde T b_{ij}^\xi({\bm{\alpha }},T^{\xi}) \frac{{\partial \varepsilon_{ij}^\xi (\bm{u}_{\Delta }^{[2\xi]})}}{\partial t}}\Big]T_{\Delta}^{[2\xi]}Hd\Omega={\int_{\Omega}}\xi \mathbb{E}({\bm{\alpha }},{\bm{\beta }},t)T_{\Delta}^{[2\xi]}Hd\Omega.
\end{aligned} \right.
\end{equation}
Integrating by parts on (3.6) and employing Green's formula, the upper equality could be simplified as follows
\begin{equation}
\left\{ 
\begin{aligned}
&{\int_{\Omega}}\rho^{\xi}({\bm{\alpha }},T^{\xi}) \frac{{{\partial ^2}{u_{\Delta i}^{[2\xi]}}}}{{\partial {t^2}}}\frac{\partial u_{\Delta i}^{[2\xi]}}{\partial t}Hd\Omega\\
&+{\int_{\Omega}}\left[ {a_{ijkl}^{\xi}({\bm{\alpha }},T^{\xi})\varepsilon^\xi_{kl}(\bm{u}_{\Delta }^{[2\xi]})-b_{ij}^\xi({\bm{\alpha }},T^{\xi})}{ T_{\Delta}^{[2\xi]}}\right ]\varepsilon^\xi_{ij}(\frac{\partial \bm{u}_{\Delta }^{[2\xi]}}{\partial t})Hd\Omega\\
&= {\int_{\Omega}}{\xi \mathbb{G}_{i}}\frac{\partial u_{\Delta i}^{[2\xi]}}{\partial t}Hd\Omega+{\int_{\cup_{\mathbf{z}\in I_\xi}\partial E_\mathbf{z}}}\sigma_{i\Theta}(\bm{u}_{\Delta }^{[2\xi]})\frac{\partial u_{\Delta i}^{[2\xi]}}{\partial t}Hd\Gamma_{\bm{\beta}},\\
&{\int_{\Omega}}\rho^{\xi}({\bm{\alpha }},T^{\xi}) c^{\xi}({\bm{\alpha }},T^{\xi})\frac{{\partial T_{\Delta}^{[2\xi]}}}{{\partial t}}T_{\Delta}^{[2\xi]}Hd\Omega\\
&-{\int_{\Omega}}{q_i^{\xi}(T_{\Delta}^{[2\xi]}) }{\Psi _i}(T_{\Delta}^{[2\xi]})Hd\Omega+{\int_{\Omega}}{\widetilde Tb_{ij}^\xi({\bm{\alpha }},T^{\xi})}\frac{\partial \varepsilon^\xi_{ij}(\bm{u}_{\Delta }^{[2\xi]})}{\partial t}T_{\Delta}^{[2\xi]}Hd\Omega\\
&={\int_{\Omega}} \xi \mathbb{E}T_{\Delta}^{[2\xi]}Hd\Omega+{\int_{\cup_{\mathbf{z}\in I_\xi}\partial E_\mathbf{z}}}\sigma_{T\Theta}(T_{\Delta}^{[2\xi]})T_{\Delta}^{[2\xi]}Hd\Gamma_{\bm{\beta}},
\end{aligned} \right.
\end{equation}
where $\sigma_{i\Theta}(\bm{u}_{\Delta }^{[2\xi]})$ and $\sigma_{T\Theta}(T_{\Delta}^{[2\xi]})$ arise from employing Green's formula on the boundary $\partial E_\mathbf{z}$.

Recalling lemma 3.1, we shall hereby get the following results for the integral terms on the boundary $\partial E_\mathbf{z}$ in (3.7).
\begin{equation}
\left\{ 
\begin{aligned}
&{\int_{\cup_{\mathbf{z}\in I_\xi}\partial E_\mathbf{z}}}\sigma_{i\Theta}(\bm{u}_{\Delta }^{[2\xi]})\frac{\partial u_{\Delta i}^{[2\xi]}}{\partial t}Hd\Gamma_{\bm{\beta}}=\sum\limits _{\mathbf{z}\in I_\xi}\int_{\partial E_\mathbf{z}}\sigma_{i\Theta}(\bm{u}^\xi-\bm{u}^{[2\xi]})\frac{\partial u_{\Delta i}^{[2\xi]}}{\partial t}d\Gamma_{\bm{\beta}}\\
&=-\sum\limits _{\mathbf{z}\in I_\xi}\int_{\partial E_\mathbf{z}}\sigma_{i\Theta}(\bm{u}^{[2\xi]})\frac{\partial u_{\Delta i}^{[2\xi]}}{\partial t}d\Gamma_{\bm{\beta}}=0,\\
&{\int_{\cup_{\mathbf{z}\in I_\xi}\partial E_\mathbf{z}}}\sigma_{T\Theta}(T_{\Delta}^{[2\xi]})T_{\Delta}^{[2\xi]}Hd\Gamma_{\bm{\beta}}=\sum\limits _{\mathbf{z}\in I_\xi}\int_{\partial E_\mathbf{z}}\sigma_{T\Theta}(T^{\xi}-T^{[2\xi]})T_{\Delta}^{[2\xi]}d\Gamma_{\bm{\beta}}
\\
&=-\sum\limits _{\mathbf{z}\in I_\xi}\int_{\partial E_\mathbf{z}}\sigma_{T\Theta}(T^{[2\xi]})T_{\Delta}^{[2\xi]}d\Gamma_{\bm{\beta}}=0.
\end{aligned} \right.
\end{equation}

In the following, two equalities are easily derived by combining (3.7) and (3.8) together as below
\begin{align}
\begin{array}{*{20}{l}}
\displaystyle\frac{1}{2}\frac{\partial}{\partial t}\Big[\int_{\Omega}\rho^{\xi}({\bm{\alpha }},T^{\xi})(\frac{\partial u_{\Delta i}^{[2\xi]}}{\partial t})^2Hd\Omega+\int_{\Omega} {a_{ijkl}^{\xi}({\bm{\alpha }},T^{\xi})\varepsilon^\xi_{kl}(\bm{u}_{\Delta }^{[2\xi]})} \varepsilon^\xi_{ij}(\bm{u}_{\Delta }^{[2\xi]})Hd\Omega\Big]\\
-\displaystyle\frac{1}{2}\int_{\Omega}\frac{\partial \rho^{\xi}({\bm{\alpha }},T^{\xi})}{\partial t}(\frac{\partial u_{\Delta i}^{[2\xi]}}{\partial t})^2Hd\Omega\\
\displaystyle-\frac{1}{2}\int_{\Omega} {\frac{\partial a_{ijkl}^{\xi}({\bm{\alpha }},T^{\xi})}{\partial t}\varepsilon^\xi_{kl}(\bm{u}_{\Delta }^{[2\xi]})} \varepsilon^\xi_{ij}(\bm{u}_{\Delta }^{[2\xi]})Hd\Omega\\
\displaystyle -\int_{\Omega}{b_{ij}^\xi({\bm{\alpha }},T^{\xi})}{ T_{\Delta}^{[2\xi]}}\varepsilon^\xi_{ij}(\frac{\partial \bm{u}_{\Delta }^{[2\xi]}}{\partial t})Hd\Omega=\int_{\Omega}{\xi \mathbb{G}_{i}}\frac{\partial u_{\Delta i}^{[2\xi]}}{\partial t}Hd\Omega.
\end{array}
\end{align}
\begin{align}
\begin{array}{*{20}{l}}
\displaystyle\frac{1}{2}\frac{\partial}{\partial t}\Big[\int_{\Omega}\rho^{\xi}({\bm{\alpha }},T^{\xi}) c^{\xi}({\bm{\alpha }},T^{\xi})(T_{\Delta}^{[2\xi]})^2Hd\Omega\Big]-\int_{\Omega}{q_i^{\xi}(T_{\Delta}^{[2\xi]}) }{\Psi _i}(T_{\Delta}^{[2\xi]})Hd\Omega \\
\displaystyle-\frac{1}{2}\int_{\Omega}\frac{\partial \rho^{\xi}({\bm{\alpha }},T^{\xi})}{\partial t} c^{\xi}({\bm{\alpha }},T^{\xi})(T_{\Delta}^{[2\xi]})^2Hd\Omega\\
\displaystyle-\frac{1}{2}\int_{\Omega}\rho^{\xi}({\bm{\alpha }},T^{\xi}) \frac{\partial c^{\xi}({\bm{\alpha }},T^{\xi})}{\partial t}(T_{\Delta}^{[2\xi]})^2Hd\Omega\\
\displaystyle+\int_{\Omega}{\widetilde Tb_{ij}^\xi({\bm{\alpha }},T^{\xi})}\frac{\partial \varepsilon^\xi_{ij}(\bm{u}_{\Delta }^{[2\xi]})}{\partial t}T_{\Delta}^{[2\xi]}Hd\Omega= \int_{\Omega} \xi \mathbb{E}T_{\Delta}^{[2\xi]}Hd\Omega.
\end{array}
\end{align}
Subsequently, implementing the calculation (3.10)+$\widetilde T\times$(3.9), we get the equality as below
\begin{align}
\begin{array}{*{20}{l}}
\displaystyle\frac{1}{2}\frac{\partial}{\partial t}\Big[\int_{\Omega}\rho^{\xi}(\frac{\partial u_{\Delta i}^{[2\xi]}}{\partial t})^2\widetilde THd\Omega+\int_{\Omega} {a_{ijkl}^{\xi}\varepsilon^\xi_{kl}(\bm{u}_{\Delta }^{[2\xi]})} \varepsilon^\xi_{ij}(\bm{u}_{\Delta }^{[2\xi]})\widetilde THd\Omega\Big]\\
-\displaystyle\frac{1}{2}\int_{\Omega}\frac{\partial \rho^{\xi}}{\partial t}(\frac{\partial u_{\Delta i}^{[2\xi]}}{\partial t})^2Hd\Omega-\frac{1}{2}\int_{\Omega} {\frac{\partial a_{ijkl}^{\xi}}{\partial t}\varepsilon^\xi_{kl}(\bm{u}_{\Delta }^{[2\xi]})} \varepsilon^\xi_{ij}(\bm{u}_{\Delta }^{[2\xi]})Hd\Omega\\
\displaystyle+\frac{1}{2}\frac{\partial}{\partial t}\Big[\int_{\Omega}\rho^{\xi} c^{\xi}(T_{\Delta}^{[2\xi]})^2Hd\Omega\Big]-\int_{\Omega}{q_i^{\xi}(T_{\Delta}^{[2\xi]}) }{\Psi _i}(T_{\Delta}^{[2\xi]})Hd\Omega\\
\displaystyle-\frac{1}{2}\int_{\Omega}\frac{\partial \rho^{\xi}}{\partial t} c^{\xi}(T_{\Delta}^{[2\xi]})^2Hd\Omega-\frac{1}{2}\int_{\Omega}\rho^{\xi} \frac{\partial c^{\xi}}{\partial t}(T_{\Delta}^{[2\xi]})^2Hd\Omega\\
= \displaystyle\int_{\Omega}{\xi\mathbb{G}_{i}}\frac{\partial u_{\Delta i}^{[2\xi]}}{\partial t}\widetilde THd\Omega \displaystyle+\int_{\Omega} \xi \mathbb{E}T_{\Delta}^{[2\xi]}Hd\Omega.
\end{array}
\end{align}
Afterwards, integrating both sides of (3.11) from $0$ to $t$ $(0<t\leq \mathcal T)$, and also putting the initial and boundary conditions of equations (3.4) into the upper equality (3.11), one can deduce the following equality
\begin{align}
\begin{array}{*{20}{l}}
\displaystyle\int_{\Omega}\rho^{\xi}(\frac{\partial u_{\Delta i}^{[2\xi]}(\bm{\alpha},t)}{\partial t})^2\widetilde THd\Omega +\int_{\Omega}\rho^{\xi} c^{\xi}(T_{\Delta}^{[2\xi]}(\bm{\alpha},t))^2Hd\Omega\\
+\displaystyle\int_{\Omega} {a_{ijkl}^{\xi}\varepsilon^\xi_{kl}(\bm{u}_{\Delta }^{[2\xi]}(\bm{\alpha},t))} \varepsilon^\xi_{ij}(\bm{u}_{\Delta }^{[2\xi]}(\bm{\alpha},t))\widetilde THd\Omega\\
\displaystyle
-2\int_0^t\int_{\Omega}{q_i^{\xi}(T_{\Delta}^{[2\xi]}(\bm{\alpha},\tau)) }{\Psi _i}(T_{\Delta}^{[2\xi]}(\bm{\alpha},\tau))Hd\Omega d\tau\\
\displaystyle=\int_{\Omega}\rho^{\xi}(\xi \widehat{\varpi}_{2i}(\bm{\alpha}))^2\widetilde THd\Omega+\int_{\Omega}\rho^{\xi} c^{\xi}(\xi \widehat{\omega}_2(\bm{\alpha}))^2Hd\Omega\\
\displaystyle+\int_{\Omega} {a_{ijkl}^{\xi}\varepsilon^\xi_{kl}(\xi \widehat{\varphi}_{2i}(\bm{\alpha}))} \varepsilon^\xi_{ij}(\xi \widehat{\varphi}_{2i}(\bm{\alpha}))\widetilde THd\Omega\\ \displaystyle+\int_0^t\int_{\Omega}\frac{\partial \rho^{\xi}}{\partial \tau}(\frac{\partial u_{\Delta i}^{[2\xi]}}{\partial \tau})^2Hd\Omega d\tau+\int_0^t\int_{\Omega} {\frac{\partial a_{ijkl}^{\xi}}{\partial \tau}\varepsilon^\xi_{kl}(\bm{u}_{\Delta }^{[2\xi]})} \varepsilon^\xi_{ij}(\bm{u}_{\Delta }^{[2\xi]})Hd\Omega d\tau\\
\displaystyle+\int_0^t\int_{\Omega}\frac{\partial \rho^{\xi}}{\partial \tau} c^{\xi}(T_{\Delta}^{[2\xi]})^2Hd\Omega d\tau+\int_0^t\int_{\Omega}\rho^{\xi} \frac{\partial c^{\xi}}{\partial \tau}(T_{\Delta}^{[2\xi]})^2Hd\Omega d\tau\\
\displaystyle+2\int_0^t\int_{\Omega}{\xi \mathbb{G}_{i}}\frac{\partial u_{\Delta i}^{[2\xi]}(\bm{\alpha},\tau)}{\partial \tau}\widetilde THd\Omega d\tau+2\int_0^t\int_{\Omega} \xi \mathbb{E}T_{\Delta}^{[2\xi]}(\bm{\alpha},\tau)Hd\Omega d\tau.
\end{array}
\end{align}
So far, we gain the crucial equality (3.12) for proving theorem 3.3. Initiating from here, the final proof will be obtained.

By the hypotheses (A) and (B), and assuming $H\in[\kappa_1,\kappa_2]$, the following inequality can be obtained by utilizing Korn's inequality and Poincar$\rm{\acute{e}}$-Friedrichs inequality in curvilinear coordinates \cite{R2,R3} on left side of equation (3.12)
\begin{equation}
\begin{aligned}
&\int_{\Omega}\rho^{\xi}(\frac{\partial u_{\Delta i}^{[2\xi]}(\bm{\alpha},t)}{\partial t})^2\widetilde THd\Omega+\int_{\Omega} {a_{ijkl}^{\xi}\varepsilon^\xi_{kl}(\bm{u}_{\Delta }^{[2\xi]}(\bm{\alpha},t))} \varepsilon^\xi_{ij}(\bm{u}_{\Delta }^{[2\xi]}(\bm{\alpha},t))\widetilde THd\Omega\\
&+\int_{\Omega}\rho^{\xi} c^{\xi}(T_{\Delta}^{[2\xi]}(\bm{\alpha},t))^2Hd\Omega
-2\int_0^t\int_{\Omega}{q_i^{\xi}(T_{\Delta}^{[2\xi]}(\bm{\alpha},\tau))}{\Psi _i}(T_{\Delta}^{[2\xi]}(\bm{\alpha},\tau))Hd\Omega d\tau\\
&\geq\rho^{*}\widetilde T\kappa_1\Big \| \frac{\partial \bm{u}_{\Delta }^{[2\xi]}(\bm{\alpha},t)}{\partial t}\Big\|_{(L^2(\Omega))^3}^2+\widetilde T\kappa_1C_1\left \| \bm{u}_{\Delta }^{[2\xi]}(\bm{\alpha},t)\right\|_{(H^1_0(\Omega))^3}^2\\
&+\rho^{*}c^{*}\kappa_1\left \| T_{\Delta}^{[2\xi]}(\bm{\alpha},t)\right\|_{L^2(\Omega)}^2
+\kappa_1C_2\int_0^t\left \| T_{\Delta}^{[2\xi]}(\bm{\alpha},\tau)\right\|_{H^1_0(\Omega)}^2d\tau\\
&\geq \lambda_1\Big(\Big \| \frac{\partial \bm{u}_{\Delta }^{[2\xi]}}{\partial t}\Big\|_{(L^2(\Omega))^3}^2+\left \| \bm{u}_{\Delta }^{[2\xi]}\right\|_{(H^1_0(\Omega))^3}^2\\
&+\left \| T_{\Delta}^{[2\xi]}\right\|_{L^2(\Omega)}^2
+\int_0^t\left \| T_{\Delta}^{[2\xi]}\right\|_{H^1_0(\Omega)}^2d\tau\Big).
\end{aligned}
\end{equation}
Here $C_1$ and $C_2$ arise from Korn's inequality and Poincar$\rm{\acute{e}}$-Friedrichs inequality in curvilinear coordinates separately, and $\lambda_1=min(\rho^{*}\widetilde T\kappa_1,\widetilde T\kappa_1C_1,\rho^{*}c^{*}\kappa_1,\kappa_1C_2)$.

To implement, employing Schwarz's inequality and Young's inequality \cite{R18,R19} on the right side of equation (3.12), a direct computation gives
\begin{equation}
\begin{aligned}
&\int_{\Omega}\rho^{\xi}(\xi \widehat{\varpi}_{2i}(\bm{\alpha}))^2\widetilde THd\Omega+\int_{\Omega} {a_{ijkl}^{\xi}\varepsilon^\xi_{kl}(\xi \widehat{\varphi}_{2i}(\bm{\alpha}))} \varepsilon^\xi_{ij}(\xi \widehat{\varphi}_{2i}(\bm{\alpha}))\widetilde THd\Omega\\
&\displaystyle+\int_{\Omega}\rho^{\xi} c^{\xi}(\xi \widehat{\omega}_2(\bm{\alpha}))^2Hd\Omega+\int_0^t\int_{\Omega}\frac{\partial \rho^{\xi}}{\partial \tau}(\frac{\partial u_{\Delta i}^{[2\xi]}}{\partial \tau})^2Hd\Omega d\tau\\
&\displaystyle+\int_0^t\int_{\Omega} {\frac{\partial a_{ijkl}^{\xi}}{\partial \tau}\varepsilon^\xi_{kl}(\bm{u}_{\Delta }^{[2\xi]})} \varepsilon^\xi_{ij}(\bm{u}_{\Delta }^{[2\xi]})Hd\Omega d\tau\\
&\displaystyle+\int_0^t\int_{\Omega}\frac{\partial \rho^{\xi}}{\partial \tau} c^{\xi}(T_{\Delta}^{[2\xi]})^2Hd\Omega d\tau+\int_0^t\int_{\Omega}\rho^{\xi} \frac{\partial c^{\xi}}{\partial \tau}(T_{\Delta}^{[2\xi]})^2Hd\Omega d\tau\\
&\displaystyle+2\int_0^t\int_{\Omega}{\xi \mathbb{G}_{i}}\frac{\partial u_{\Delta i}^{[2\xi]}(\bm{\alpha},\tau)}{\partial \tau}\widetilde THd\Omega d\tau+2\int_0^t\int_{\Omega} \xi \mathbb{E}T_{\Delta}^{[2\xi]}(\bm{\alpha},\tau)Hd\Omega d\tau\\
&\leq C(\Omega)\xi^2+C(\Omega)\Big[\int_0^t\Big\| \frac{\partial \bm{u}_{\Delta }^{[2\xi]}}{\partial \tau}\Big\|_{(L^2(\Omega))^3}^2d\tau+\int_0^t\left \| \bm{u}_{\Delta }^{[2\xi]}\right\|_{(H^1_0(\Omega))^3}^2d\tau\\ &+\int_0^t\left \| T_{\Delta}^{[2\xi]}\right\|_{L^2(\Omega)}^2d\tau+\int_0^t\int_0^\tau\left \| T_{\Delta}^{[2\xi]}\right\|_{H^1_0(\Omega)}^2dsd\tau \Big].
\end{aligned}
\end{equation}
Here $C(\Omega)$ defines a positive constant, which is dependent of $\Omega$. Moreover, it should be highlighted that $C(\Omega)$ will have distinct values in different positions in this study.

Through combining (3.13) and (3.14) together, we then get $\lambda_1\digamma(t)\leq C(\Omega)\xi^2+C(\Omega)\mathlarger{\int}_0^t\digamma(\tau)d\tau$ by setting $\displaystyle\digamma(t)=\big \| \frac{\partial \bm{u}_{\Delta }^{[2\xi]}}{\partial t}\big\|_{(L^2(\Omega))^3}^2+\big \| \bm{u}_{\Delta }^{[2\xi]}\big\|_{(H^1_0(\Omega))^3}^2+\big \| T_{\Delta}^{[2\xi]}\big\|_{L^2(\Omega)}^2
+\mathlarger{\int}_0^t\big \| T_{\Delta}^{[2\xi]}\big\|_{H^1_0(\Omega)}^2d\tau$. Without loss of generality, denoting $\displaystyle C(\Omega)=C(\Omega)/\lambda_1$, it follows from Gronwall's inequality \cite{R18} that $\digamma(t)\leq C(\Omega,\mathcal T)\xi^2$. Subsequently, there holds the following inequality
\begin{equation}
\begin{aligned}
&\Big \| \frac{\partial \bm{u}_{\Delta }^{[2\xi]}}{\partial t}\Big\|_{(L^2(\Omega))^3}^2+\big \| \bm{u}_{\Delta }^{[2\xi]}\big\|_{(H^1_0(\Omega))^3}^2\\
&+\big \| T_{\Delta}^{[2\xi]}\big\|_{L^2(\Omega)}^2
+{\int}_0^t\big \| T_{\Delta}^{[2\xi]}\big\|_{H^1_0(\Omega)}^2d\tau\leq C(\Omega,\mathcal T)\xi^2.
\end{aligned}
\end{equation}
To proceed, employing the AM-GM inequality $\displaystyle\frac{a+b+c+d}{4}\leq\sqrt{\frac{a^2+b^2+c^2+d^2}{4}}$ to the left side of (3.15) and then squaring root on both sides of (3.15), we can deduce that
\begin{equation}
\begin{aligned}
&\Big \| \frac{\partial \bm{u}_{\Delta }^{[2\xi]}}{\partial t}\Big\|_{(L^2(\Omega))^3}+\big \| \bm{u}_{\Delta }^{[2\xi]}\big\|_{(H^1_0(\Omega))^3}\\
&+\big \| T_{\Delta}^{[2\xi]}\big\|_{L^2(\Omega)}
+\big\| T_{\Delta}^{[2\xi]}\big\|_{L^2(0,t;H^1_0(\Omega))}\leq C(\Omega,\mathcal T)\xi.
\end{aligned}
\end{equation}
By the arbitrariness of time variable $t$, employing (3.16) yields the final error estimation as follows
\begin{equation}
{\begin{aligned}
&\Big\| \frac{\partial\bm{u}_{\Delta }^{[2\xi]}(\bm{\alpha},t)}{\partial t}\Big\|_{L^\infty(0,\mathcal T;(L^2(\Omega))^3)}+
\big\| \bm{u}_{\Delta }^{[2\xi]}(\bm{\alpha},t)\big\|_{L^\infty(0,\mathcal T;(H^1_0(\Omega))^3)}\\
&+\big\| T_{\Delta}^{[2\xi]}(\bm{\alpha},t)\big\|_{L^\infty(0,\mathcal T;L^2(\Omega))}+
\big\| T_{\Delta}^{[2\xi]}(\bm{\alpha},t)\big\|_{L^2(0,\mathcal T;H^1_0(\Omega))}\leq C(\Omega,\mathcal T)\xi,
\end{aligned}}
\end{equation}
in which $C(\Omega,\mathcal T)$ is a constant irrespective of $\xi$, but dependent of $\Omega$ and $\mathcal T$.
\section{Space-time numerical algorithm}
This section presents the details of numerical algorithm to simulate the time-dependent nonlinear thermo-mechanical problem (1.1). Note that the auxiliary cell functions (2.15)-(2.17) and (2.21)-(2.32) all depend on macroscopic geometric parameter $\bm{H}$ and macroscopic temperature $T^{[0]}$, we thereby can prove the continuous properties of microscopic cell functions in Appendix A. Employing the continuous properties of auxiliary cell functions, we just need to evaluate auxiliary cell functions corresponding to distinct representative macroscopic parameters, and then employ the interpolation technique to obtain the auxiliary cell functions that are concerned in simulation process \cite{R28,R44,R48}. In the following, we present the following space-time numerical algorithm comprising of off-line and on-line stages for multi-scale nonlinear equations (1.1), as shown in Fig.\hspace{1mm}1.
\begin{figure}[!htb]
\centering
\begin{minipage}[c]{0.9\textwidth}
  \centering
  \includegraphics[width=1.0\linewidth,totalheight=2.3in]{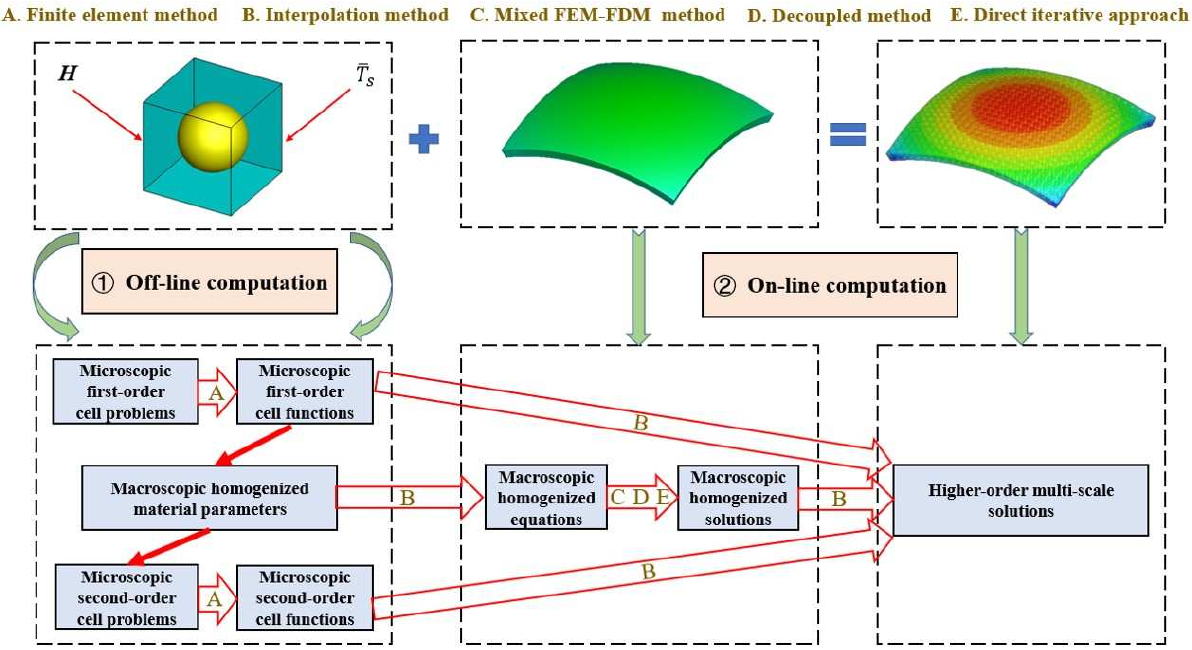}
\end{minipage}
\caption{The flowchart of space-time multi-scale numerical algorithm.}
\end{figure}
\subsection{Off-line computation for microscopic cell problems}
\begin{enumerate}
\item[(1)]
Determine the geometric configuration of PUC $\Theta=(0,1)^3$ and create tetrahedra finite element mesh family $T_{h_1}=\{K\}$ for PUC $\Theta$, where $h_1=$max$_K\{h_K\}$. Whereupon denote the linear conforming finite element space $S_{h_1}(\Theta)=\{\nu\in C^0(\bar{\Theta}):\nu\mid_{\partial \Theta}=0,\nu\mid_{K}\in P_1(K)\}\subset H^1(\Theta)$ for auxiliary cell problems.
\item[(2)]
Set computational temperature range $[T_{min},T_{max}]$ and choose a certain number of representative macroscopic parameters $\bar T_{s}$ in computational temperature range. Next, define computational domain $\Omega$ and choose a certain number of representative material points $\bm{\alpha}_I$ in computational domain. Then, evaluate the first-order cell functions defined by (2.15)-(2.17) on $S_{h_1}(\Theta)$ corresponding to distinct representative macro-scale parameters $\big(\bm{H}(\bm{\alpha}_I),\bar T_{s}\big)$. Note that classical periodic boundary condition of auxiliary cell problems is displaced by novel homogeneous Dirichlet boundary condition for practical numerical implementation \cite{R4,R30,R56,R57}. The concrete FEM scheme for solving auxiliary cell functions defined by (2.17) is proposed as below
\begin{equation}
-\int_{\Theta}{k_{ij}^{[0]}}{{\widetilde \Psi }_j}({M^m}){\widetilde \Psi }_i( \upsilon^{h_1})d\Theta=\int_{\Theta} k_{im}^{[0]}{\widetilde \Psi }_i( \upsilon^{h_1})d\Theta,\;\forall\upsilon^{h_1}\in S_{h_1}(\Theta).
\end{equation}
Then FEM is employed to solve other first-order cell functions analogously.
\item[(3)]
The macroscopic homogenized material parameters $\widehat{\rho}(T^{[0]})$, $\widehat{a}_{ijkl}(\bm{H},T^{[0]})$, $\widehat{b} _{ij}(\bm{H},T^{[0]})$, $\widehat{S}(\bm{H},T^{[0]})$, $\widehat{k}_{ij}(\bm{H},T^{[0]})$ and $\widehat\vartheta_{ij}(\bm{H},T^{[0]})$ are evaluated by formula (2.19) corresponding to distinct macroscopic parameters $\big(\bm{H}(\bm{\alpha}_I),\bar T_{s}\big)$.
\item[(4)]
Employing the same mesh as first-order auxiliary cell functions and the similar FEM scheme to (4.1), second-order auxiliary cell functions defined by (2.21)-(2.32), that correspond to distinct representative parameters $\big(\bm{H}(\bm{\alpha}_I),\bar T_{s}\big)$ in macro-scale, are solved respectively on $S_{h_1}(\Theta)$.
\end{enumerate}
\subsection{On-line computation for macroscopic homogenized problem}
\begin{enumerate}
\item[(1)]
Let $T_{h_0}=\{e\}$ be a tetrahedra finite element mesh family of the macroscopic region $\Omega$, where $h_0=$max$_e\{h_e\}$. Then define the linear conforming finite element spaces $S_{h_0}(\Omega)=\{\nu\in C^0(\bar{\Omega}):\nu\mid_{\partial\Omega_{T}}=0,\nu\mid_{e}\in P_1(e)\}\subset H^1(\Omega)$ and $S_{h_0}^{*}(\Omega)=\{\nu\in C^0(\bar{\Omega}):\nu\mid_{\partial\Omega_{u}}=0,\nu\mid_{e}\in P_1(e)\}\subset H^1(\Omega)$ for the macroscopic homogenized equations (2.18), the homogenized material parameters can be calculated by interpolation approach on each nodes $\bm{\alpha}$ of $S_{h_0}(\Omega)$ and $S_{h_0}^{*}(\Omega)$.
\item[(2)]
Solve the macroscopic homogenized equations (2.18) by mixed FDM-FEM on a coarse mesh and with a large time step on the computational domain $\Omega\times(0,\mathcal T)$. Using the equidistant time step $\displaystyle\Delta t={\mathcal T}/{M}$ to discretize time-domain $(0,\mathcal T)$ as $0=t_0<t_1<\cdots<t_M=\mathcal T$ and $t_N=N\Delta t(N=0,\cdots,M)$, then we define $u_{i}^{N}=u_i^{[0]}(\bm{\alpha},t_N)$ and $T^{N}=T^{[0]}(\bm{\alpha},t_N)$. For maintaining the numerical stability of space-time numerical scheme, the Newmark scheme and the $\delta$-scheme are utilized for the elastic-dynamics equations  and transient heat conduction equation in time domain, respectively \cite{R27}. The fully discrete scheme for dynamic homogenized thermo-mechanical equations (2.18) by mixed FDM-FEM are given as below
\begin{equation}
\left\{\begin{aligned}
&\int_{\Omega}\widehat \rho(T^{N+1}) \Big[\frac{u_i^{N+1}-u_i^{N}}{\omega(\Delta t)^2}-\frac{1}{\omega\Delta t}\frac{\partial u_i^{N}}{\partial t}-(\frac{1}{2\omega}-1)\frac{\partial^2 u_i^{N}}{\partial t^2}\Big]\nu_i^{h_0}Hd\Omega\\
&+\int_{\Omega}{\widehat{a}_{ijkl}(\bm{H},T^{N+1})\breve\varepsilon_{kl}^{[0]}(\bm{u}^{N+1})\breve\varepsilon _{ij}^{[0]}(\bm{\nu}^{h_0})}Hd\Omega\\
&-\int_{\Omega} {\widehat{b}_{ij}(\bm{H},T^{N+1})}{ (T^{N+1}-\widetilde T)}\breve\varepsilon _{ij}^{[0]}(\bm{\nu}^{h_0})Hd\Omega\\
&=\int_{\Omega}{ f_{i}({\bm{\alpha }},t_{N+1})}\nu_i^{h_0}Hd\Omega,\;\;\forall\bm{\nu}^{h_0}\in \big(S_{h_0}^{*}(\Omega)\big)^3,\\
&\int_{\Omega}{\widehat S(\bm{H},T^{N+1}) \frac{{T^{N+1}-T^{N} }}{{\Delta t}}}\widetilde{\varphi}^{h_0}Hd\Omega\\
&+\delta\int_{\Omega}\widehat{k}_{ij}(\bm{H},T^{N+1})\Psi_j(T^{N+1})\Psi_i(\widetilde{\varphi}^{h_0})Hd\Omega\\
&+(1-\delta)\int_{\Omega}\widehat{k}_{ij}(\bm{H},T^{N+1})\Psi_j(T^{N})\Psi_i(\widetilde{\varphi}^{h_0})Hd\Omega\\
&+\int_{\Omega}\widehat \vartheta _{ij}(\bm{H},T^{N+1})\frac{{\breve\varepsilon_{ij}^{[0]}(\bm{u}^{N+1})-\breve\varepsilon_{ij}^{[0]}}(\bm{u}^{N})}{{\Delta t}}\widetilde{\varphi}^{h_0}Hd\Omega\\
&=\delta\int_{\Omega} h({\bm{\alpha }},t_{N+1})\widetilde{\varphi}^{h_0}Hd\Omega\\
&+(1-\delta)\int_{\Omega} h({\bm{\alpha }},t_{N})\widetilde{\varphi}^{h_0}Hd\Omega,\;\;\forall\widetilde{\varphi}^{h_0}\in S_{h_0}(\Omega),\\
&\bm{u}^{[0]}({\bm{\alpha }},t)=\widehat{\bm{u}}({\bm{\alpha }},t),\;\;\text{on}\;\;\partial\Omega_{u},\\
&T^{[0]}({\bm{\alpha }},t)=\widehat{T}({\bm{\alpha }},t),\;\;\text{on}\;\;\partial\Omega_{T}.
\end{aligned}\right.
\end{equation}
\item[(3)]
Note that the system (4.2) still is a strongly coupled system, we here propose a decoupling approach to decompose the coupled system (4.2) into two sub-problems, and maintain the original stability and accuracy \cite{R28}. We set that
\begin{equation}
u_i^{N+1}=u_i^{N}+\varpi(u_i^{N}-u_i^{N-1}),
\end{equation}
where $\varpi$ denotes a corrected coefficient. Thereupon, two sub-systems are obtained by employing (4.3) to decompose the coupled system (4.2) as below
\begin{equation}
\left\{\begin{aligned}
&\int_{\Omega}{\widehat S(\bm{H},T^{N+1}) \frac{{T^{N+1}-T^{N} }}{{\Delta t}}}\widetilde{\varphi}^{h_0}Hd\Omega\\
&+\delta\int_{\Omega}\widehat{k}_{ij}(\bm{H},T^{N+1})\Psi_j(T^{N+1})\Psi_i(\widetilde{\varphi}^{h_0})Hd\Omega\\
&+(1-\delta)\int_{\Omega}\widehat{k}_{ij}(\bm{H},T^{N+1})\Psi_j(T^{N})\Psi_i(\widetilde{\varphi}^{h_0})Hd\Omega\\
&+\varpi\int_{\Omega}\widehat \vartheta _{ij}(\bm{H},T^{N+1})\frac{{\breve\varepsilon_{ij}^{[0]}(\bm{u}^{N})-\breve\varepsilon_{ij}^{[0]}}(\bm{u}^{N-1})}{{\Delta t}}\widetilde{\varphi}^{h_0}Hd\Omega\\
&=\delta\int_{\Omega} h({\bm{\alpha }},t_{N+1})\widetilde{\varphi}^{h_0}Hd\Omega\\
&+(1-\delta)\int_{\Omega} h({\bm{\alpha }},t_{N})\widetilde{\varphi}^{h_0}Hd\Omega,\;\;\forall\widetilde{\varphi}^{h_0}\in S_{h_0}(\Omega),\\
&T^{[0]}({\bm{\alpha }},t)=\widehat{T}({\bm{\alpha }},t),\;\;\text{on}\;\;\partial\Omega_{T}.
\end{aligned}\right.
\end{equation}
\begin{equation}
\left\{\begin{aligned}
&\int_{\Omega}\widehat \rho(T^{N+1}) \Big[\frac{u_i^{N+1}-u_i^{N}}{\omega(\Delta t)^2}-\frac{1}{\omega\Delta t}\frac{\partial u_i^{N}}{\partial t}-(\frac{1}{2\omega}-1)\frac{\partial^2 u_i^{N}}{\partial t^2}\Big]\nu_i^{h_0}Hd\Omega\\
&+\int_{\Omega}{\widehat{a}_{ijkl}(\bm{H},T^{N+1})\breve\varepsilon_{kl}^{[0]}(\bm{u}^{N+1})\breve\varepsilon _{ij}^{[0]}(\bm{\nu}^{h_0})}Hd\Omega\\
&-\int_{\Omega} {\widehat{b}_{ij}(\bm{H},T^{N+1})}{ (T^{N+1}-\widetilde T)}\breve\varepsilon _{ij}^{[0]}(\bm{\nu}^{h_0})Hd\Omega\\
&=\int_{\Omega}{ f_{i}({\bm{\alpha }},t_{N+1})}\nu_i^{h_0}Hd\Omega,\;\;\forall\bm{\nu}^{h_0}\in \big(S_{h_0}^{*}(\Omega)\big)^3,\\
&\bm{u}^{[0]}({\bm{\alpha }},t)=\widehat{\bm{u}}({\bm{\alpha }},t),\;\;\text{on}\;\;\partial\Omega_{u}.
\end{aligned}\right.
\end{equation}
As a result, we can solve the macroscopic temperature and displacement fields at each time step via (4.4) and (4.5) by turn. As is known to us, the $\delta$-scheme and the Newmark scheme can be unconditionally stable by setting their parameters $0\leq\delta\leq1$, $\gamma\geq0.5$ and $\omega\geq0.25(\gamma+0.5)^2$. In the present work, we presume $\delta=1$, $\gamma = 0.6$ and $\omega = 0.35$.
\item[(4)]
It should be highlighted that, decoupled equations (4.5) and (4.6) are a nonlinear system which can not be computed directly. Herein, the direct iterative approach is utilized to solve this nonlinear system.
\end{enumerate}
\subsection{On-line computation for higher-order multi-scale solutions}
\begin{enumerate}
\item[(1)]
For arbitrary point $(\bm{\alpha },t)=(\alpha_1,\alpha_2,\alpha_3,t)\in \Omega\times(0,\mathcal T)$, we employ the interpolation technique to solve the corresponding values of first-order cell functions, second-order cell functions and homogenized solutions.
\item[(2)]
In formulas (2.35) and (2.36), the average technique on relative elements \cite{R59} is employed to evaluate the spatial derivatives $\breve\varepsilon^{[0]}_{mn}$, $\Psi_j(\breve\varepsilon^{[0]}_{mn})$, $\Psi_{m}(T^{[0]})$ and $\Psi_{m}\Psi_{n}(T^{[0]})$, and the difference schemes are employed to evaluate the temporal derivatives $\displaystyle\frac{{{\partial ^2}u_j^{[0]}}}{{\partial {t^2}}}$, $\displaystyle\frac{{{\partial}T^{[0]}}}{{\partial {t}}}$ and $\displaystyle\frac{{\partial\breve\varepsilon _{mn}^{[0]}}}{{\partial t}}$ at every time steps.
\item[(3)]
Finally, the displacement field $\bm{u}^{[2\xi]}(\bm{\alpha},t)$ and temperature field $T^{[2\xi]}(\bm{\alpha},t)$ can be evaluated by the formulas (2.35) and (2.36) separately. Furthermore, we shall employ the higher-order interpolation and post-processing techniques to obtain the high-accuracy HOMS solutions.
\end{enumerate}
\begin{rmk}
In the present work, the equidistant interpolating points are employed as representative points. Furthermore, we recommend that at least one representative material point should be arranged at inclusion and matrix components of heterogeneous shells in practical engineering computation.
\end{rmk}
\section{Numerical results and discussions}
In this section, four numerical examples are presented to validate the feasibility of the proposed HOMS computational model and space-time numerical algorithm. The numerical experiments are implemented on a HP desktop workstation, which has a Intel(R) Core(TM) i7-8750H processor (2.20 GHz) and internal memory (16.0 GB). Moreover, all numerical simulations are implemented on the basis of experimental material parameters \cite{R5,R15,R16}.
\subsection{Multi-scale nonlinear simulation of heterogeneous block}
This example researches the nonlinear thermo-mechanical behaviors of a heterogeneous block, that is comprised of two distinct materials with matrix SiC and inclusion C. The Lam$\rm{\acute{e}}$ coefficients of this heterogeneous block are $H_1=1$, $H_2=1$ and $H_3=1$. Also, the whole domain $\Omega$ and PUC $\Theta$ are illustrated in Fig.\hspace{1mm}2, where $\Omega=(\alpha_1,\alpha_2,\alpha_3)=(x,y,z)=[0,1]\times[0,1]\times[0,1]cm^3$ and $\displaystyle\xi=1/5$.
\begin{figure}[!htb]
\centering
\begin{minipage}[c]{0.46\textwidth}
  \centering
  \includegraphics[width=0.7\linewidth,totalheight=1.5in]{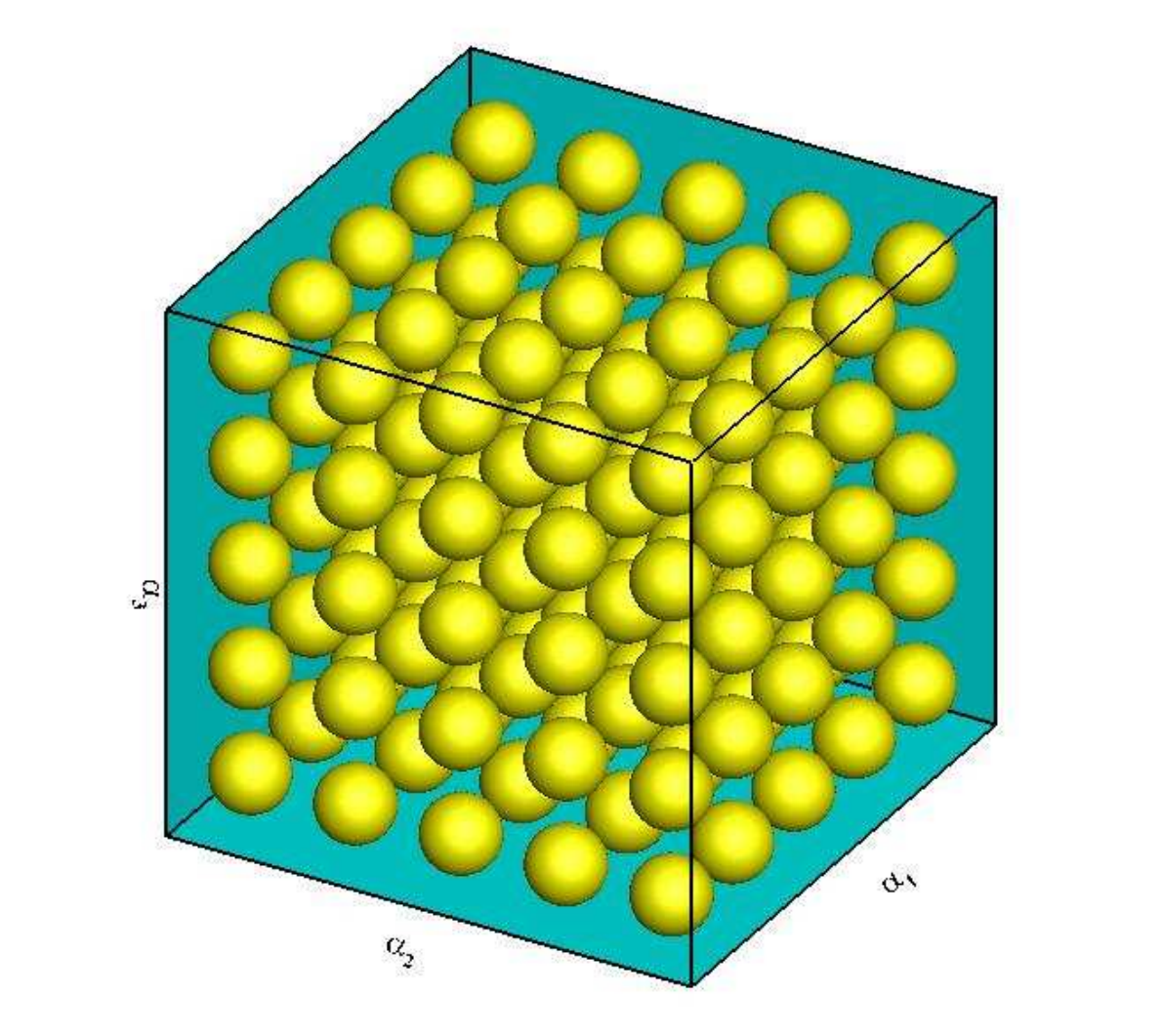} \\
  (a)
\end{minipage}
\begin{minipage}[c]{0.46\textwidth}
  \centering
  \includegraphics[width=0.7\linewidth,totalheight=1.5in]{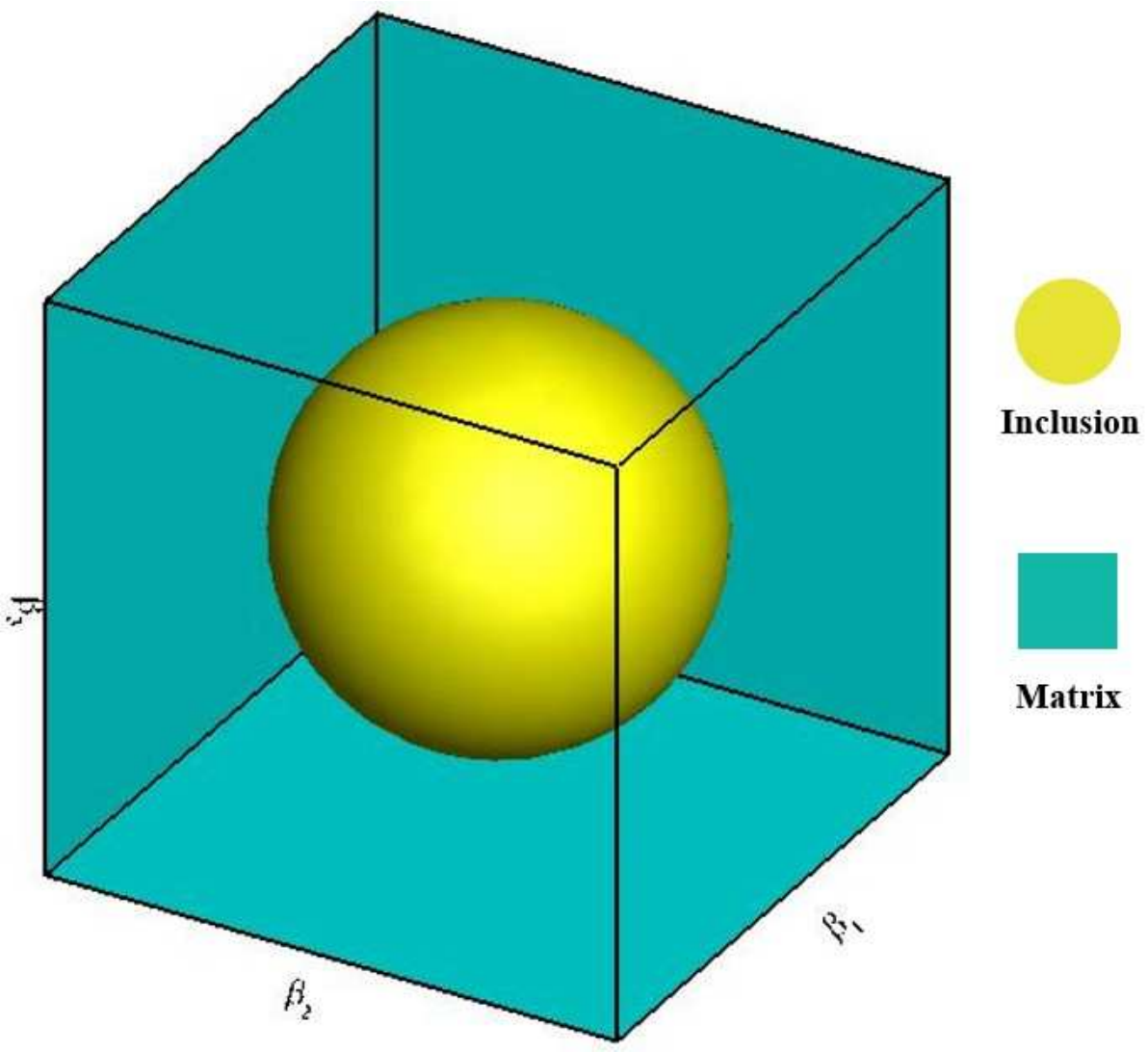} \\
  (b)
\end{minipage}
\caption{(a) The whole domain $\Omega$ of heterogeneous block; (b) the unit cell $\Theta$ (radius $r=0.35$).}
\end{figure}

The heterogeneous block $\Omega$ is clamped on its four side surfaces, that are perpendicular to $\alpha_1$-axis and $\alpha_2$-axis. The initial temperature is set as $373.15K$ and $773.15K$ at the bottom and top surfaces, respectively. Additionally, this heterogeneous block is attached by heat source $h = 5000J/(cm^3\bm{\cdot}s)$ and body forces $(f_1,f_2,f_3) = (0, 0,-10000)N/cm^3$. Furthermore, the material parameters of SiC/C heterogeneous block are exhibited in Table 1.
\begin{table}[h]{\caption{Material parameters of heterogeneous material SiC/C.}}
\centering
\begin{tabular}{cc}
\hline
Material parameters & Matrix SiC \\
\hline
Mass density $(kg/m^3)$ & 3210  \\
Young's modulus $(GPa)$ & 350.0-3.04$\times10^{-2}T$ \\
Poisson's ratio & 0.25   \\
Thermal modulus $(Pa/K)$ & 3.50  \\
Specific heat $(J/(kg\bm{\cdot} K))$ & 660.0+1.915$T$-1.491$\times10^{-3}T^2$  \\
Thermal conductivity $(W/(m\bm{\cdot} K))$ & 250.0+0.02728$T$ \\
Two-way modulus $(Pa/K)$ & 1306.025 \\
\hline
Material properties & Inclusion C\\
\hline
Mass density $(kg/m^3)$ & 1760 \\
Young's modulus $(GPa)$ & 220.0-1.10$\times10^{-4}T$  \\
Poisson's ratio & 0.20\\
Thermal modulus $(Pa/K)$ & 9273.0-57.53$\times T$+0.0817$\times T^2$\\
Specific heat $(J/(kg\bm{\cdot} K))$ & 710.0+$1.781T$-1.976$\times10^{-3}T^2$ \\
Thermal conductivity $(W/(m\bm{\cdot} K))$ & 8.0+0.02535$T$ \\
Two-way modulus $(Pa/K)$ & 3.46$\times10^{6}$-2.147$\times10^{4}T$+30.486$\times T^2$\\
\hline
\end{tabular}
\end{table}

Obviously, it is very difficult to obtain the exact solutions of this example, we replace $T^\xi(\bm{\alpha},t)$ and $\bm{u}^\xi(\bm{\alpha},t)$ with direct finite element solutions $T_{Fe}^\xi(\bm{\alpha},t)$ and $\bm{u}_{Fe}^\xi(\bm{\alpha},t)$ on a very fine mesh. Herein, some error models are elaborated as below
\begin{equation}
\text{Terr$_0$}=\frac{||T^{[0]}-T_{Fe}^\xi||_{L^2}}{||T_{Fe}^\xi||_{L^2}},
\text{Terr$_1$}=\frac{||T^{[1\xi]}-T_{Fe}^\xi||_{L^2}}{||T_{Fe}^\xi||_{L^2}},
\text{Terr$_2$}=\frac{||T^{[2\xi]}-T_{Fe}^\xi||_{L^2}}{||T_{Fe}^\xi||_{L^2}}.
\end{equation}
\begin{equation}
\text{TErr$_0$}=\frac{|T^{[0]}-T_{Fe}^\xi|_{H^1}}{|T_{Fe}^\xi|_{H^1}},
\text{TErr$_1$}=\frac{|T^{[1\xi]}-T_{Fe}^\xi|_{H^1}}{|T_{Fe}^\xi|_{H^1}},
\text{TErr$_2$}=\frac{|T^{[2\xi]}-T_{Fe}^\xi|_{H^1}}{|T_{Fe}^\xi|_{H^1}}.
\end{equation}
\begin{equation}
\text{Derr$_0$}=\frac{||\bm{u}^{[0]}-\bm{u}^\xi_{Fe}||_{L^2}}{||\bm{u}^\xi_{Fe}||_{L^2}},
\text{Derr$_1$}=\frac{||\bm{u}^{[1\xi]}-\bm{u}^\xi_{Fe}||_{L^2}}{||\bm{u}^\xi_{Fe}||_{L^2}},
\text{Derr$_2$}=\frac{||\bm{u}^{[2\xi]}-\bm{u}^\xi_{Fe}||_{L^2}}{||\bm{u}^\xi_{Fe}||_{L^2}}.
\end{equation}
\begin{equation}
\text{DErr$_0$}=\frac{|\bm{u}^{[0]}-\bm{u}^\xi_{Fe}|_{H^1}}{|\bm{u}^\xi_{Fe}|_{H^1}},
\text{DErr$_1$}=\frac{|\bm{u}^{[1\xi]}-\bm{u}^\xi_{Fe}|_{H^1}}{|\bm{u}^\xi_{Fe}|_{H^1}},
\text{DErr$_2$}=\frac{|\bm{u}^{[2\xi]}-\bm{u}^\xi_{Fe}|_{H^1}}{|\bm{u}^\xi_{Fe}|_{H^1}}.
\end{equation}
In foregoing formulas (5.1)-(5.4), $|\bm{u}^{[0]}-\bm{u}^\xi_{Fe}|_{H^1(\Omega)}=\Big(\sum\limits_{i,j=1}^{3}  \big|\big|\varepsilon_{ij}(\bm{u}^{[0]}-\bm{u}^\xi_{Fe})\big|\big|_{L^2(\Omega)}\Big)^{\frac{1}{2}}$.

Next, we create tetrahedra finite element mesh for multi-scale nonlinear problem (1.1), auxiliary cell problems and corresponding homogenized problem (2.18). This example needs to off-line solve 7075 times auxiliary cell problems totally, in which the quantity of first-order cell functions and second-order cell functions is respectively 33 and 250 with 25 macroscopic interpolation temperature. It is noteworthy that although auxiliary cell problems are costly calculated many times, they are all off-line computation before on-line multi-scale computation and can be employed in different heterogeneous structures. Then the specific numbers of FEM elements and nodes, and the computational times spent for direct finite element and multi-scale simulations are presented in Table 2.
\begin{table}[!htb]{\caption{Summary of computational cost.}}
\centering
\begin{tabular}{cccc}
\hline
 & Multi-scale eqs. & Cell eqs. & Homogenized eqs. \\
\hline
FEM elements & 983173 & 75466 & 93750\\
FEM nodes    & 165393 & 13062 & 17576\\
\hline
& DNS & off-line stage & on-line stage \\
\hline
Computational time & 42.52h  & 9.17h & 5.56h\\
\hline
\end{tabular}
\end{table}

After that, multi-scale nonlinear equations (1.1) and macroscopic homogenized equations (2.18) are on-line solved separately, where the temporal step is set as $\Delta t = 0.01$s. The nonlinear thermo-mechanical coupling problem of heterogeneous block structure is simulated in the time interval $t\in[0,1]$s.

Figs.\hspace{1mm}3-5 display the simulative results for solutions $T^{[0]}$, $T^{[1\xi]}$, $T^{[2\xi]}$, $T_{Fe}^\xi$, and $u_1^{[0]}$, $u_1^{[1\xi]}$, $u_1^{[2\xi]}$, $u_{1,Fe}^\xi$, and $u_3^{[0]}$, $u_3^{[1\xi]}$, $u_3^{[2\xi]}$, $u_{3,Fe}^\xi$ at the final moment $t$=1.0s, respectively.
\begin{figure}[!htb]
\centering
\begin{minipage}[c]{0.4\textwidth}
  \centering
  \includegraphics[width=50mm]{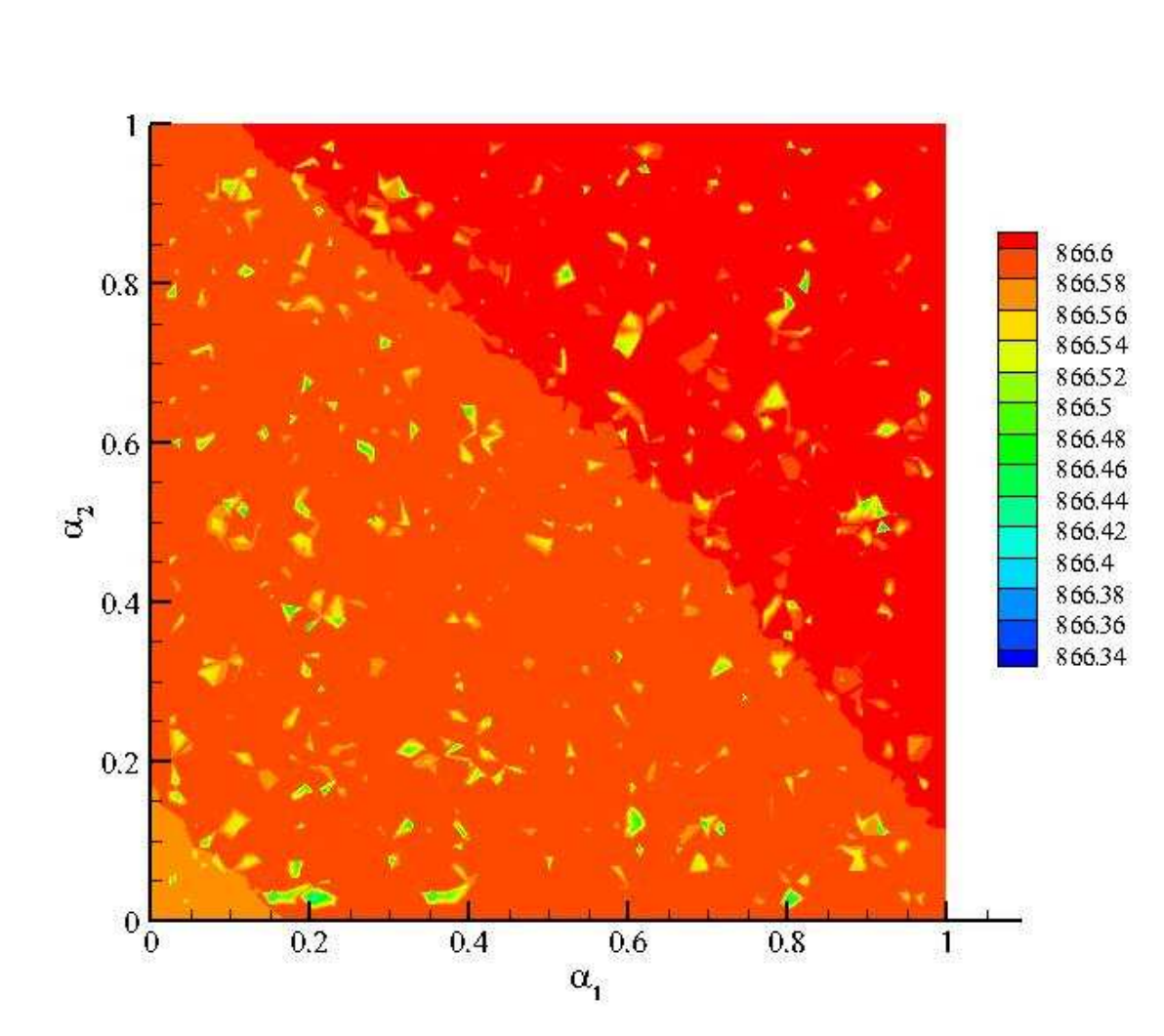}\\
  (a)
\end{minipage}
\begin{minipage}[c]{0.4\textwidth}
  \centering
  \includegraphics[width=50mm]{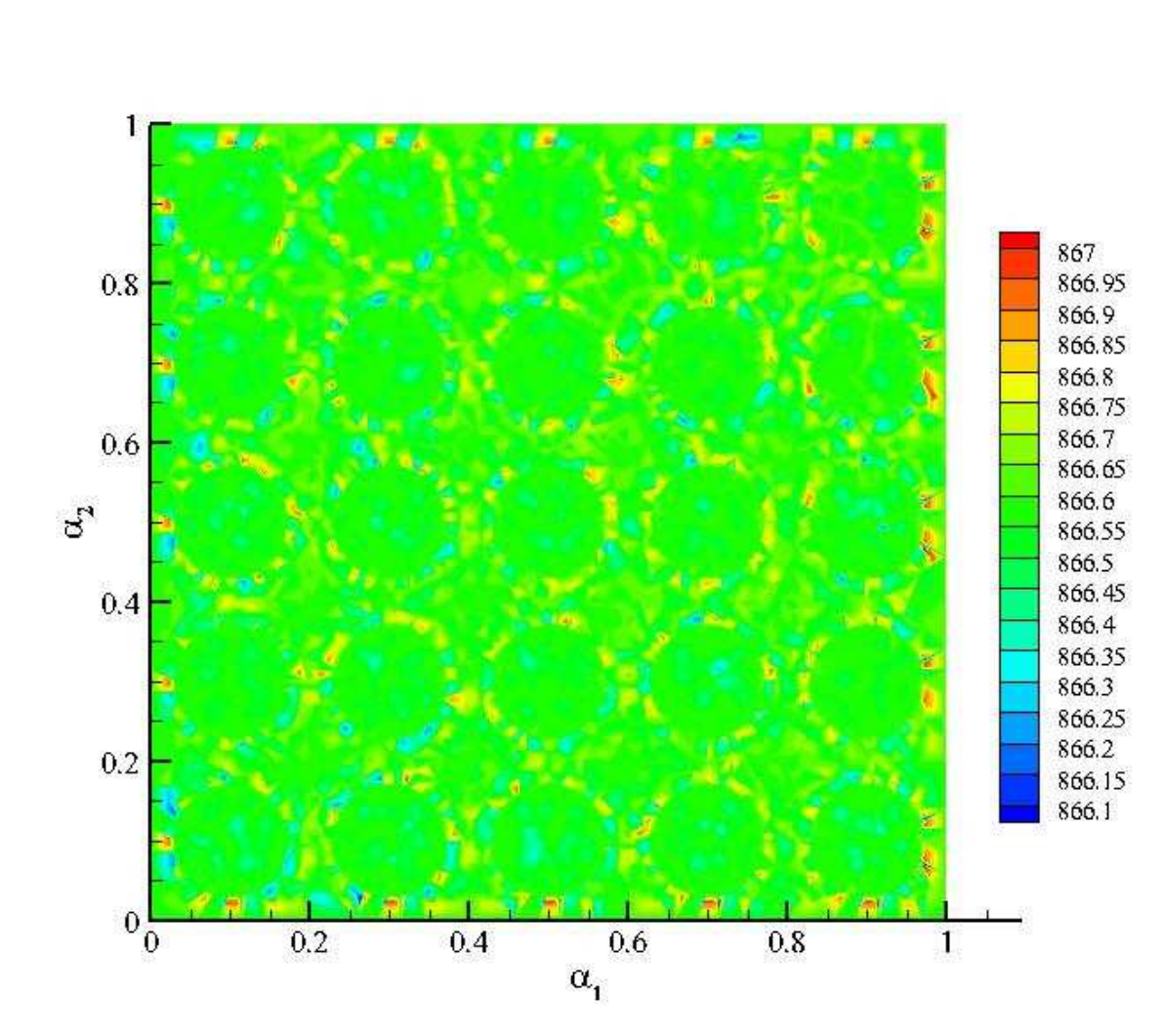}\\
  (b)
\end{minipage}
\begin{minipage}[c]{0.4\textwidth}
  \centering
  \includegraphics[width=50mm]{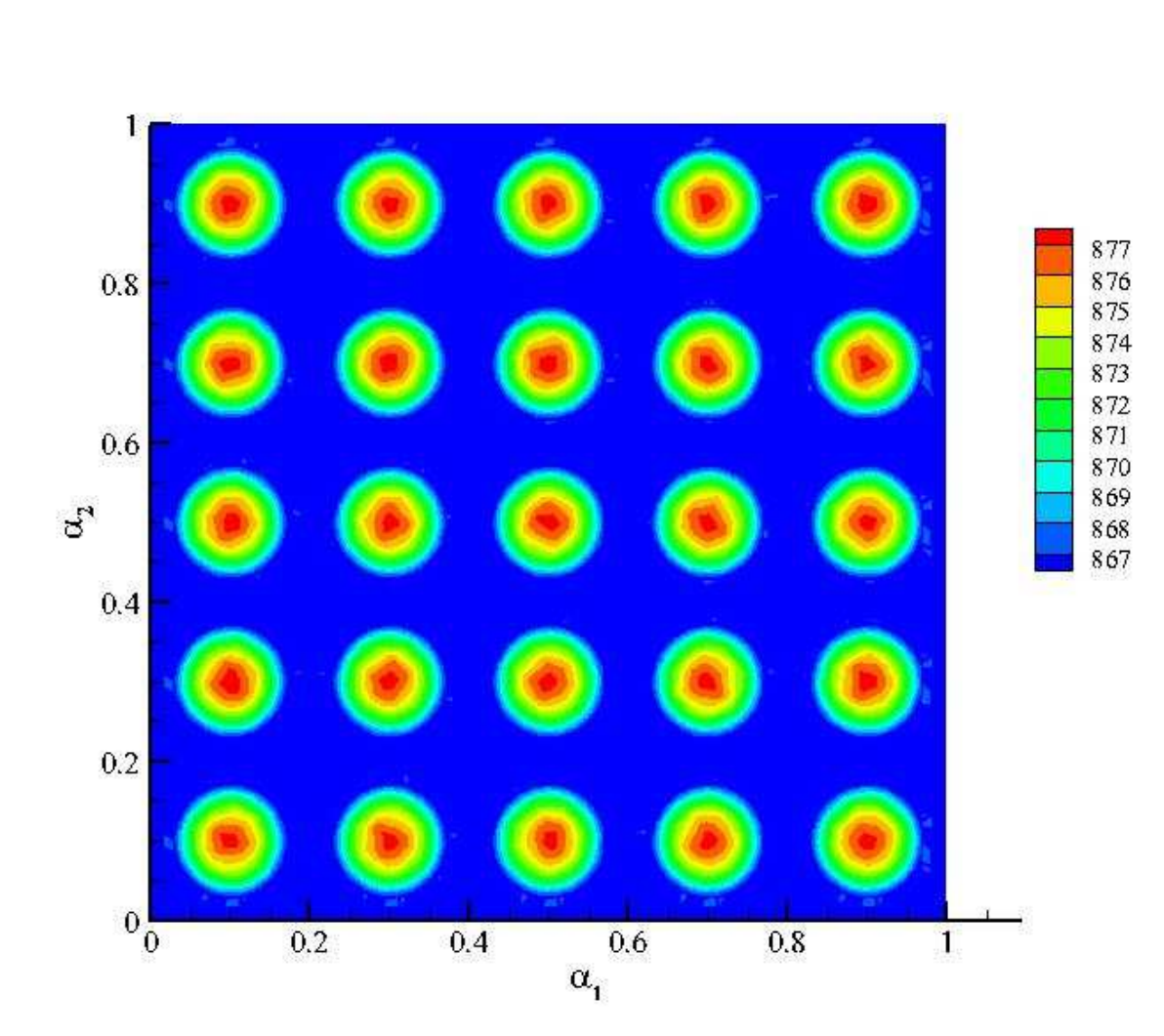}\\
  (c)
\end{minipage}
\begin{minipage}[c]{0.4\textwidth}
  \centering
  \includegraphics[width=50mm]{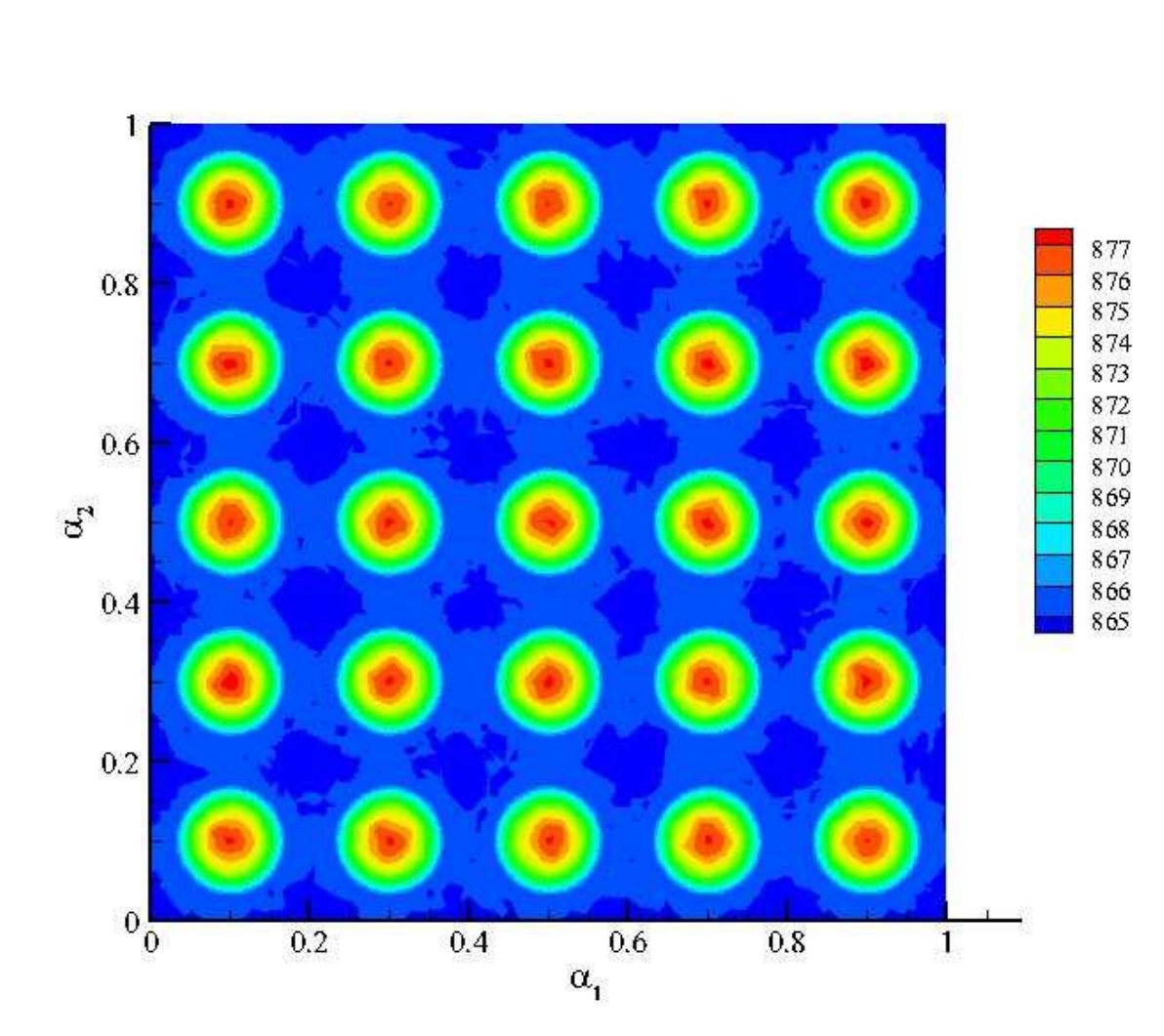}\\
  (d)
\end{minipage}
\caption{The temperature field in cross section $\alpha_3=0.5cm$ at time $t$=1.0s: (a) $T^{[0]}$; (b) $T^{[1\xi]}$; (c) $T^{[2\xi]}$; (d) $T_{Fe}^\xi$.}
\end{figure}
\begin{figure}[!htb]
\centering
\begin{minipage}[c]{0.4\textwidth}
  \centering
  \includegraphics[width=50mm]{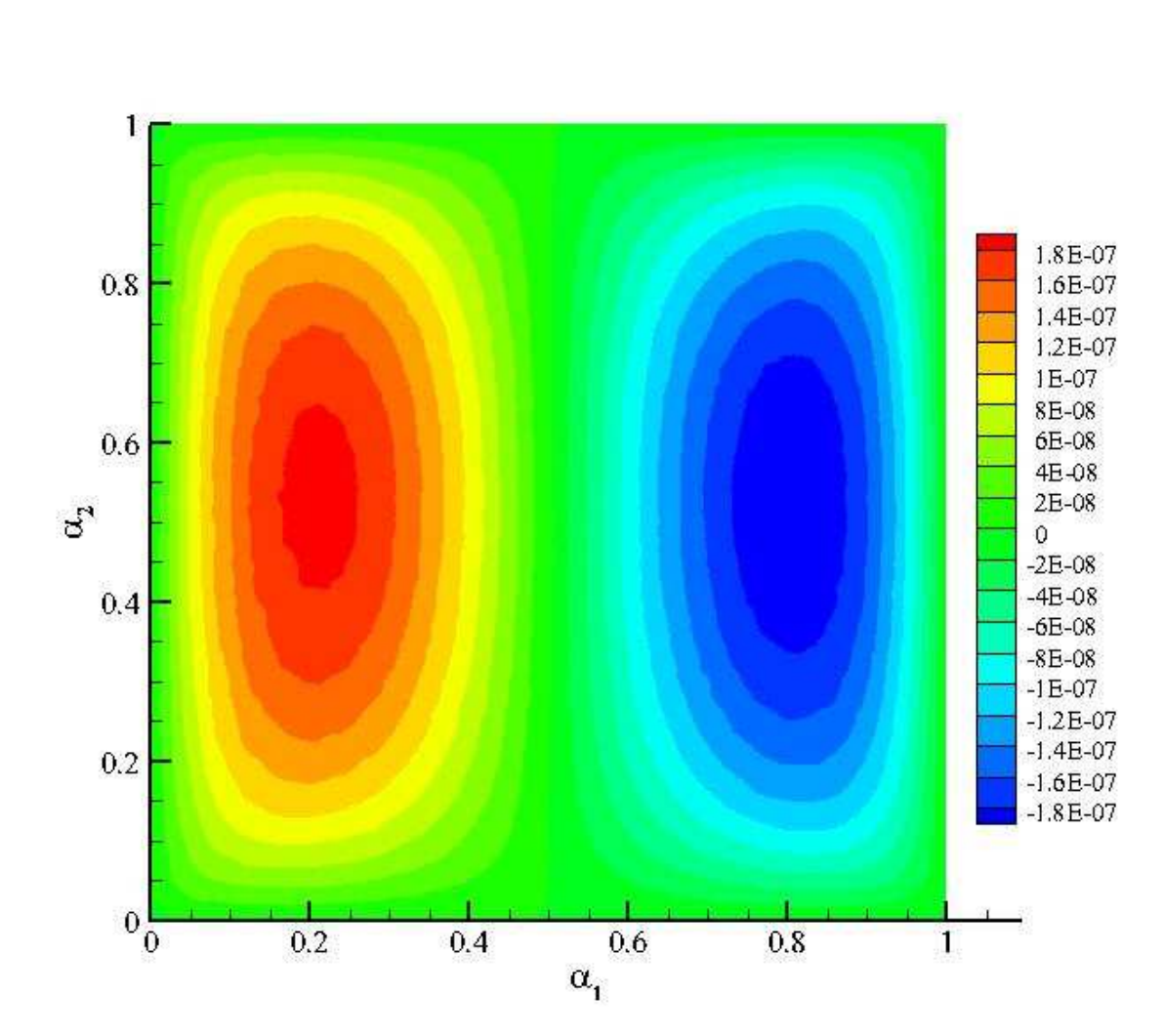}\\
  (a)
\end{minipage}
\begin{minipage}[c]{0.4\textwidth}
  \centering
  \includegraphics[width=50mm]{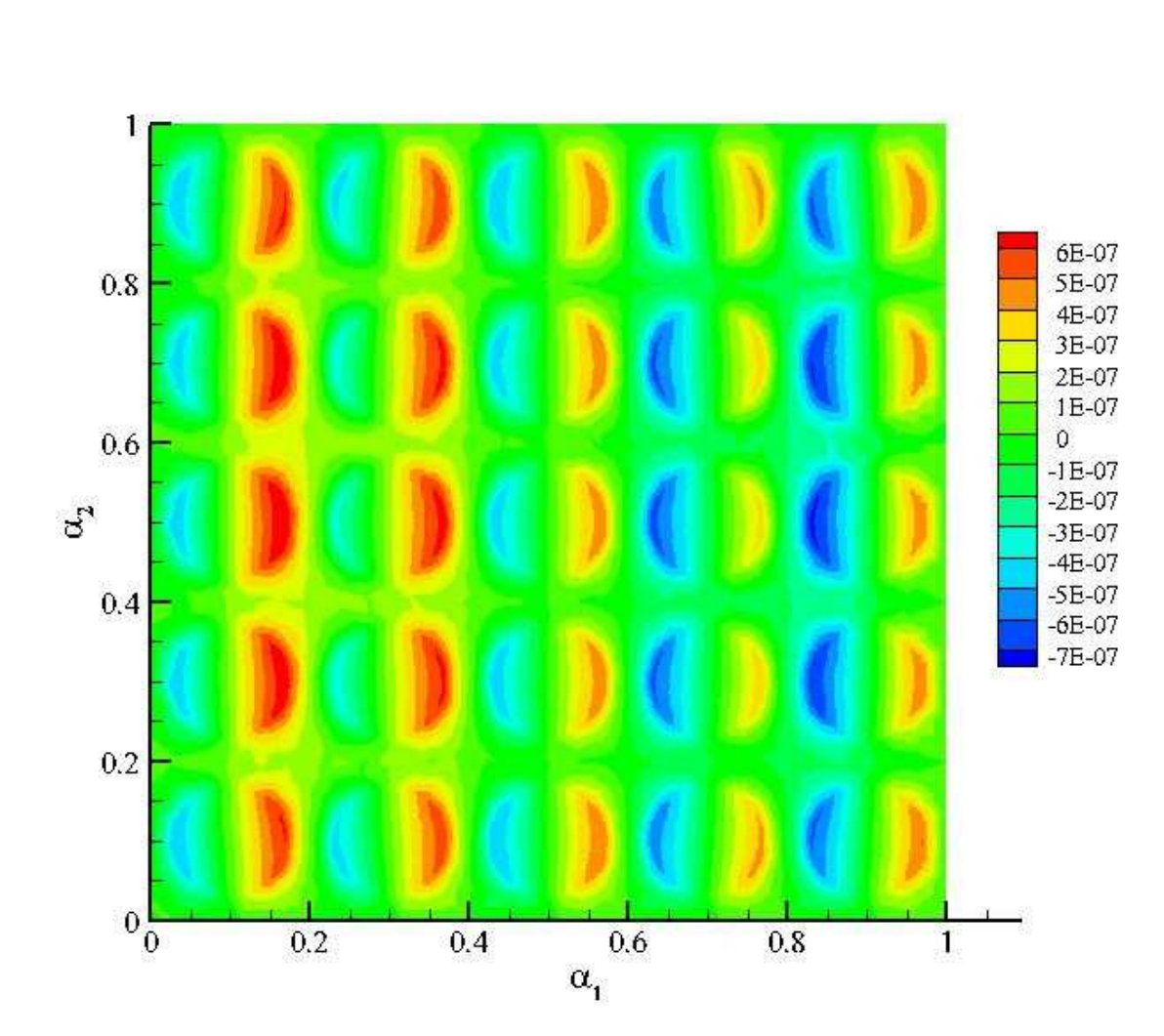}\\
  (b)
\end{minipage}
\begin{minipage}[c]{0.4\textwidth}
  \centering
  \includegraphics[width=50mm]{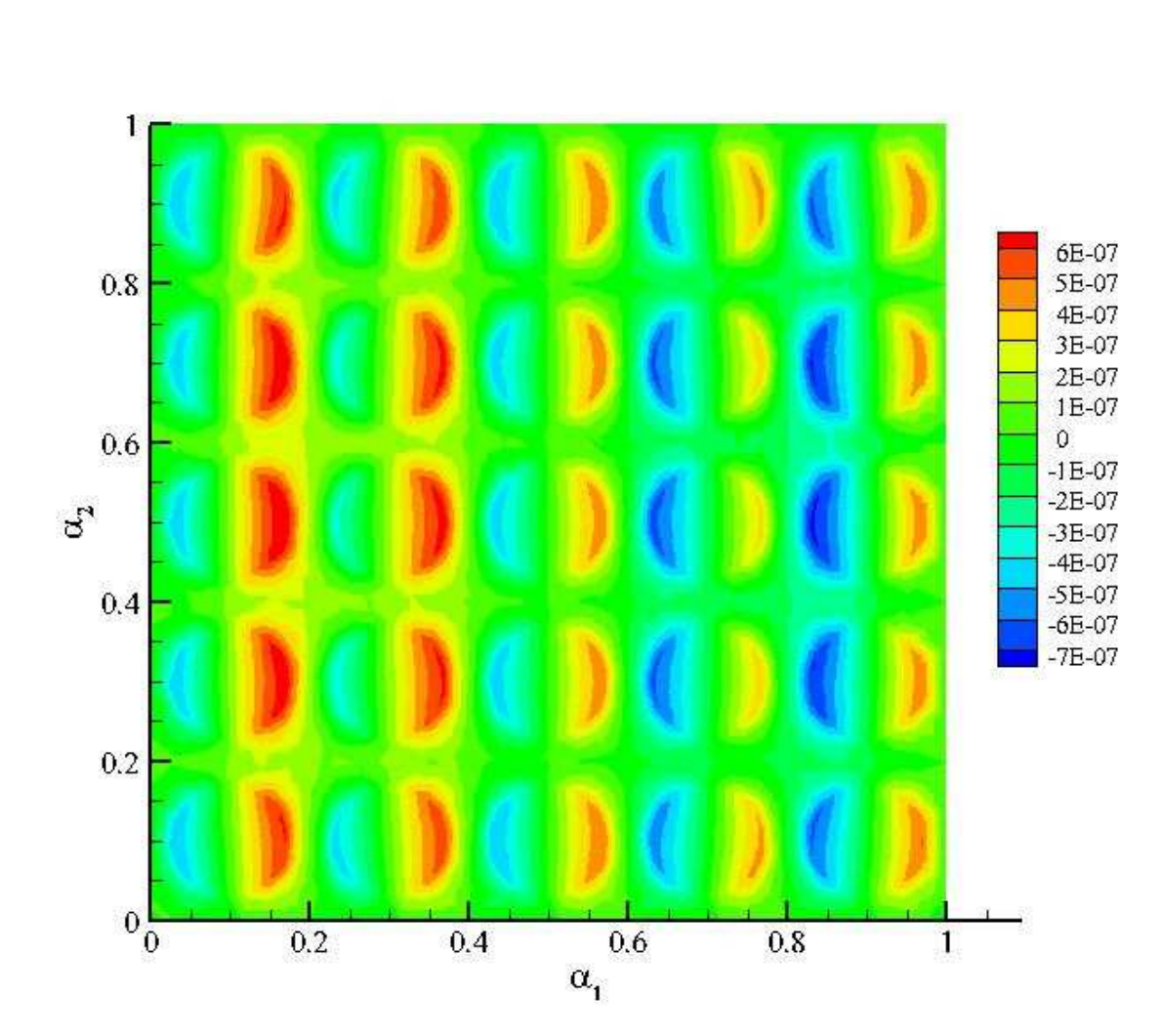}\\
  (c)
\end{minipage}
\begin{minipage}[c]{0.4\textwidth}
  \centering
  \includegraphics[width=50mm]{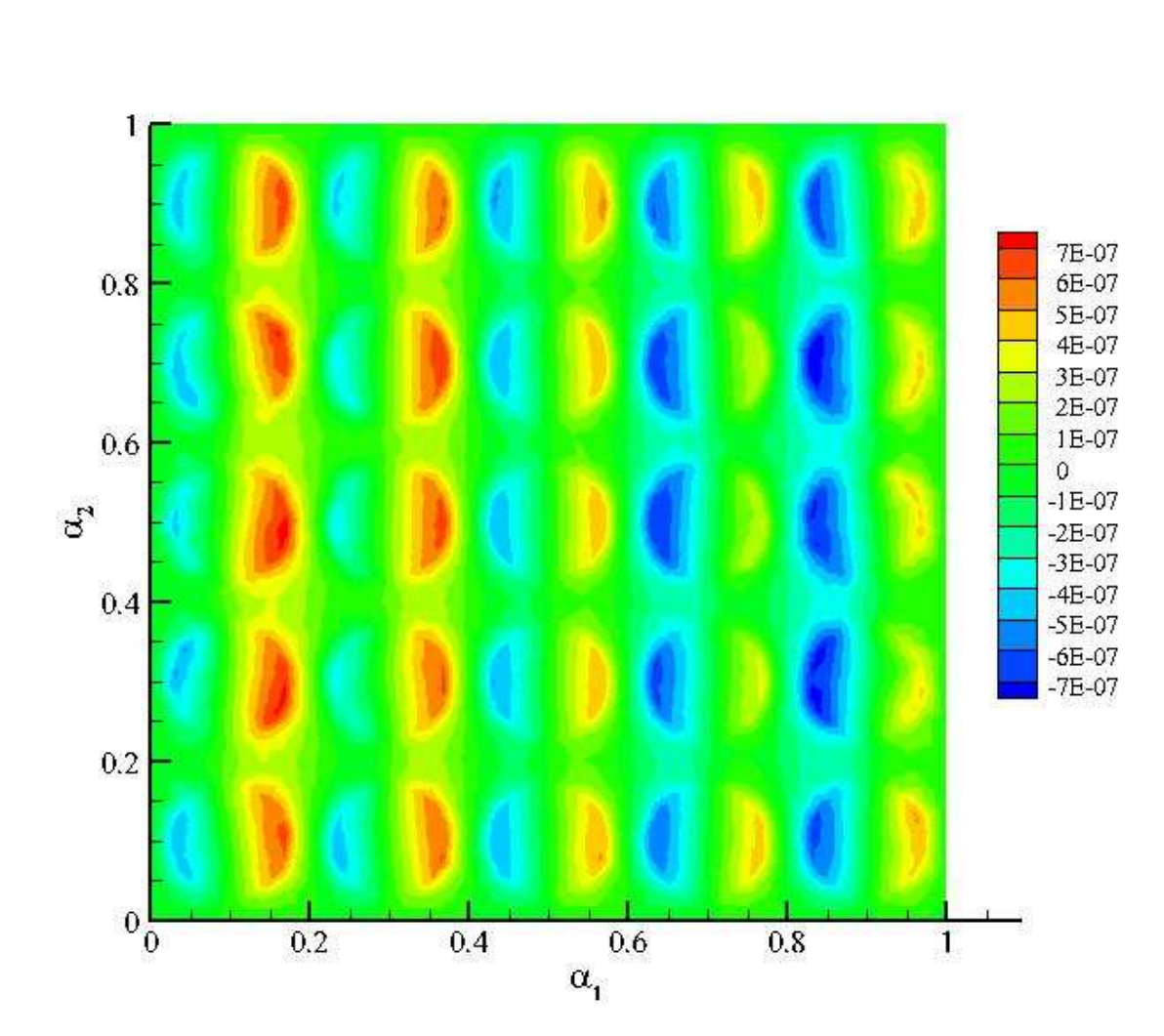}\\
  (d)
\end{minipage}
\caption{The first displacement component in cross section $\alpha_3=0.5cm$ at time $t$=1.0s: (a) $u_1^{[0]}$; (b) $u_1^{[1\xi]}$; (c) $u_1^{[2\xi]}$; (d) $u_{1,Fe}^\xi$.}
\end{figure}
\begin{figure}[!htb]
\centering
\begin{minipage}[c]{0.4\textwidth}
  \centering
  \includegraphics[width=50mm]{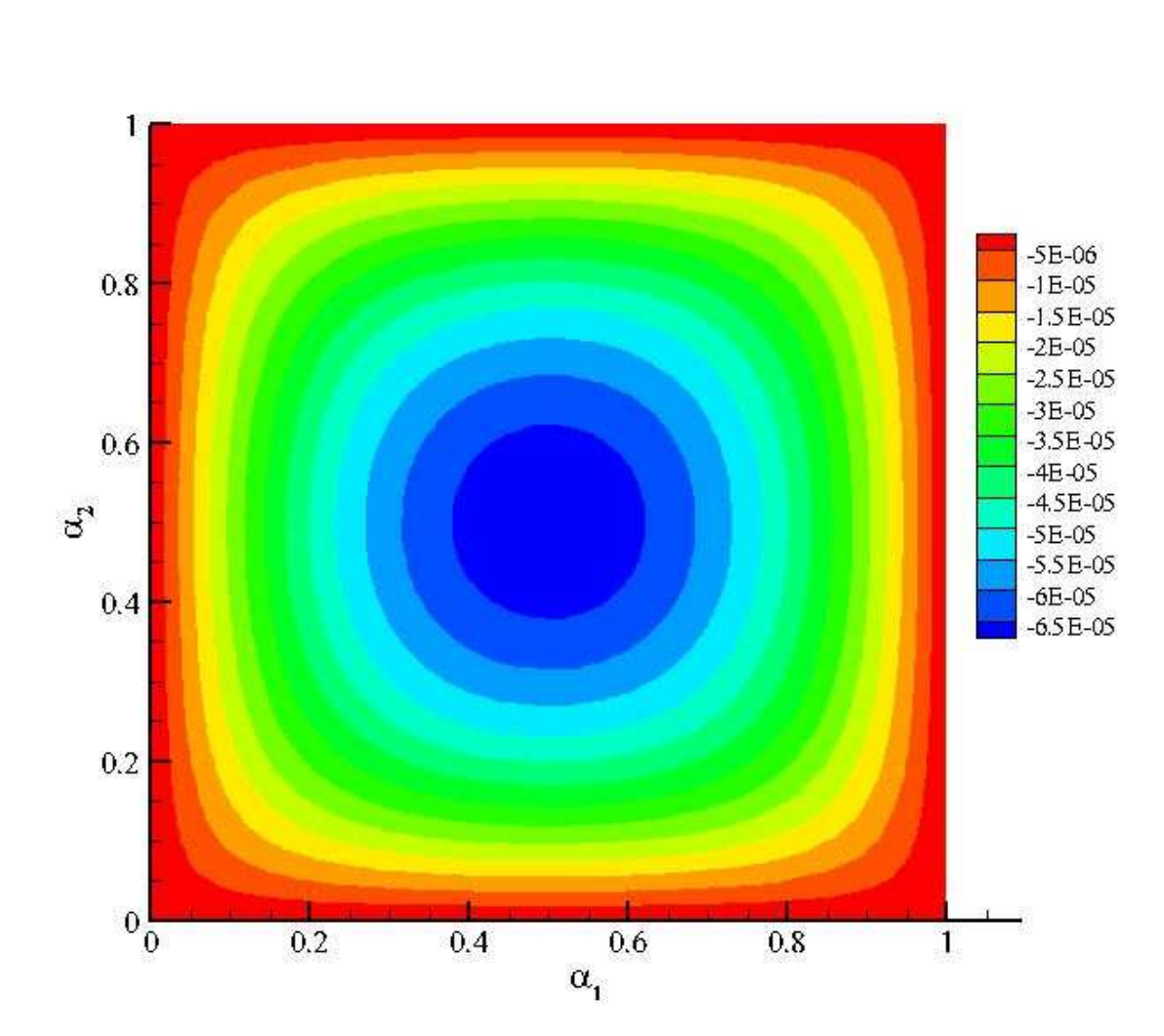}\\
  (a)
\end{minipage}
\begin{minipage}[c]{0.4\textwidth}
  \centering
  \includegraphics[width=50mm]{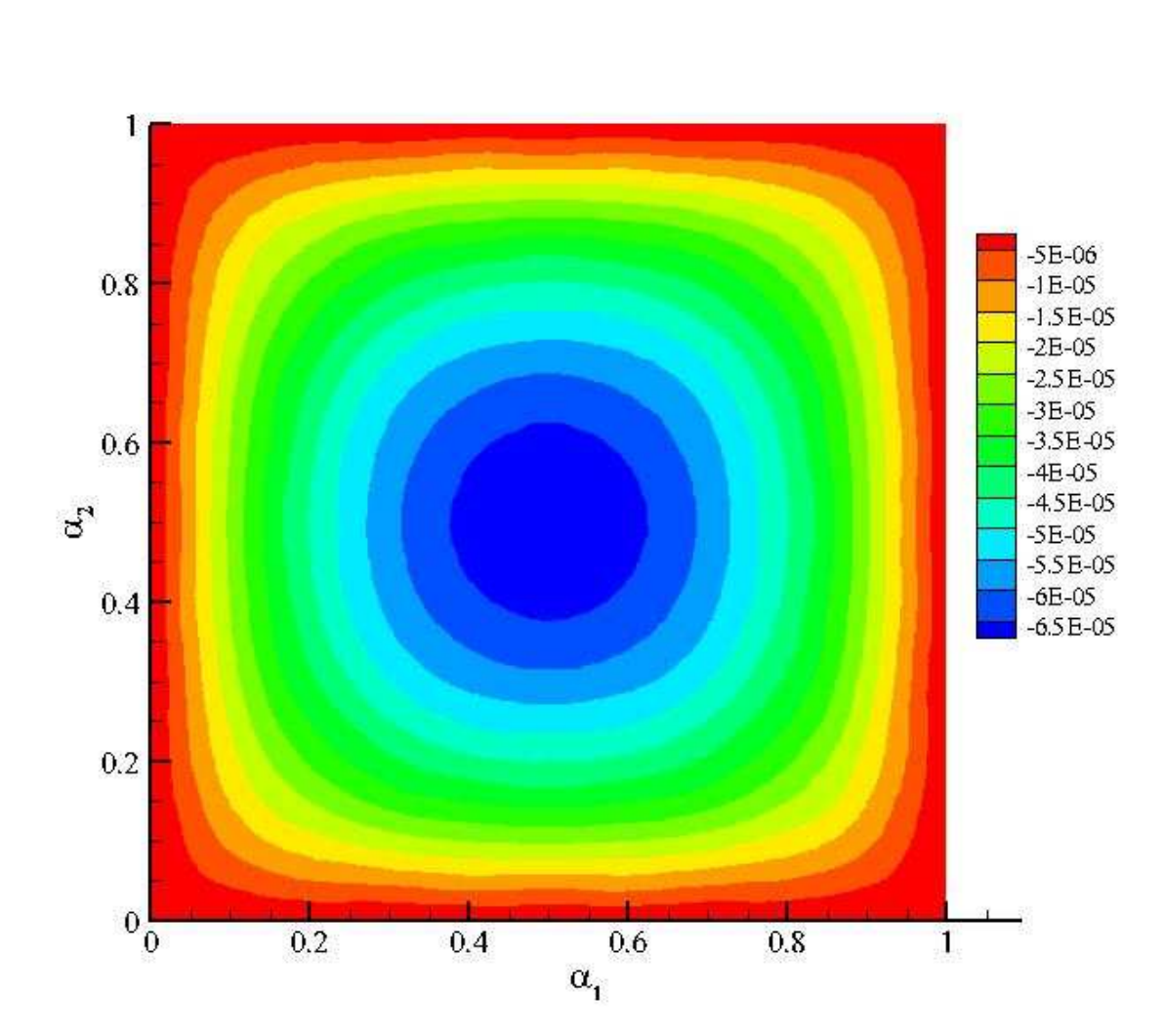}\\
  (b)
\end{minipage}
\begin{minipage}[c]{0.4\textwidth}
  \centering
  \includegraphics[width=50mm]{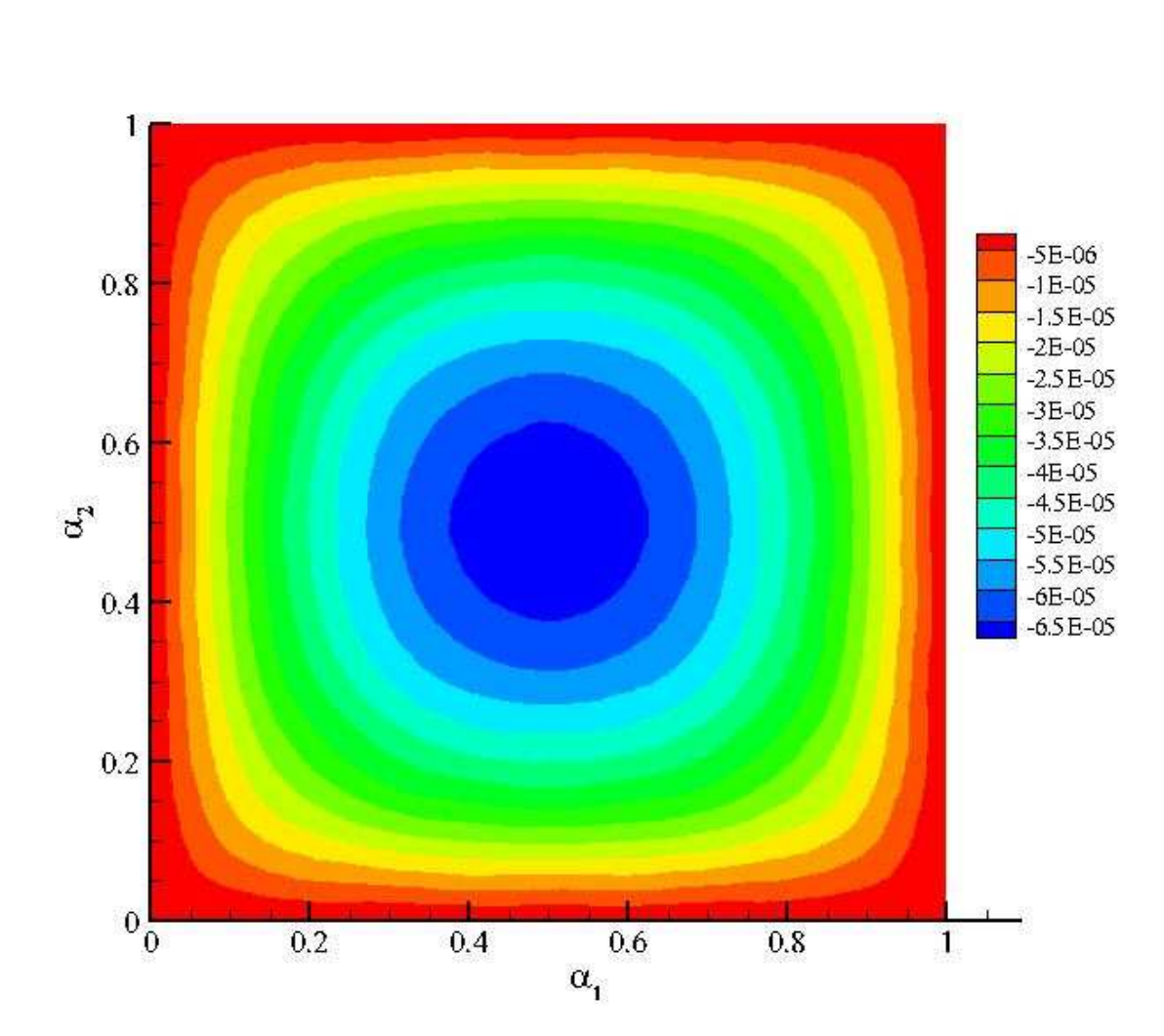}\\
  (c)
\end{minipage}
\begin{minipage}[c]{0.4\textwidth}
  \centering
  \includegraphics[width=50mm]{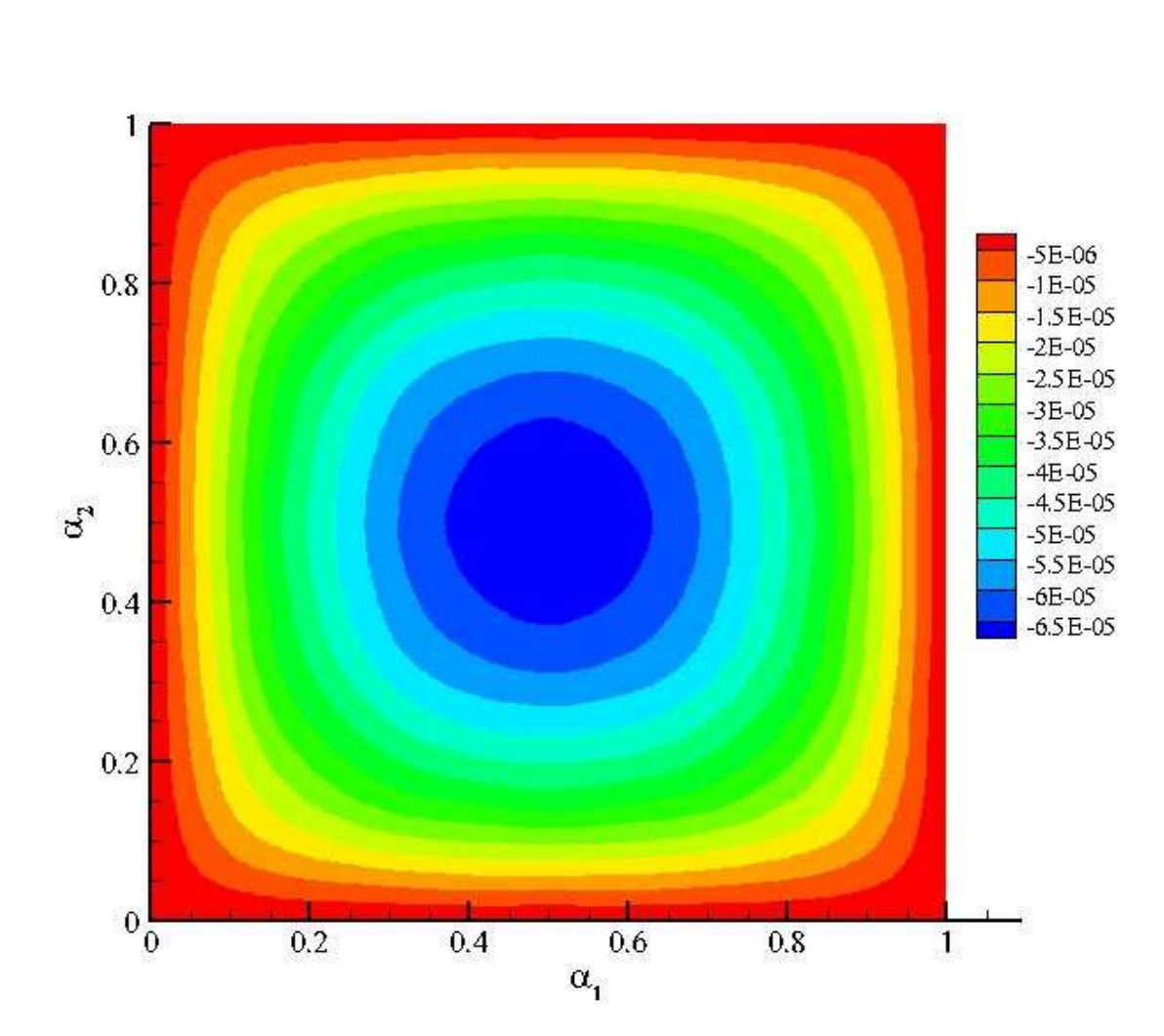}\\
  (d)
\end{minipage}
\caption{The third displacement component in cross section $\alpha_3=0.5cm$ at time $t$=1.0s: (a) $u_3^{[0]}$; (b) $u_3^{[1\xi]}$; (c) $u_3^{[2\xi]}$; (d) $u_{3,Fe}^\xi$.}
\end{figure}

Fig.\hspace{1mm}6 displays the evolutive relative errors of temperature and displacement fields.
\begin{figure}[!htb]
\centering
\begin{minipage}[c]{0.4\textwidth}
  \centering
  \includegraphics[width=50mm]{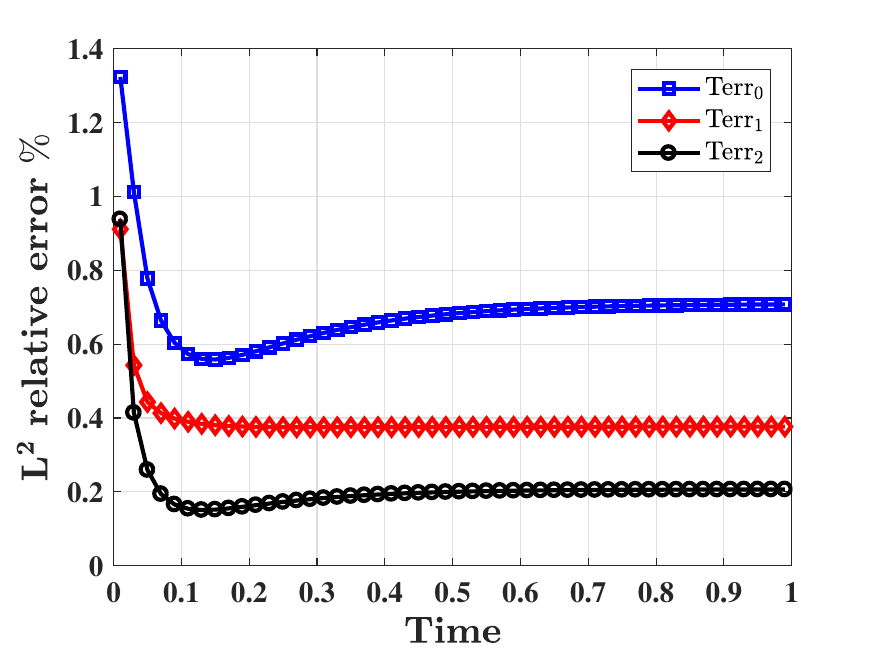}
  (a)
\end{minipage}
\begin{minipage}[c]{0.4\textwidth}
  \centering
  \includegraphics[width=50mm]{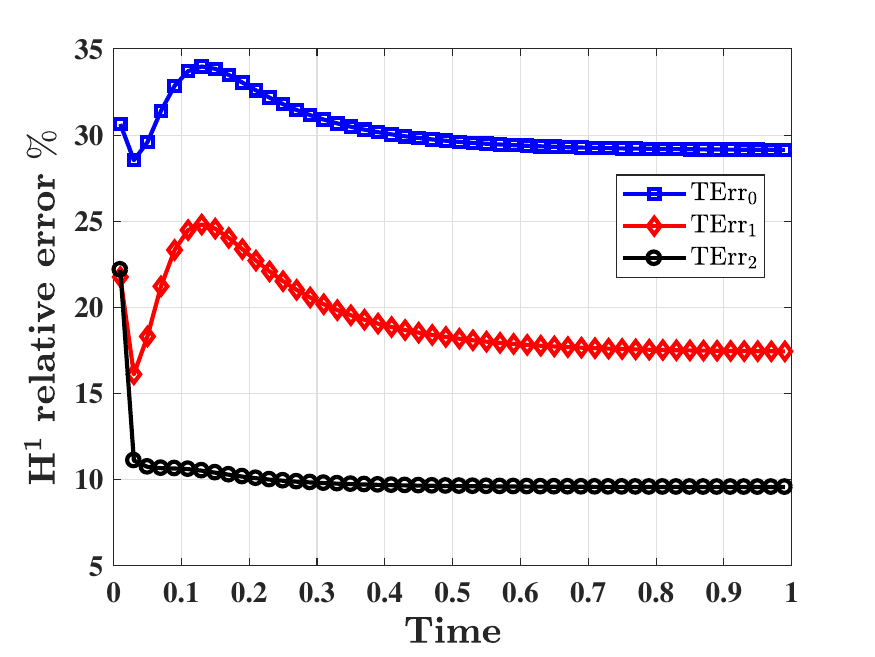}
  (b)
\end{minipage}
\begin{minipage}[c]{0.4\textwidth}
  \centering
  \includegraphics[width=50mm]{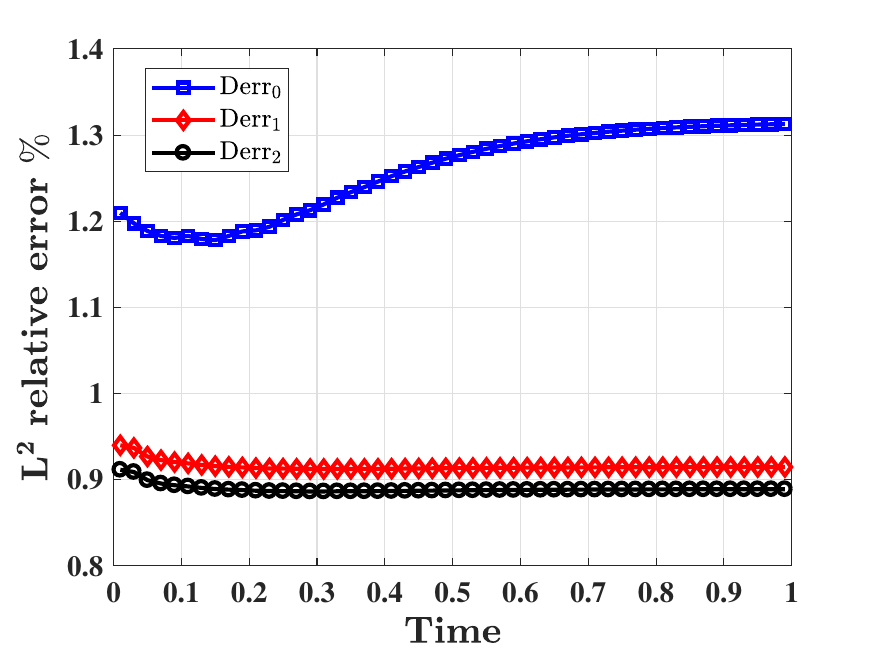}
  (c)
\end{minipage}
\begin{minipage}[c]{0.4\textwidth}
  \centering
  \includegraphics[width=50mm]{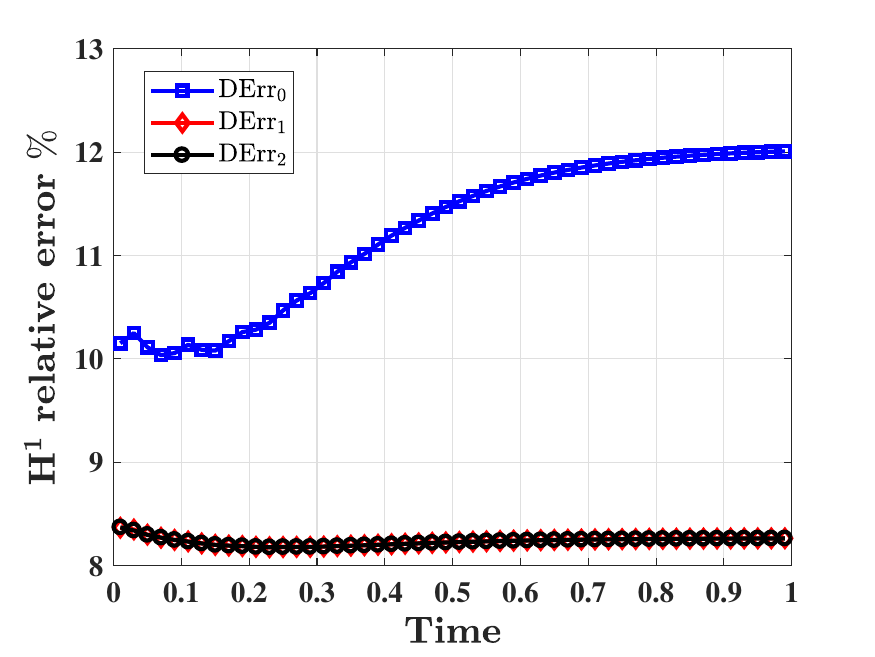}
  (d)
\end{minipage}
\caption{The evolutive relative errors of temperature and displacement fields: (a) \rm{Terr}; (b) \rm{TErr}; (c) \rm{Derr}; (d) \rm{DErr}.}
\end{figure}

As exhibited in Table 2, we can deduce that the computational resource overhead of HOMS approach is far less than direct finite element simulation both in computer space and time resources. The superiority of the proposed HOMS method over the full-scale simulation associated with direct numerical simulation is obvious. When directly simulating the multi-scale nonlinear thermo-mechanical equations, a highly fine mesh is demanded to catch the microscopic oscillatory behaviors in this heterogeneous block. After consuming enormous computing time, the convergence of direct numerical simulation is still hard to guarantee. Comparatively, although its off-line stage requires relatively more computation time, the result is a vastly accelerated simulation by the proposed HOMS approach for multi-scale nonlinear thermo-mechanical problem, which can economize $65.36\%$ computational time.

Also, Figs.\hspace{1mm}3-5 demonstrate that only HOMS solutions, which show good agreement with full-scale simulation solutions, can precisely catch the steeply oscillatory information resulting from microscopic heterogeneities in inhomogeneous block. As contrasted in Figs.\hspace{1mm}3-5, homogenized and LOMS solutions can not provide enough numerical accuracy for simulating the dynamic nonlinear thermo-mechanical problem of inhomogeneous block especially in the $H^1$ semi-norm sense. In addition, we can clearly determine that the proposed space-time multi-scale method is stable for efficient simulation on long times from Fig.\hspace{1mm}6.
\subsection{Multi-scale nonlinear simulation of heterogeneous thin plate}
This example investigates a heterogeneous thin plate with $\displaystyle\xi=1/25$, where the whole domain $\Omega=(\alpha_1,\alpha_2,\alpha_3)=(x,y,z)=[0,1]\times[0,1]\times[0,0.2]cm^3$ and PUC $\Theta$ are displayed in Fig.\hspace{1mm}7(a) and Fig.\hspace{1mm}7(b) respectively. It is a remarkable fact that the HOMS solutions of this example own the same forms as first example because they all have periodic configurations in cartesian coordinates.
\begin{figure}[!htb]
\centering
\begin{minipage}[c]{0.46\textwidth}
  \centering
  \includegraphics[width=0.7\linewidth,totalheight=1.5in]{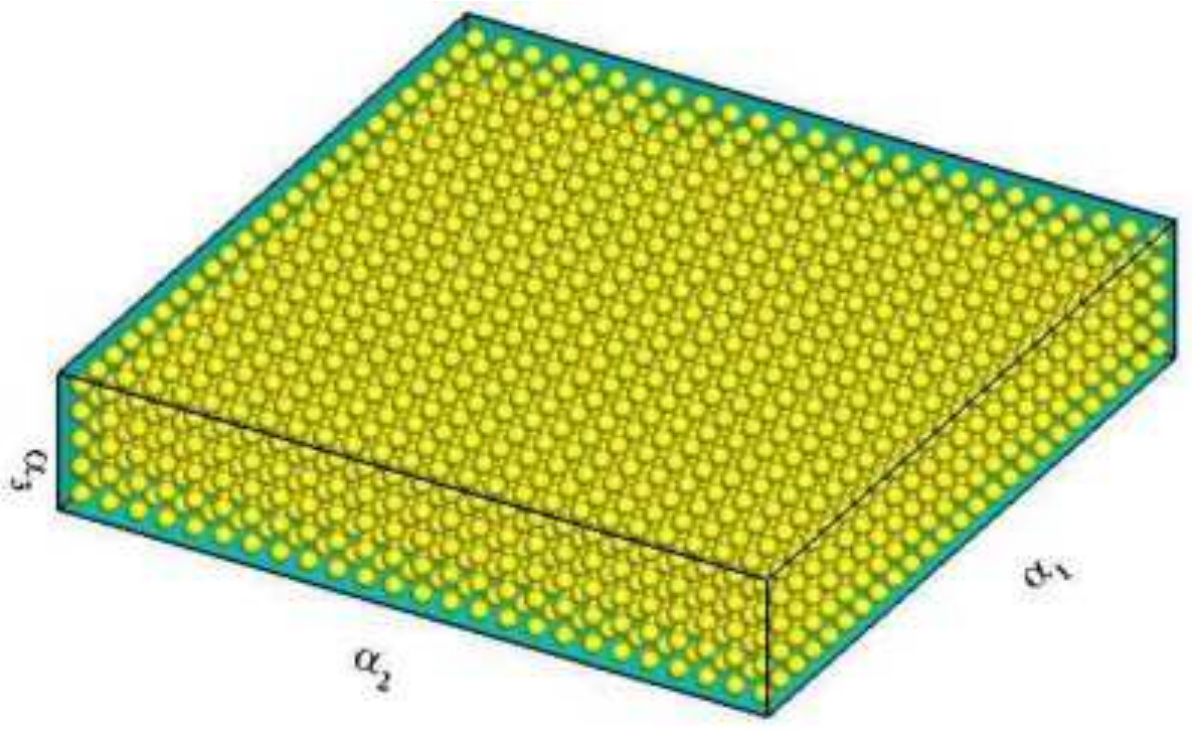}\\
  (a)
\end{minipage}
\begin{minipage}[c]{0.46\textwidth}
  \centering
  \includegraphics[width=0.7\linewidth,totalheight=1.5in]{EX4Nn-new-eps-converted-to.pdf} \\
  (b)
\end{minipage}
\caption{(a) The whole domain $\Omega$ of heterogeneous thin plate; (b) the unit cell $\Theta$.}
\end{figure}

The investigated heterogeneous thin plate in this example comprise a large number of microscopic unit cells. Direct finite element simulation needs very fine meshes and its convergent solutions can not be obtained. Thereby, we do not give the reference solutions $T_{Fe}^\xi(\bm{\alpha},t)$ and $\bm{u}_{Fe}^\xi(\bm{\alpha},t)$ for this example. Moreover, the material parameters of this heterogeneous thin plate are defined with the same values as first example. The heterogeneous thin plate $\Omega$ is clamped on its four side surfaces, that are perpendicular to $\alpha_1$-axis and $\alpha_2$-axis. The initial temperature is set as $373.15K$ and $773.15K$ at the bottom surface and the top surface, respectively. Additionally, we set $h = 10000J/(cm^3\bm{\cdot}s)$ and $(f_1,f_2,f_3) = (0,0,-20000)N/cm^3$. The computational overhead of this example is displayed in Table 3.
\begin{table}[!htb]{\caption{Summary of computational cost.}}
\centering
\begin{tabular}{cccc}
\hline
 & Multi-scale eqs. & Cell eqs. & Homogenized eqs. \\
\hline
FEM elements & 9600000 & 75466 & 150000\\
FEM nodes    & 1656441 & 13062 & 28611\\
\hline
\end{tabular}
\end{table}

The nonlinear thermo-mechanical coupling behaviors of the inhomogeneous thin plate are researched in the time interval $t\in[0,1]$s. Then we set temporal step as $\Delta t = 0.01$s, the multi-scale nonlinear problem (1.1) and macroscopic homogenized problem (2.18) are simulated separately. Employing the off-line computation results of first example, we go on with on-line computation of nonlinear thermo-mechanical coupling performances of the investigated heterogeneous plate. Then, Figs.\hspace{1mm}8-11 display the computational results $T^{[0]}$, $T^{[1\xi]}$, $T^{[2\xi]}$ and $u_3^{[0]}$, $u_3^{[1\xi]}$, $u_3^{[2\xi]}$ in the middle plane $\alpha_3=0.1cm$ and cross section $\alpha_3=0.18cm$ of the inhomogeneous thin plate, respectively.
\begin{figure}[!htb]
\centering
\begin{minipage}[c]{0.29\textwidth}
  \centering
  \includegraphics[width=37mm]{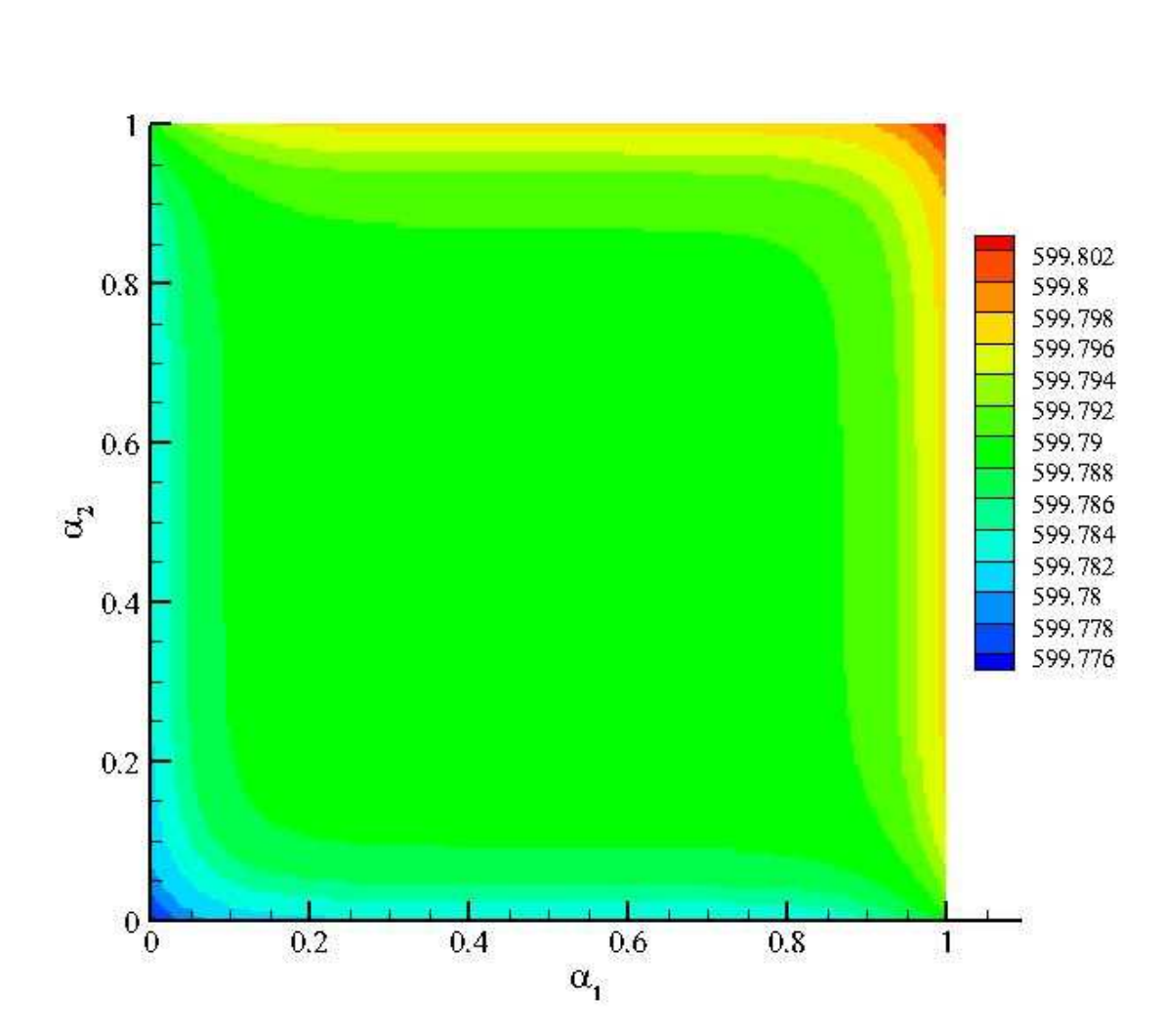}\\
  (a)
\end{minipage}
\begin{minipage}[c]{0.29\textwidth}
  \centering
  \includegraphics[width=37mm]{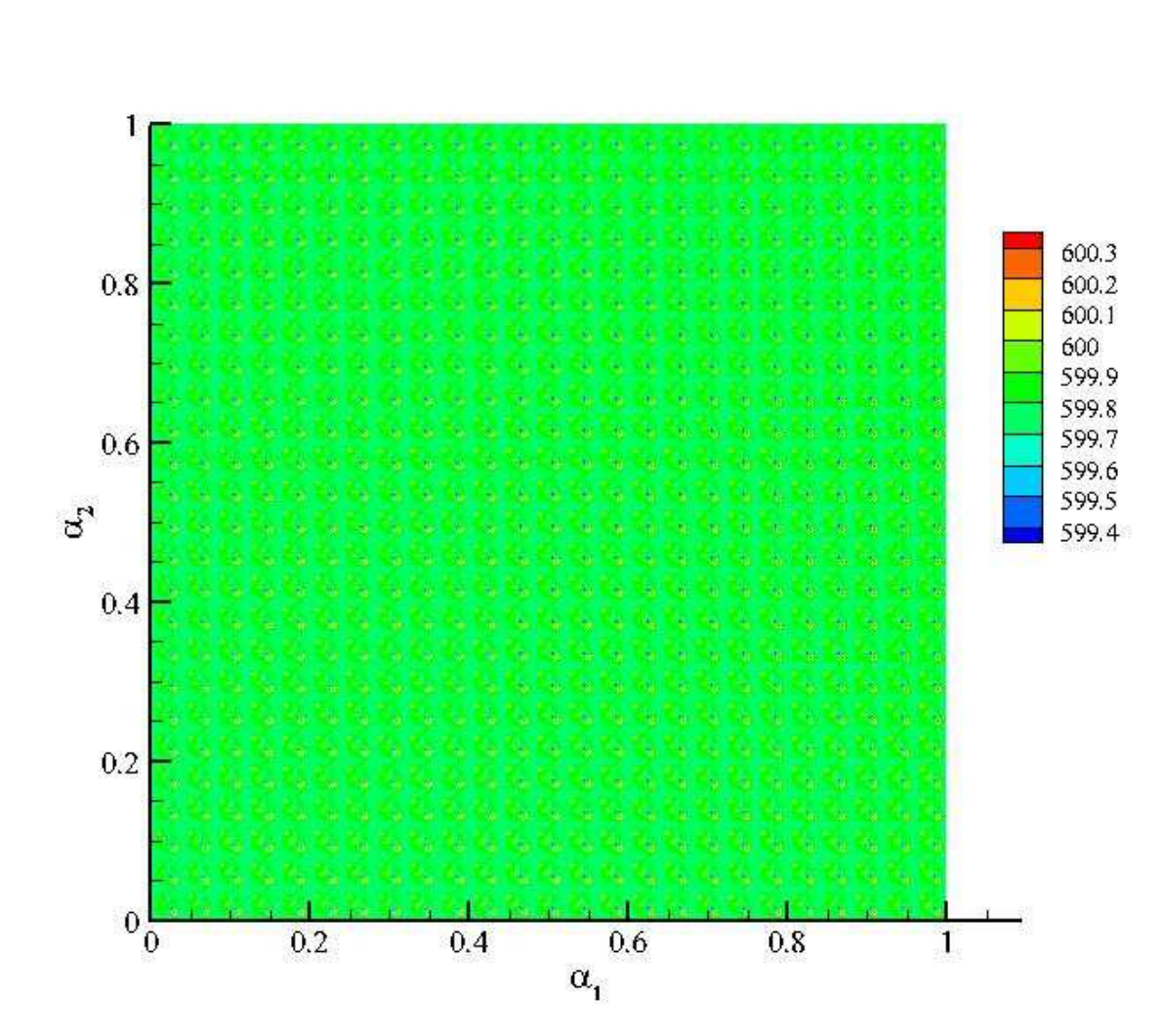}\\
  (b)
\end{minipage}
\begin{minipage}[c]{0.29\textwidth}
  \centering
  \includegraphics[width=37mm]{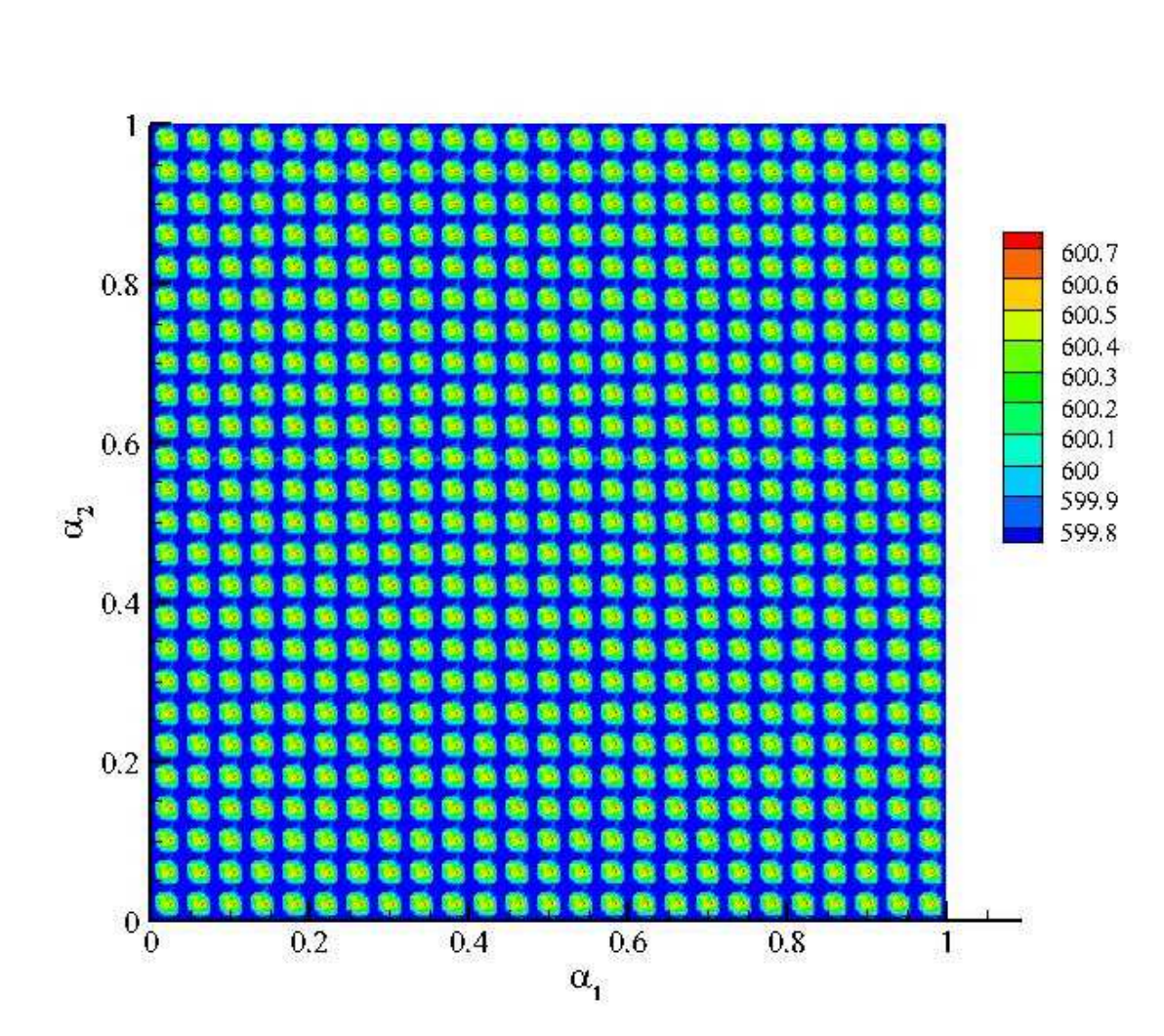}\\
  (c)
\end{minipage}
\caption{The temperature field in cross section $\alpha_3=0.1cm$ at time $t$=1.0s: (a) $T^{[0]}$; (b) $T^{[1\xi]}$; (c) $T^{[2\xi]}$.}
\end{figure}
\begin{figure}[!htb]
\centering
\begin{minipage}[c]{0.29\textwidth}
  \centering
  \includegraphics[width=37mm]{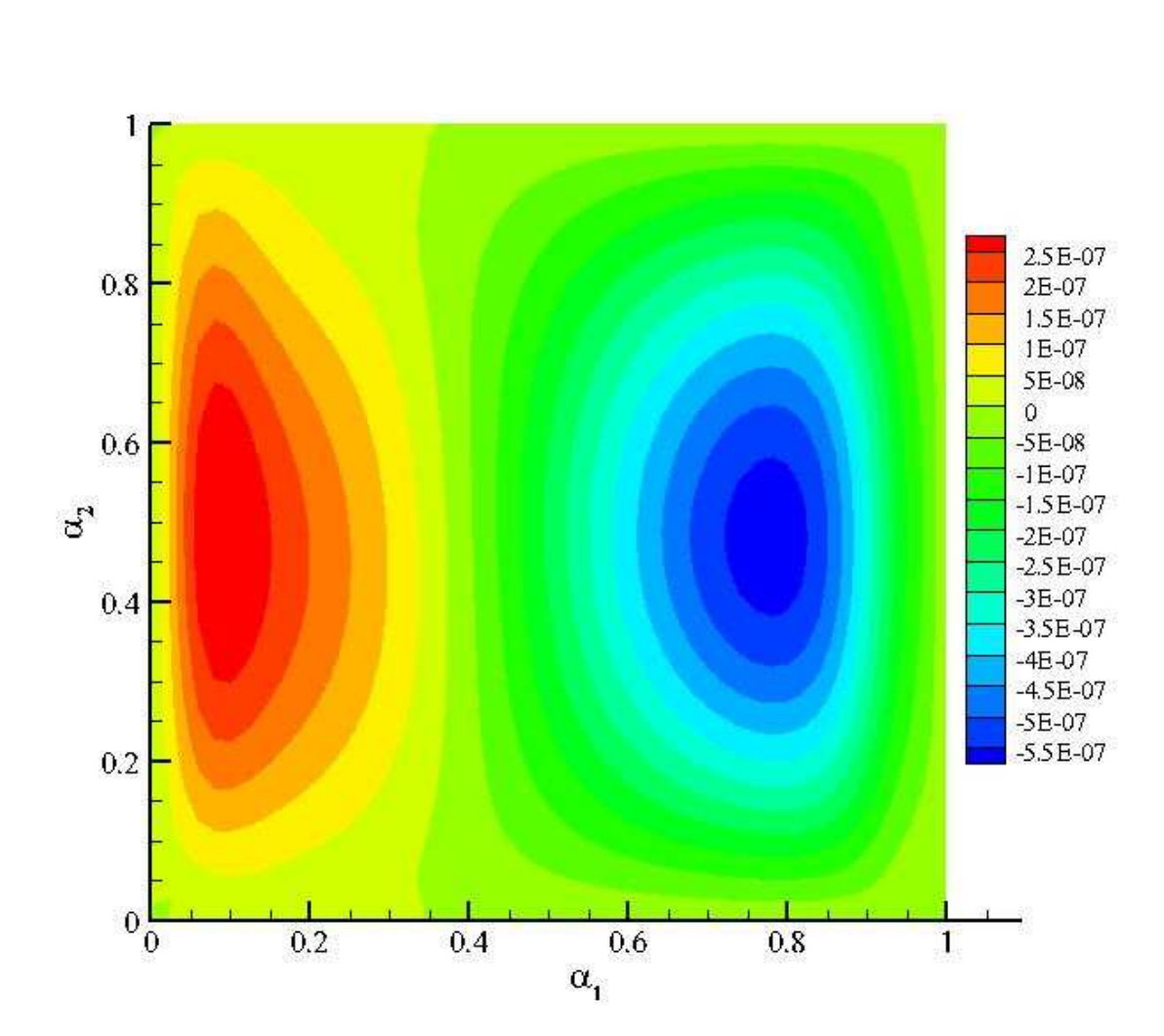}\\
  (a)
\end{minipage}
\begin{minipage}[c]{0.29\textwidth}
  \centering
  \includegraphics[width=37mm]{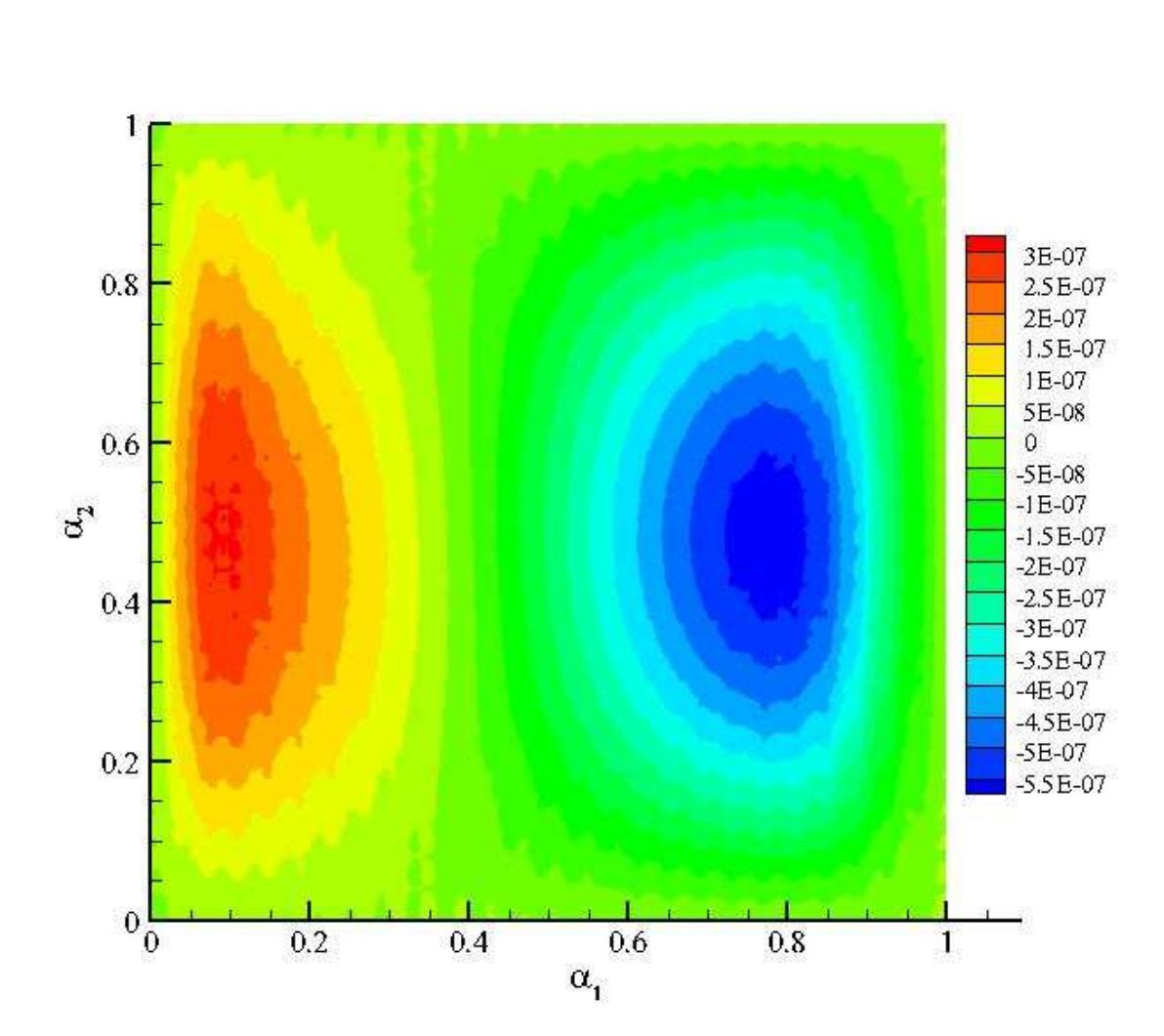}\\
  (b)
\end{minipage}
\begin{minipage}[c]{0.29\textwidth}
  \centering
  \includegraphics[width=37mm]{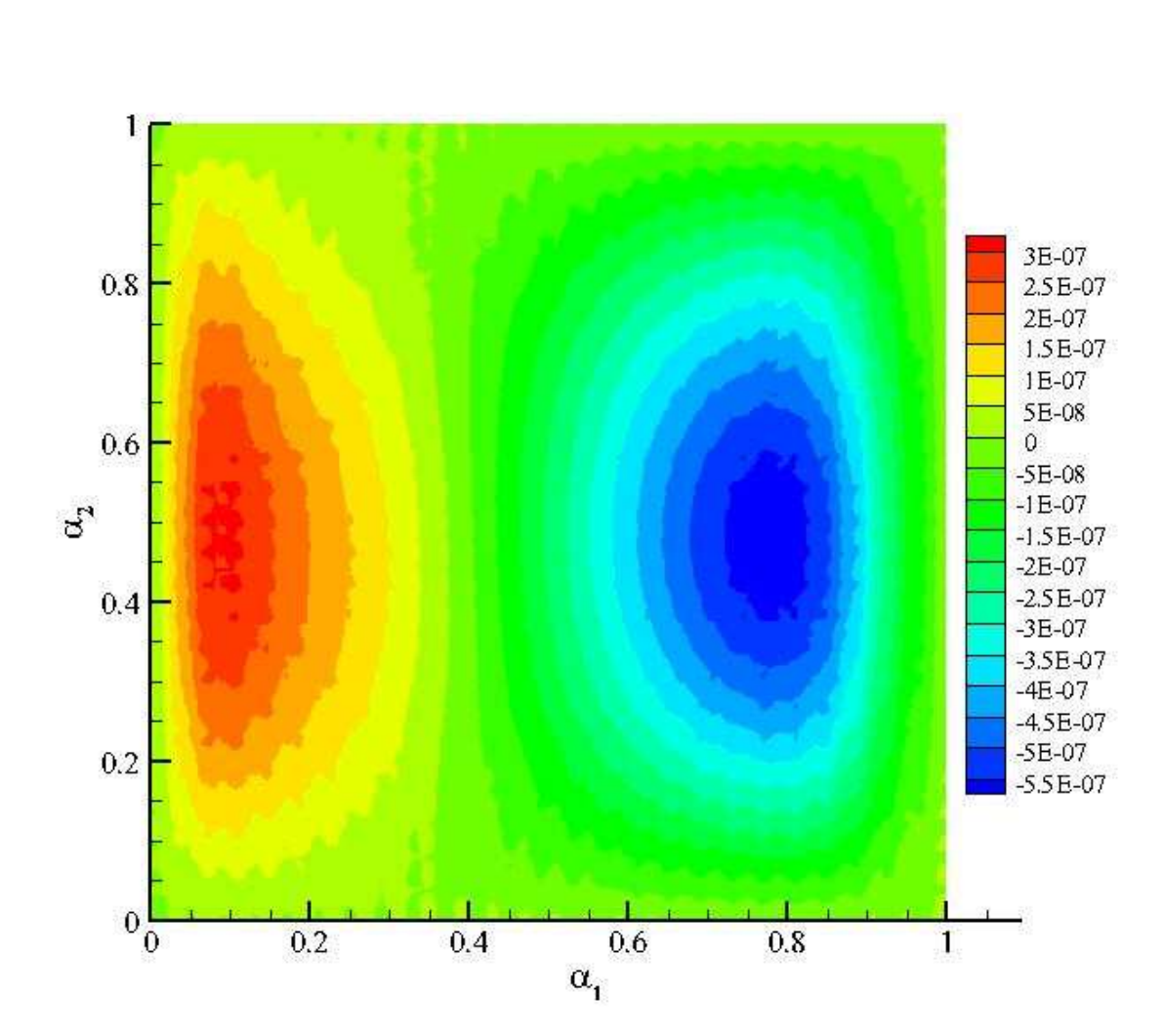}\\
  (c)
\end{minipage}
\caption{The first displacement component in cross section $\alpha_3=0.1cm$ at time $t$=1.0s: (a) $u_1^{[0]}$; (b) $u_1^{[1\xi]}$; (c) $u_1^{[2\xi]}$.}
\end{figure}

\begin{figure}[!htb]
\centering
\begin{minipage}[c]{0.29\textwidth}
  \centering
  \includegraphics[width=37mm]{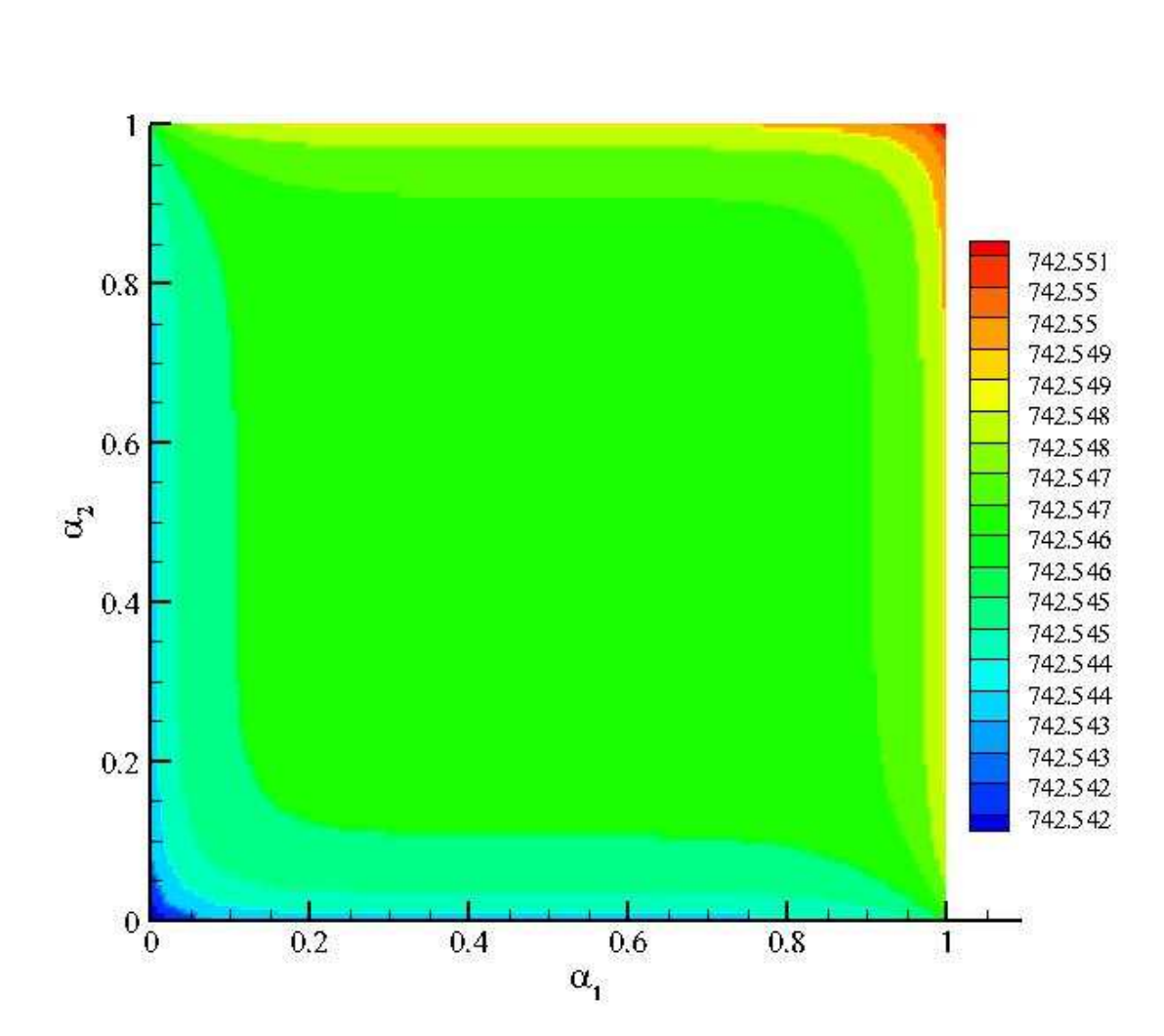}\\
  (a)
\end{minipage}
\begin{minipage}[c]{0.29\textwidth}
  \centering
  \includegraphics[width=37mm]{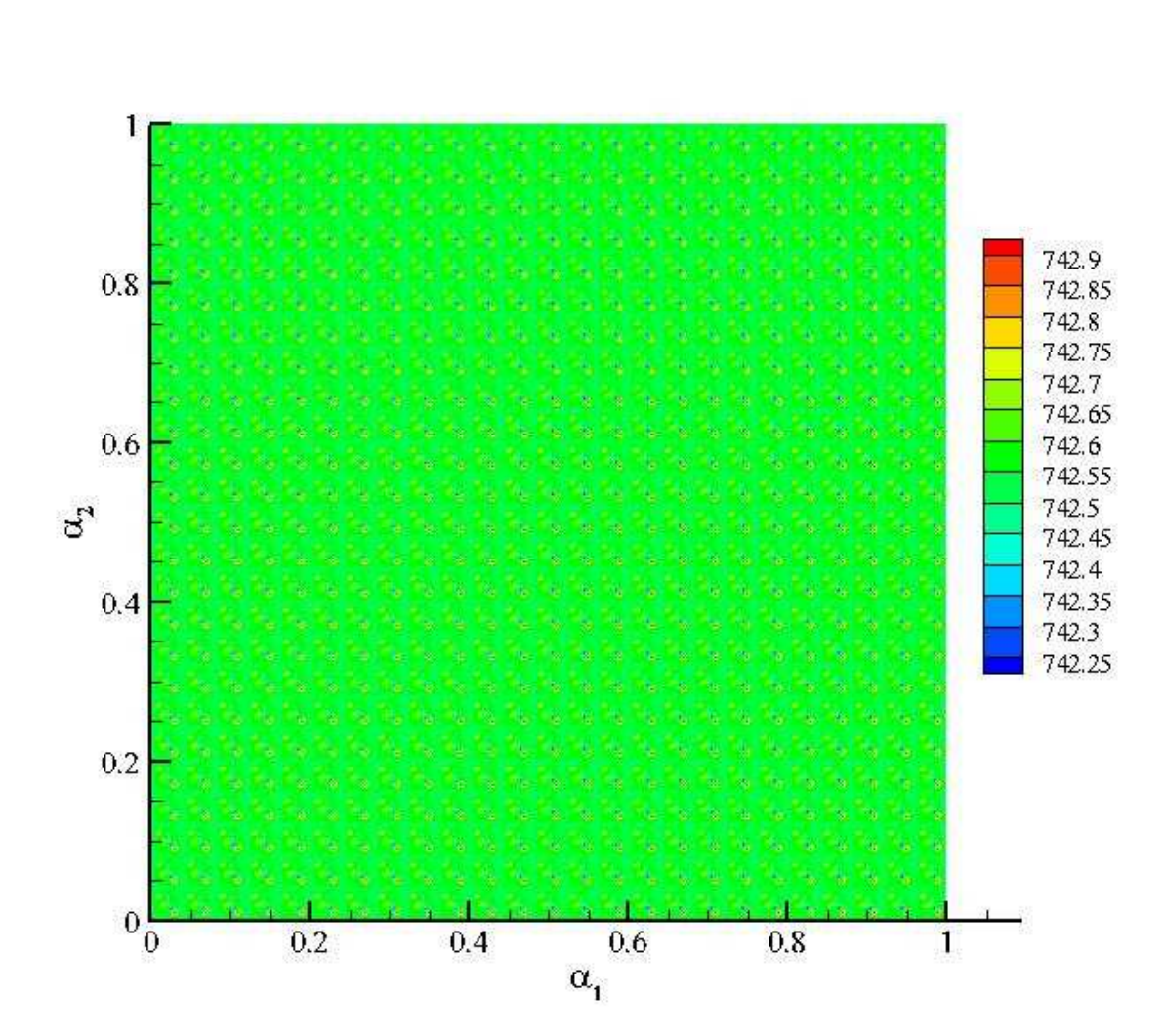}\\
  (b)
\end{minipage}
\begin{minipage}[c]{0.29\textwidth}
  \centering
  \includegraphics[width=37mm]{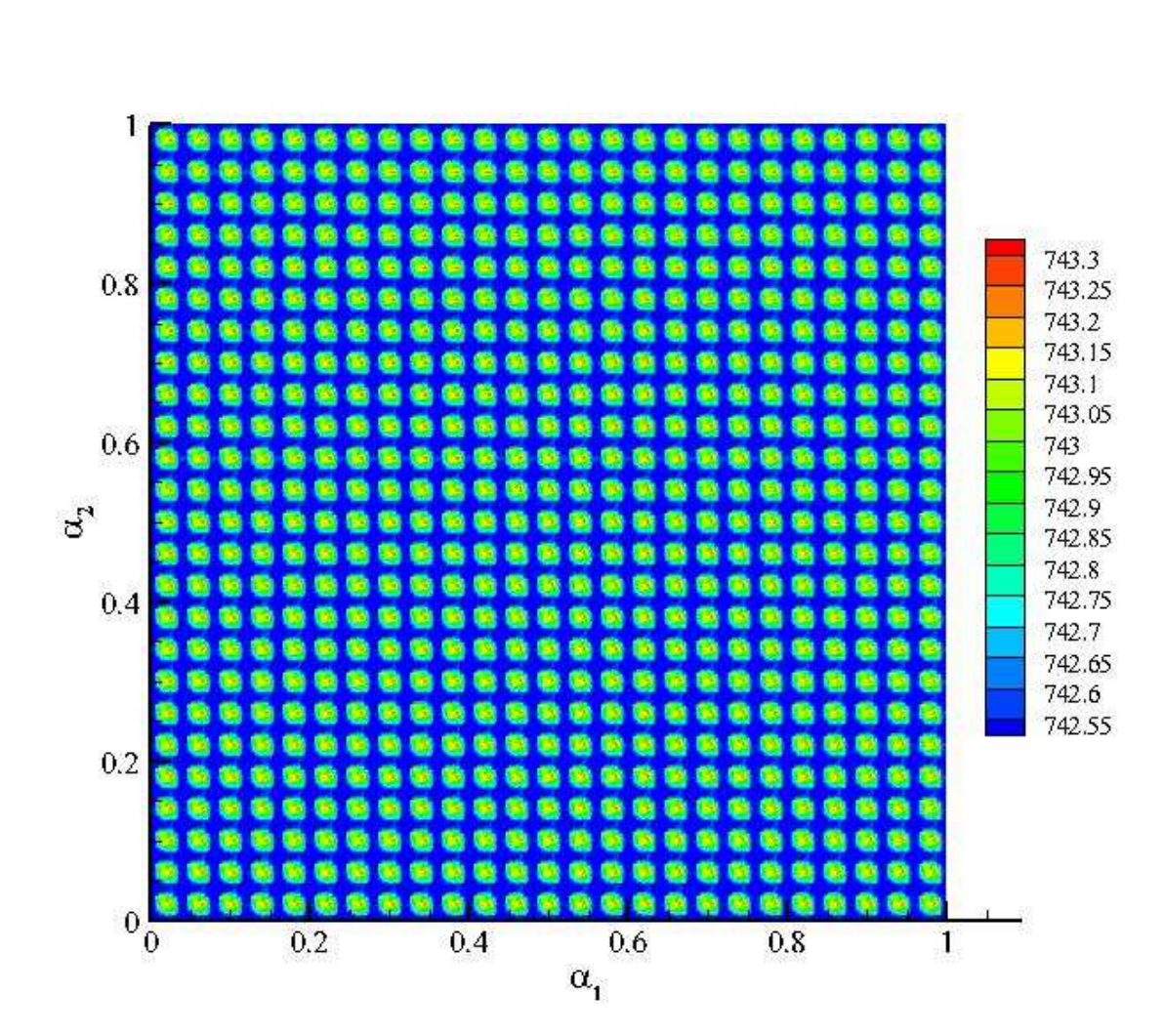}\\
  (c)
\end{minipage}
\caption{The temperature field in cross section $\alpha_3=0.18cm$ at time $t$=1.0s: (a) $T^{[0]}$; (b) $T^{[1\xi]}$; (c) $T^{[2\xi]}$.}
\end{figure}
\begin{figure}[!htb]
\centering
\begin{minipage}[c]{0.29\textwidth}
  \centering
  \includegraphics[width=37mm]{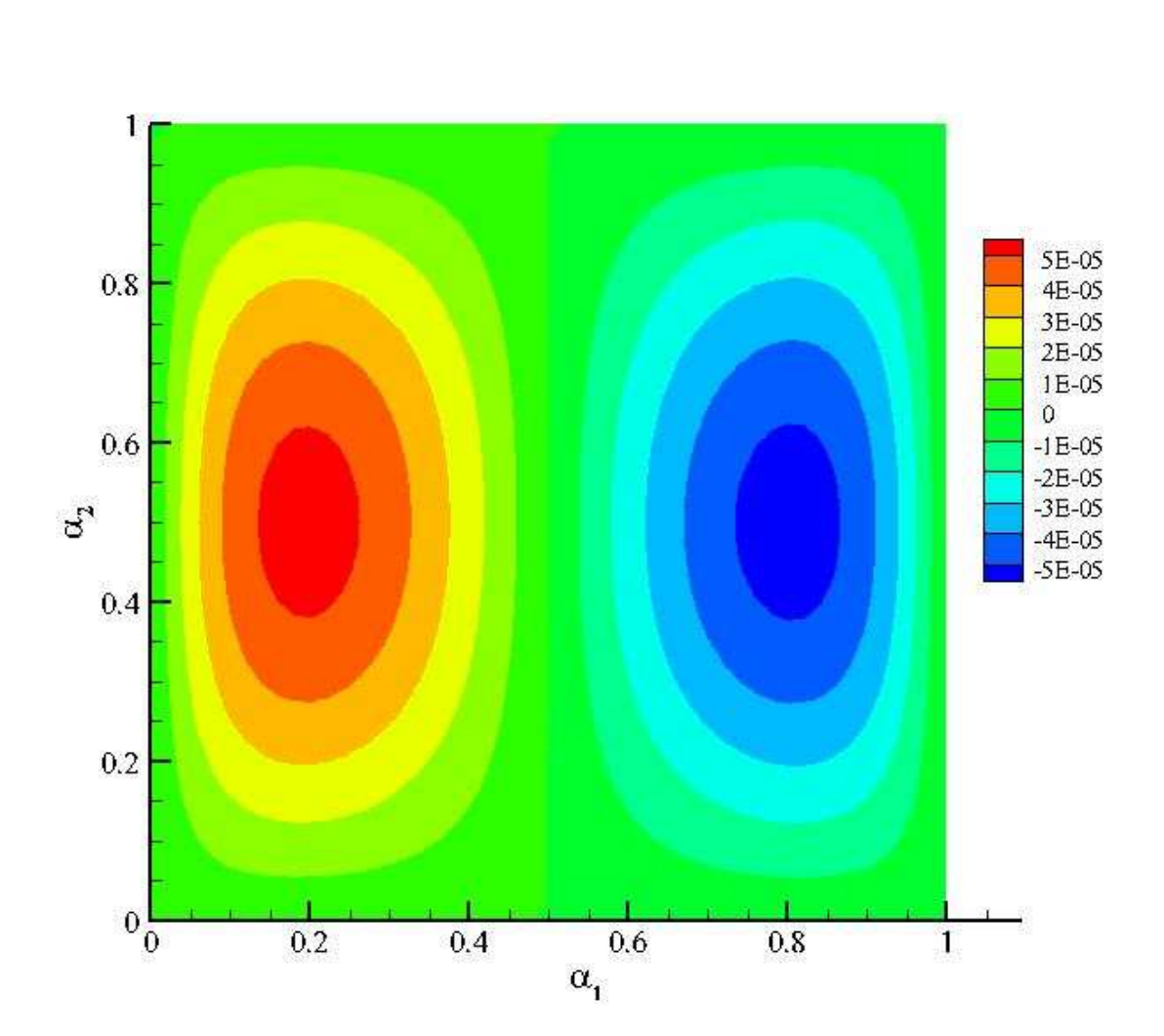}\\
  (a)
\end{minipage}
\begin{minipage}[c]{0.29\textwidth}
  \centering
  \includegraphics[width=37mm]{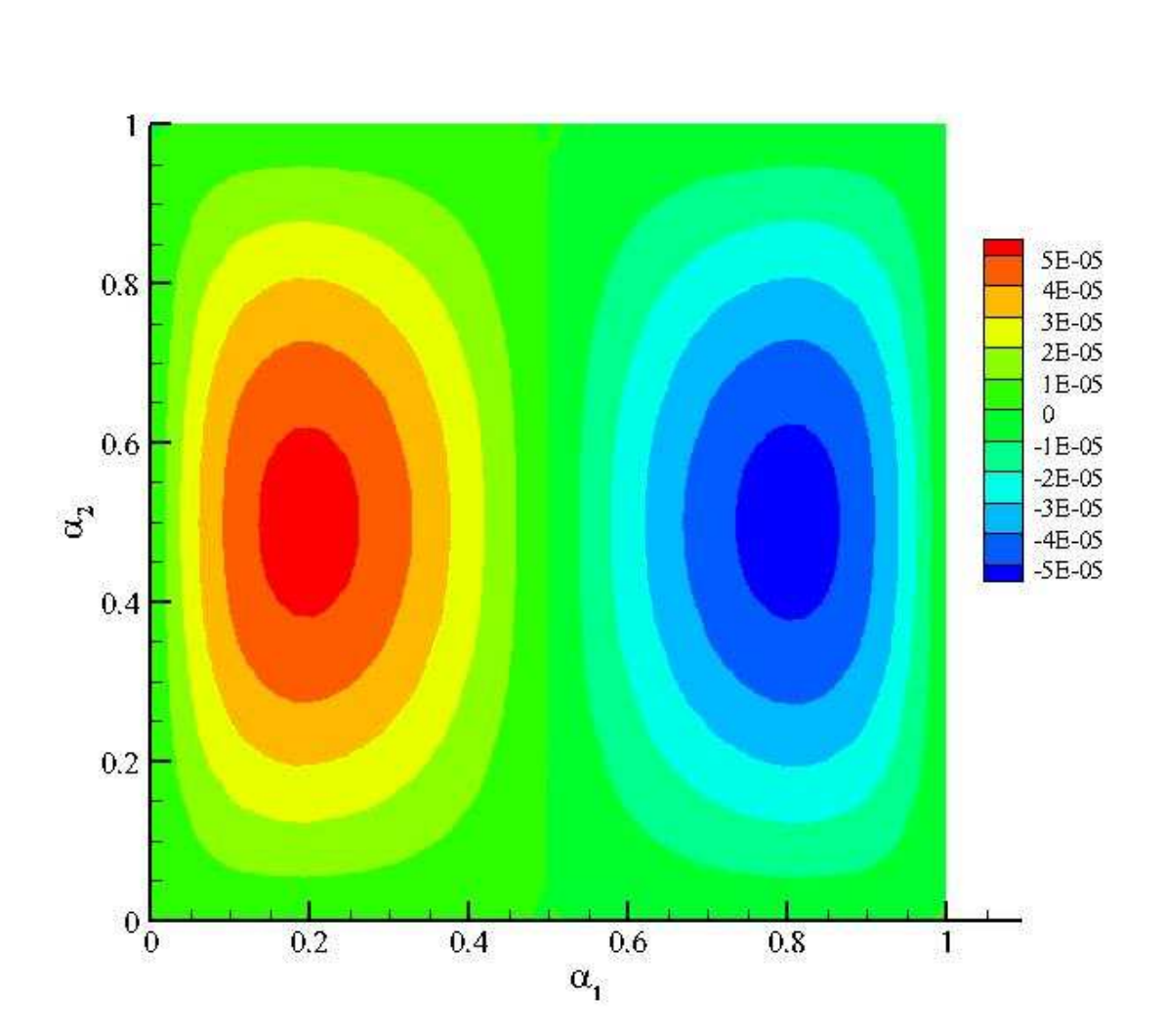}\\
  (b)
\end{minipage}
\begin{minipage}[c]{0.29\textwidth}
  \centering
  \includegraphics[width=37mm]{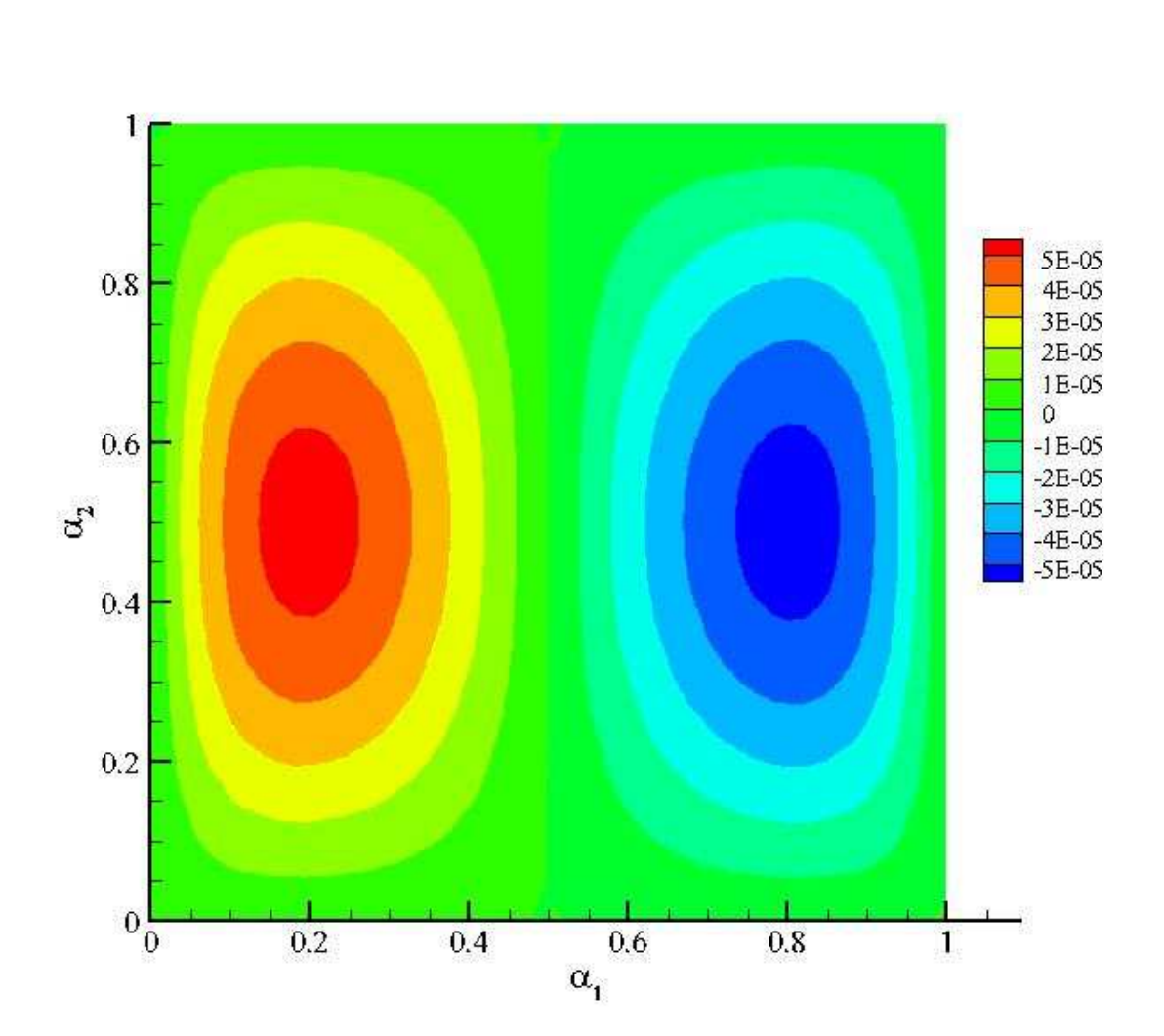}\\
  (c)
\end{minipage}
\caption{The first displacement component in cross section $\alpha_3=0.18cm$ at time $t$=1.0s: (a) $u_1^{[0]}$; (b) $u_1^{[1\xi]}$; (c) $u_1^{[2\xi]}$.}
\end{figure}

Observing the Table 3, we can deduce that the proposed HOMS approach is pretty economical to compute the nonlinear thermo-mechanical problem of heterogeneous thin plate contrast to full-scale DNS, that can notably conserve computer storage and time without losing precision. Figs.\hspace{1mm}8-11 testify that only HOMS solutions can precisely catch the microscopic steep fluctuations caused by the material heterogeneities in inhomogeneous thin plate. Furthermore, the presented HOMS approach is still effective for a relatively small parameter $\xi$, namely existing a great number of microscopic PUCs in inhomogeneous thin plate. For the investigated large-scale problem, the convergence of high-resolution FEM simulation can not be obtained. But the presented HOMS approach can offer the dependable numerical solutions for the time-variant nonlinear thermo-mechanical problem of inhomogeneous thin plate after long-time numerical computation. The foregoing merits of the proposed HOMS approach are very important in large industrial finite element applications.
\subsection{Multi-scale nonlinear simulation of heterogeneous cylindrical shell}
This example studies the nonlinear thermo-mechanical behaviors of a inhomogeneous cylindrical shell with $\displaystyle\xi=\pi/36$ as displayed in Fig.\hspace{1mm}12(a) and Fig.\hspace{1mm}12(b). The Lam$\rm{\acute{e}}$ coefficients of this cylindrical shell are $H_1=1$, $H_2=R_2+\alpha_3$ and $H_3=1$, where $R_2$ is the curvature radius of shell's middle plane in $\alpha_2$ direction with $R_2=\pi$. And the whole domain $\displaystyle\Omega=(\alpha_1,\alpha_2,\alpha_3)=(z,\theta,r-R_2)=[0,3\pi/4]cm\times[-\pi/6,\pi/6]\times[-\pi/12,\pi/12]cm$ and PUC $\Theta$ are displayed in Fig.\hspace{1mm}12(c) and Fig.\hspace{1mm}12(d). It is noteworthy that $z$, $\theta$ and $r$ are the coordinate variables of cylindrical coordinates.
\begin{figure}[!htb]
\centering
\begin{minipage}[c]{0.46\textwidth}
  \centering
  \includegraphics[width=0.6\linewidth,totalheight=1.3in]{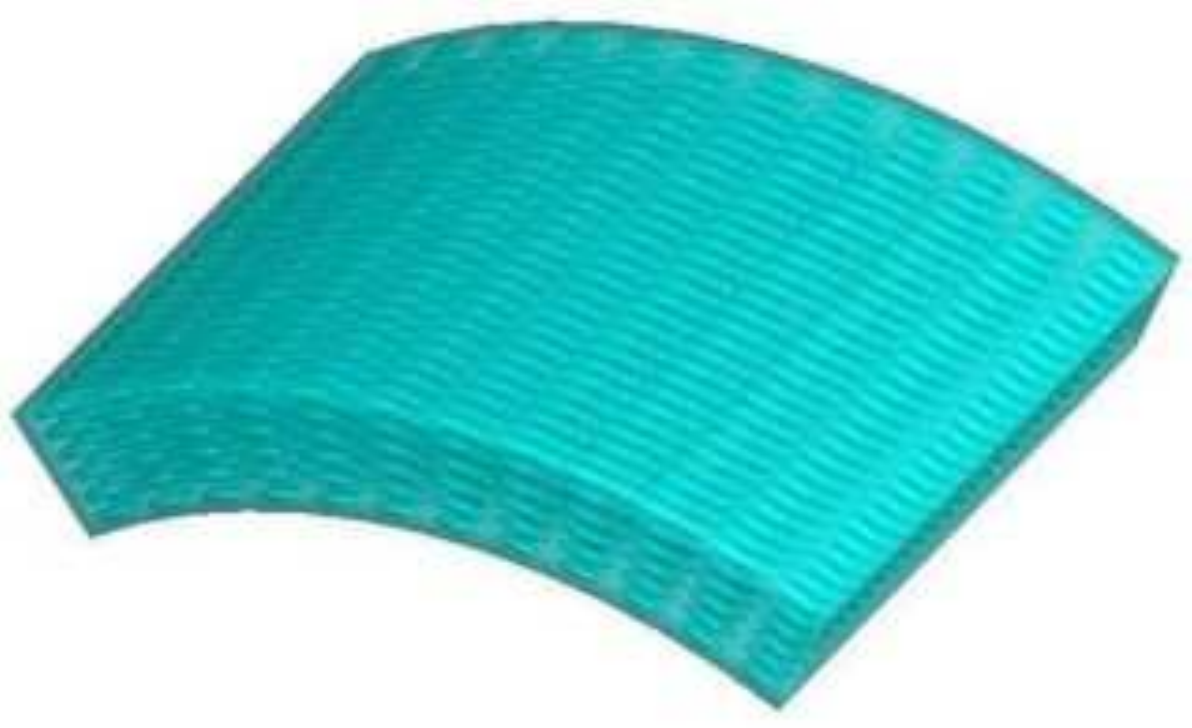} \\
  (a)
\end{minipage}
\begin{minipage}[c]{0.46\textwidth}
  \centering
  \includegraphics[width=0.6\linewidth,totalheight=1.3in]{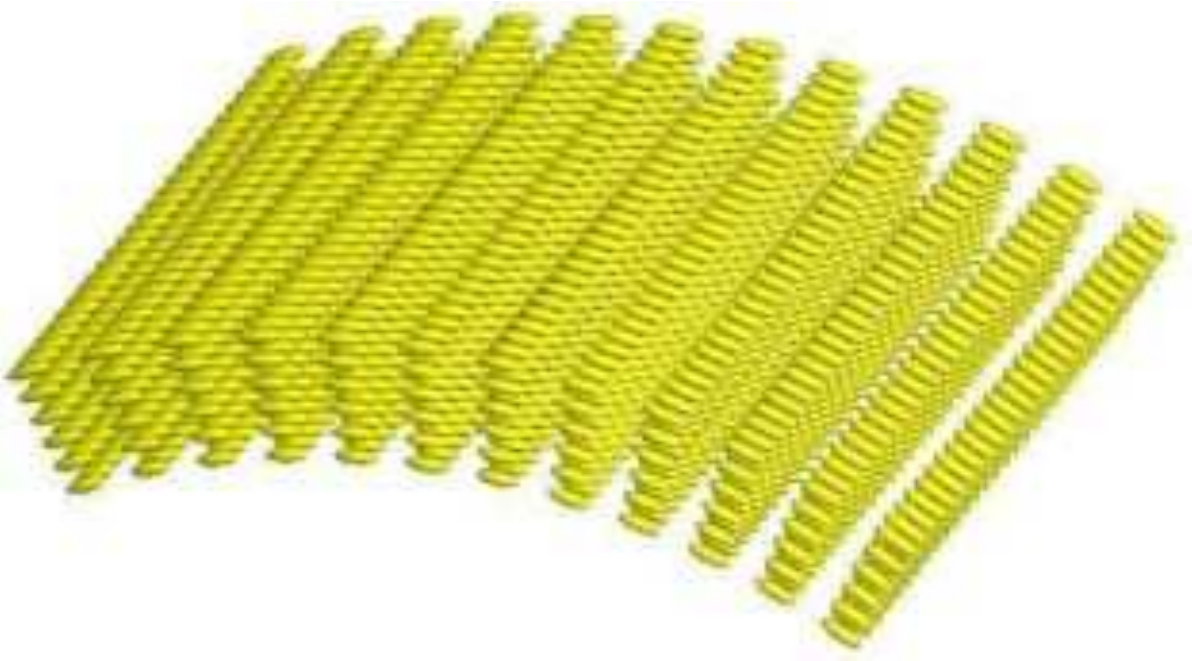} \\
  (b)
\end{minipage}
\begin{minipage}[c]{0.46\textwidth}
  \centering
  \includegraphics[width=0.6\linewidth,totalheight=1.3in]{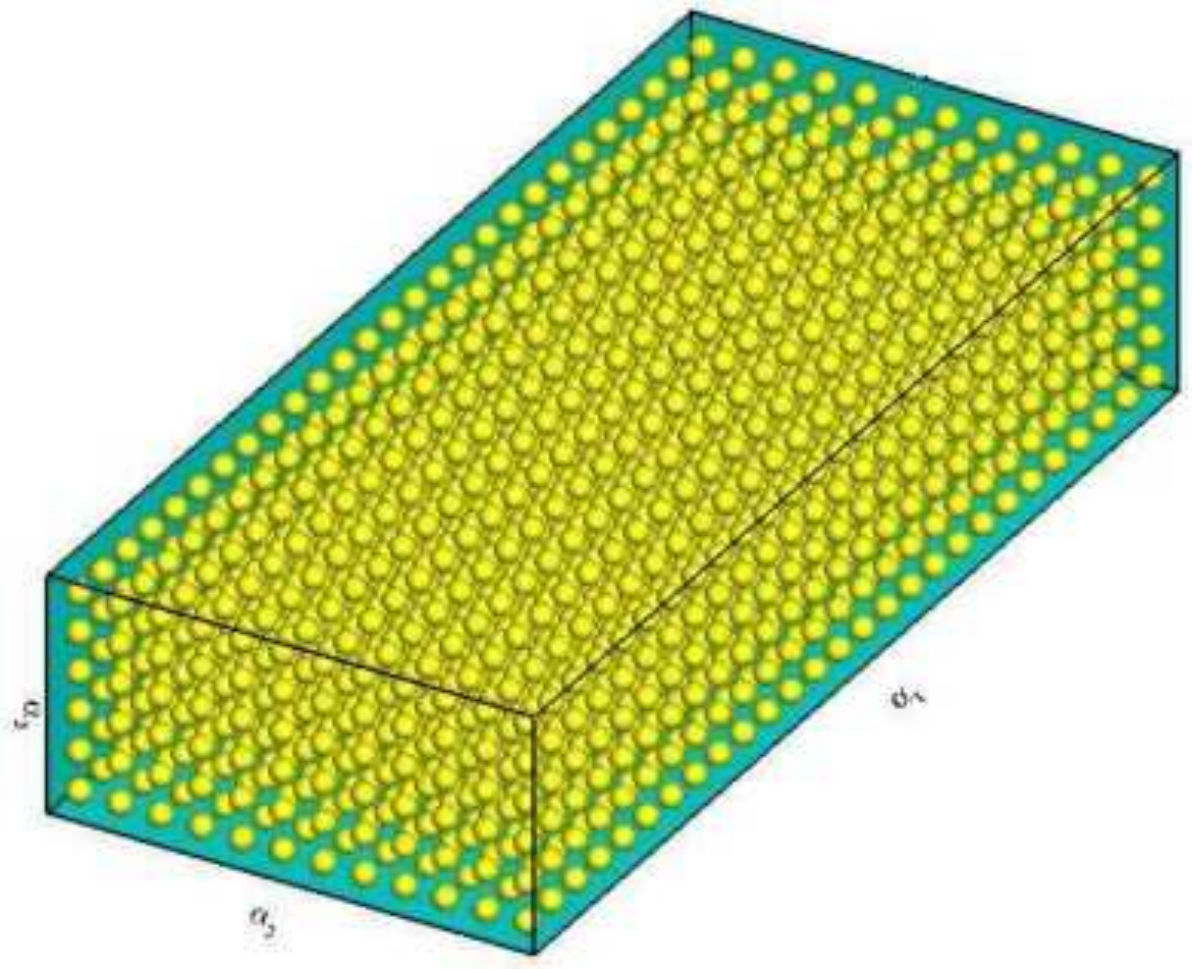} \\
  (c)
\end{minipage}
\begin{minipage}[c]{0.46\textwidth}
  \centering
  \includegraphics[width=0.6\linewidth,totalheight=1.3in]{EX4Nn-new-eps-converted-to.pdf} \\
  (d)
\end{minipage}
\caption{(a) The macrostructure of heterogeneous cylindrical shell; (b) inclusion distribution in heterogeneous cylindrical shell; (c) the whole domain $\Omega$; (d) the unit cell $\Theta$.}
\end{figure}

This example does not give the reference solutions $T_{Fe}^\xi(\bm{\alpha},t)$ and $\bm{u}_{Fe}^\xi(\bm{\alpha},t)$ owing to the same reason as the second example. Moreover, the material parameters of this heterogeneous cylindrical shell are defined with the same values as first example. The investigated cylindrical shell $\Omega$ is clamped on its four side surfaces, that are perpendicular to $\alpha_2$-axis. The initial temperature is remained at $373.15K$ on the bottom surface. Furthermore, we set $h = 5000J/(cm^3\bm{\cdot}s)$ and $(f_1,f_2,f_3) = (0,0,-40000)N/cm^3$. Next, the information of FEM meshes is illustrated in Table 4.
\begin{table}[!htb]{\caption{Summary of computational cost.}\label{t2}}
\centering
\begin{tabular}{cccc}
\hline
 & Multi-scale eqs. & Cell eqs. & Homogenized eqs. \\
\hline
FEM elements & $\approx$1944$\times$75466 & 75466 & 746496\\
FEM nodes    & $\approx$1944$\times$13062 & 13062 & 133525\\
\hline
\end{tabular}
\end{table}

In this example, we demand to off-line compute 17500 times auxiliary cell problems in total, where the number of first-order cell functions and second-order cell functions is 33 and 250 respectively. Also, we assign 5 equidistant interpolating points toward the thickness direction and 14 macroscopic interpolation temperature in each unit cell. Whereupon, the nonlinear thermo-mechanical coupling performance of the heterogeneous cylindrical shell is studied in the time interval $t\in[0,1]$s. Setting the temporal step $\Delta t = 0.02$s, the multi-scale equations (1.1) and corresponding homogenized equations (2.18) are on-line computed separately. Afterwards, the numerical results at $t=0.2$s and $t=1.0$s are separately displayed in Figs.\hspace{1mm}13 and 14, and Figs.\hspace{1mm}15 and 16.
\begin{figure}[!htb]
\centering
\begin{minipage}[c]{0.29\textwidth}
  \centering
  \includegraphics[width=37mm]{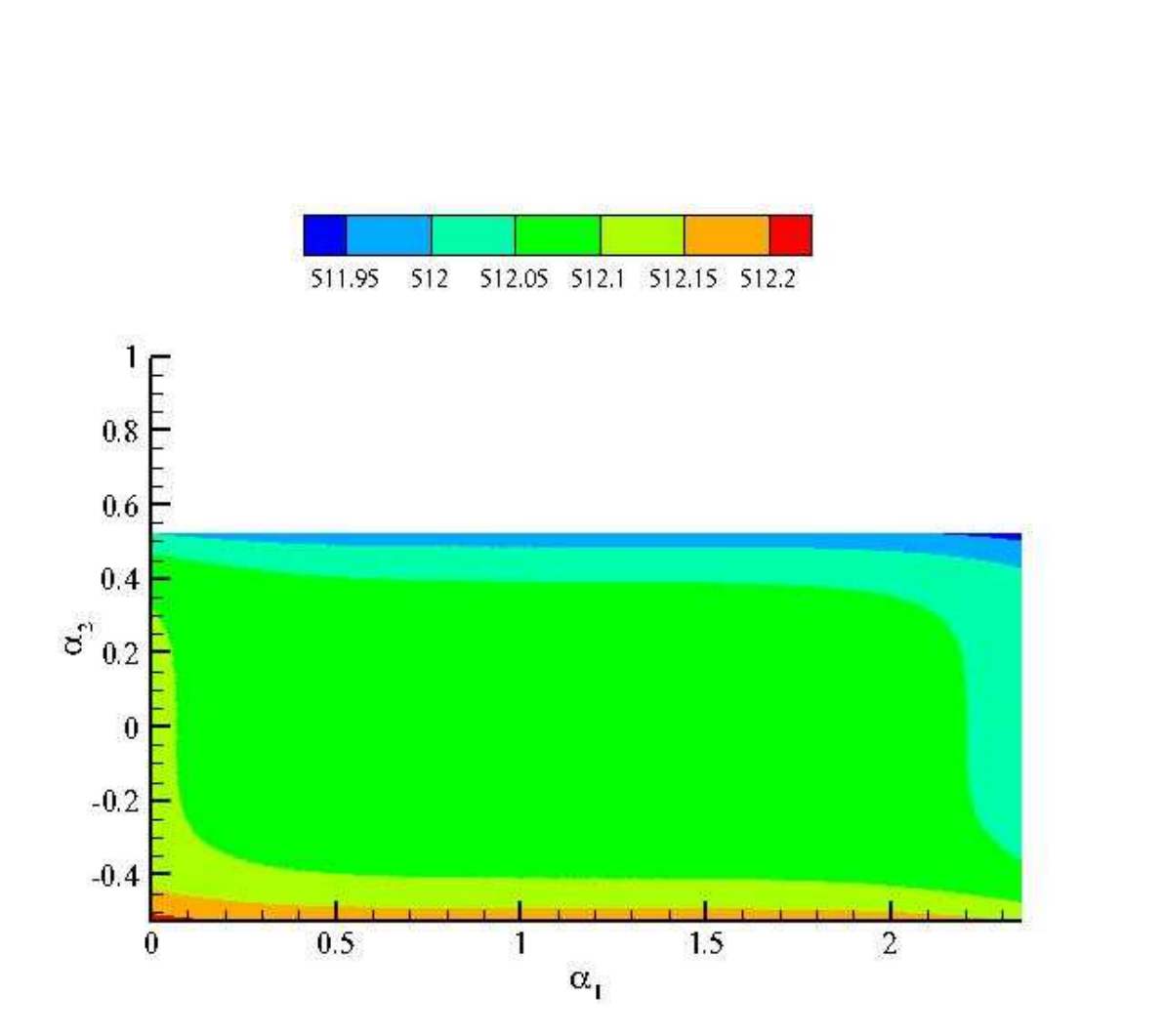}\\
  (a)
\end{minipage}
\begin{minipage}[c]{0.29\textwidth}
  \centering
  \includegraphics[width=37mm]{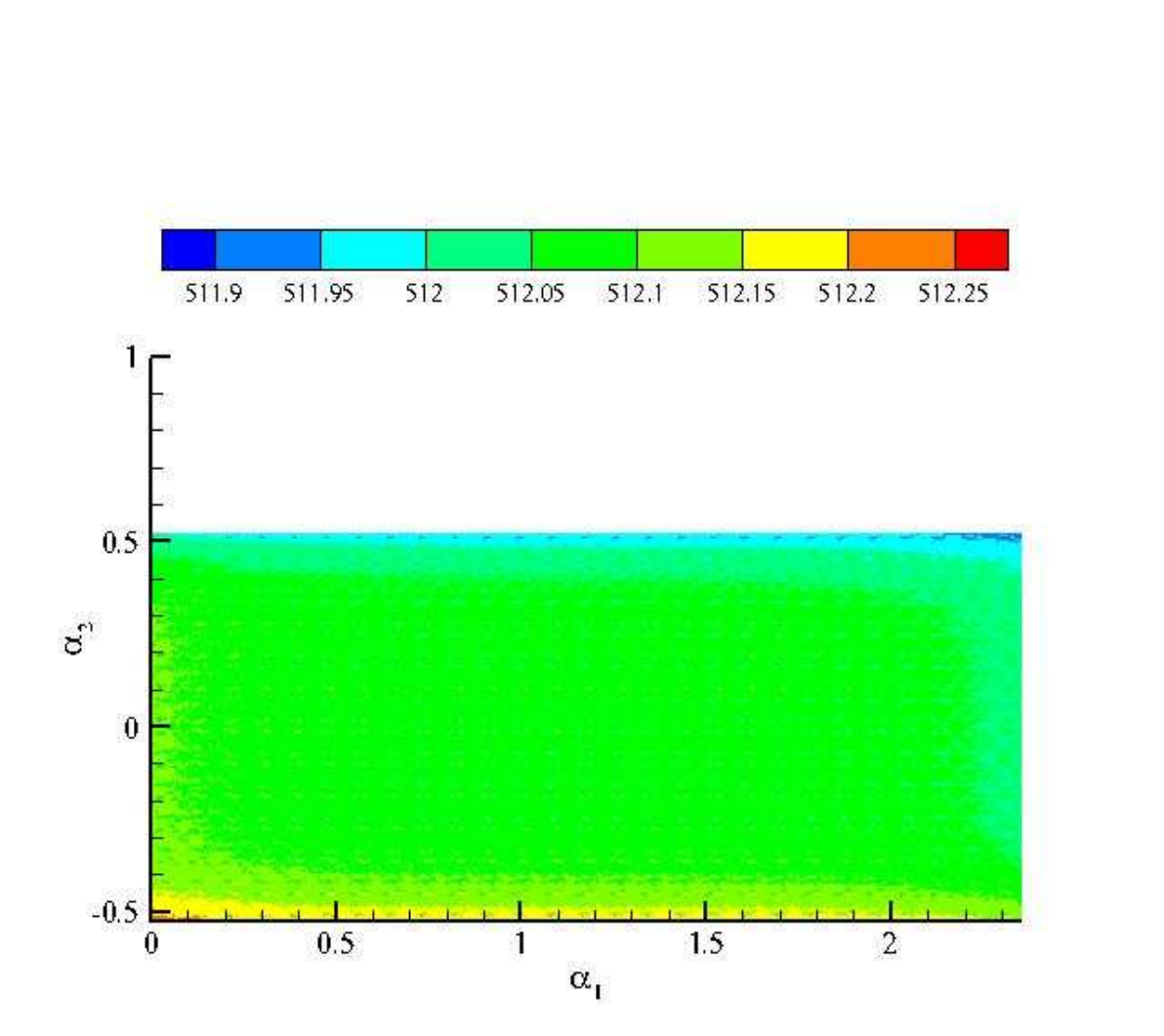}\\
  (b)
\end{minipage}
\begin{minipage}[c]{0.29\textwidth}
  \centering
  \includegraphics[width=37mm]{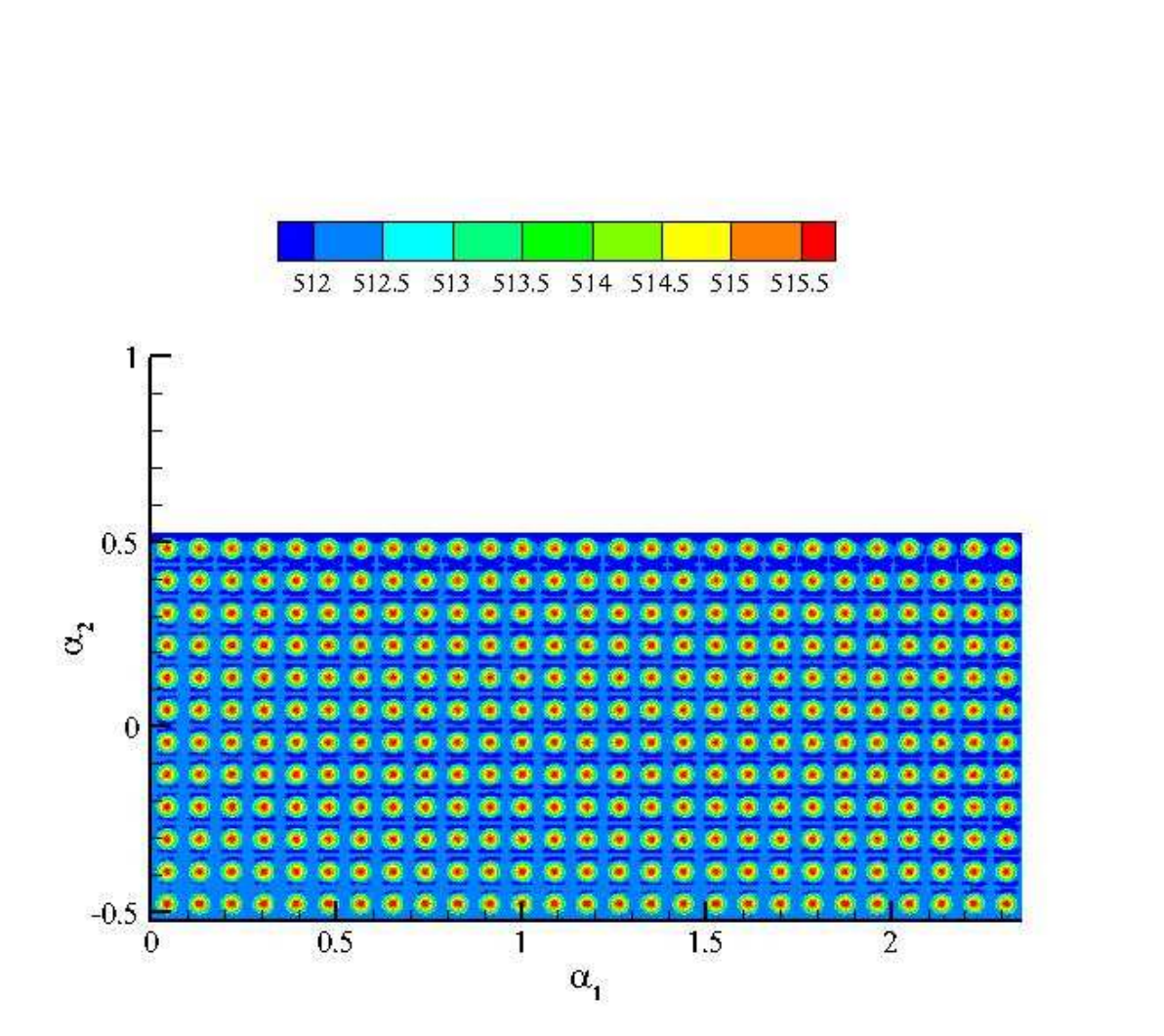}\\
  (c)
\end{minipage}
\caption{The temperature field in cross section $\alpha_3={\pi}/{72}cm$ at time $t$=0.2s: (a) $T^{[0]}$; (b) $T^{[1\xi]}$; (c) $T^{[2\xi]}$.}
\end{figure}
\begin{figure}[!htb]
\centering
\begin{minipage}[c]{0.29\textwidth}
  \centering
  \includegraphics[width=37mm]{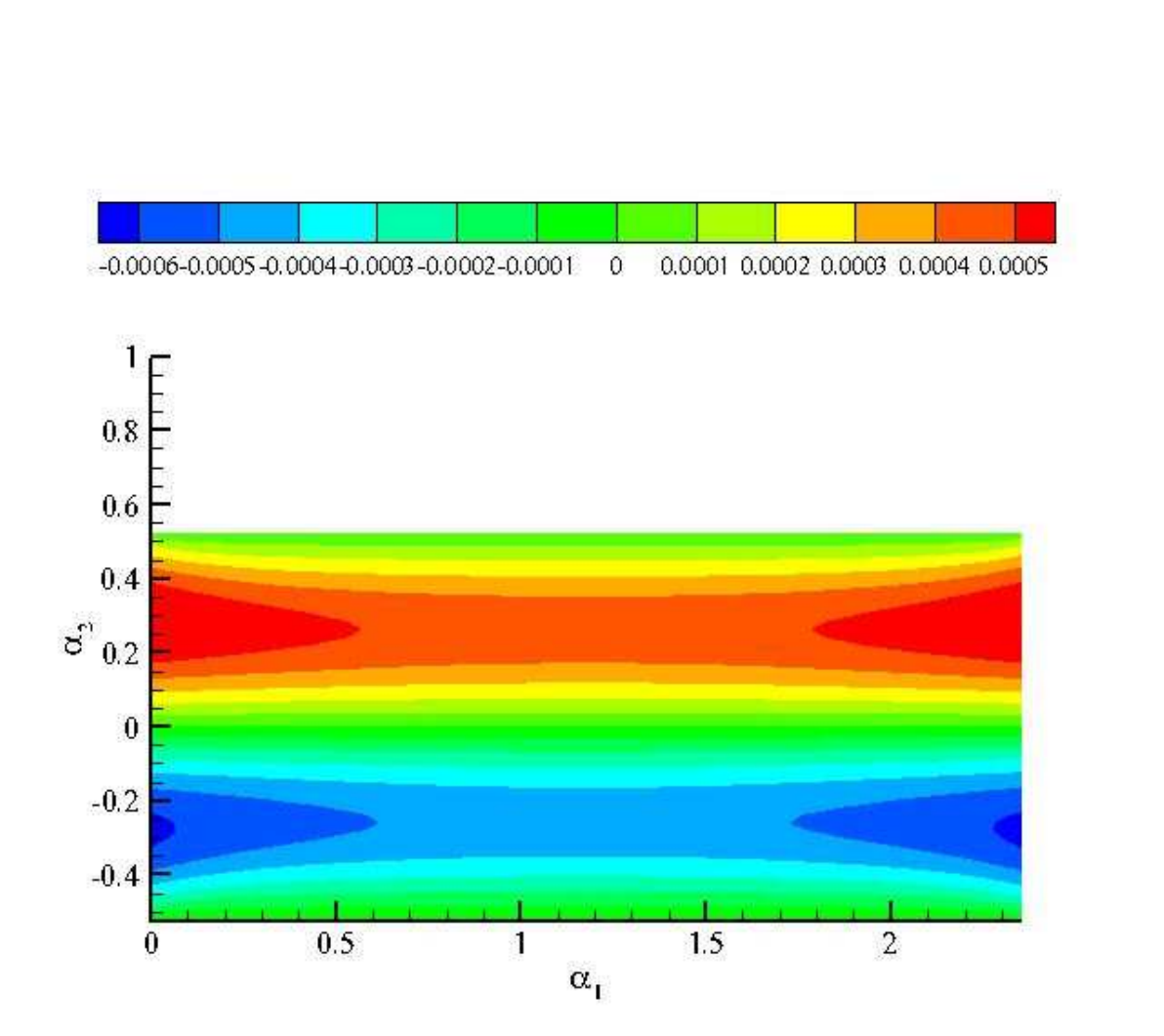}\\
  (a)
\end{minipage}
\begin{minipage}[c]{0.29\textwidth}
  \centering
  \includegraphics[width=37mm]{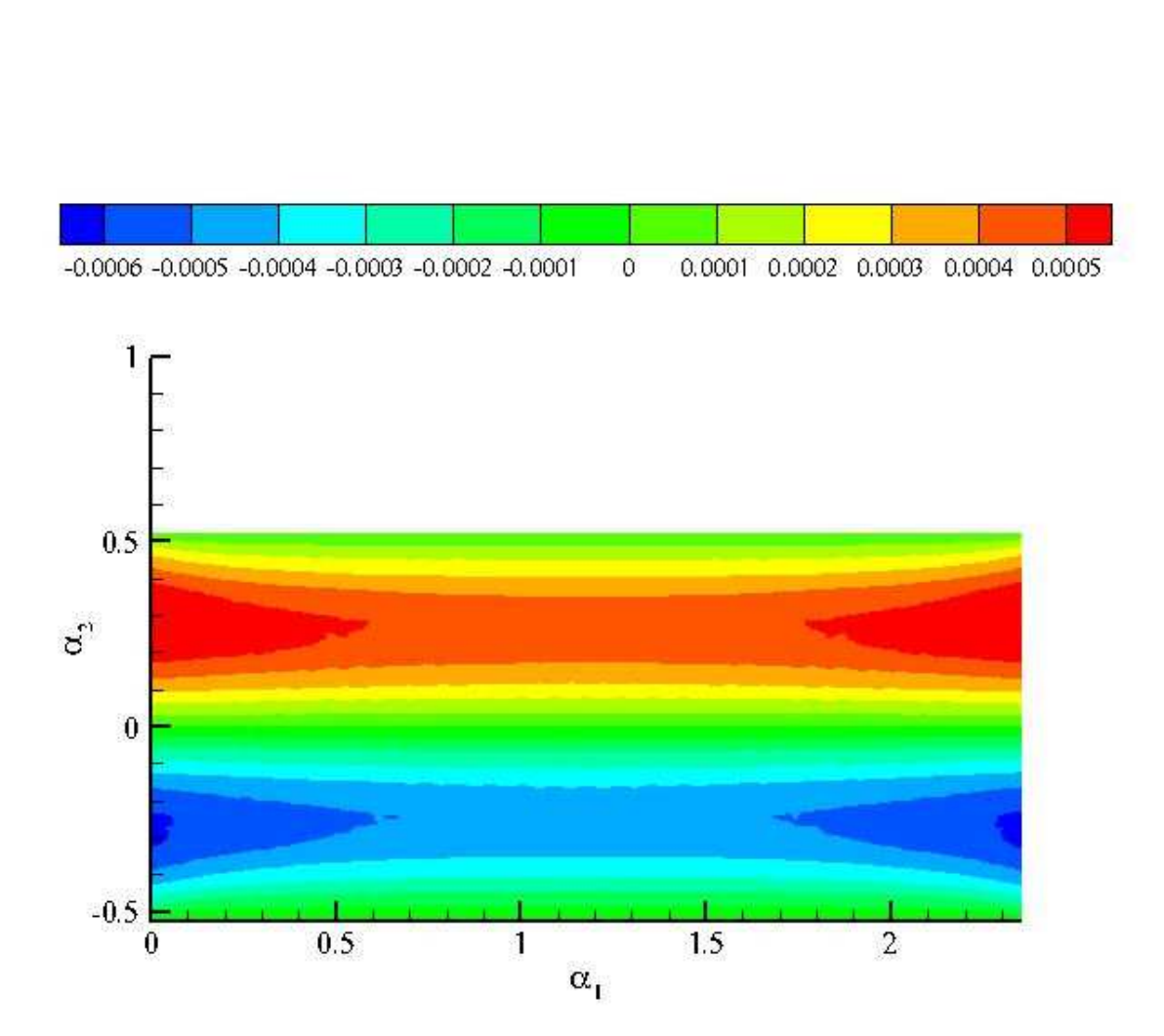}\\
  (b)
\end{minipage}
\begin{minipage}[c]{0.29\textwidth}
  \centering
  \includegraphics[width=37mm]{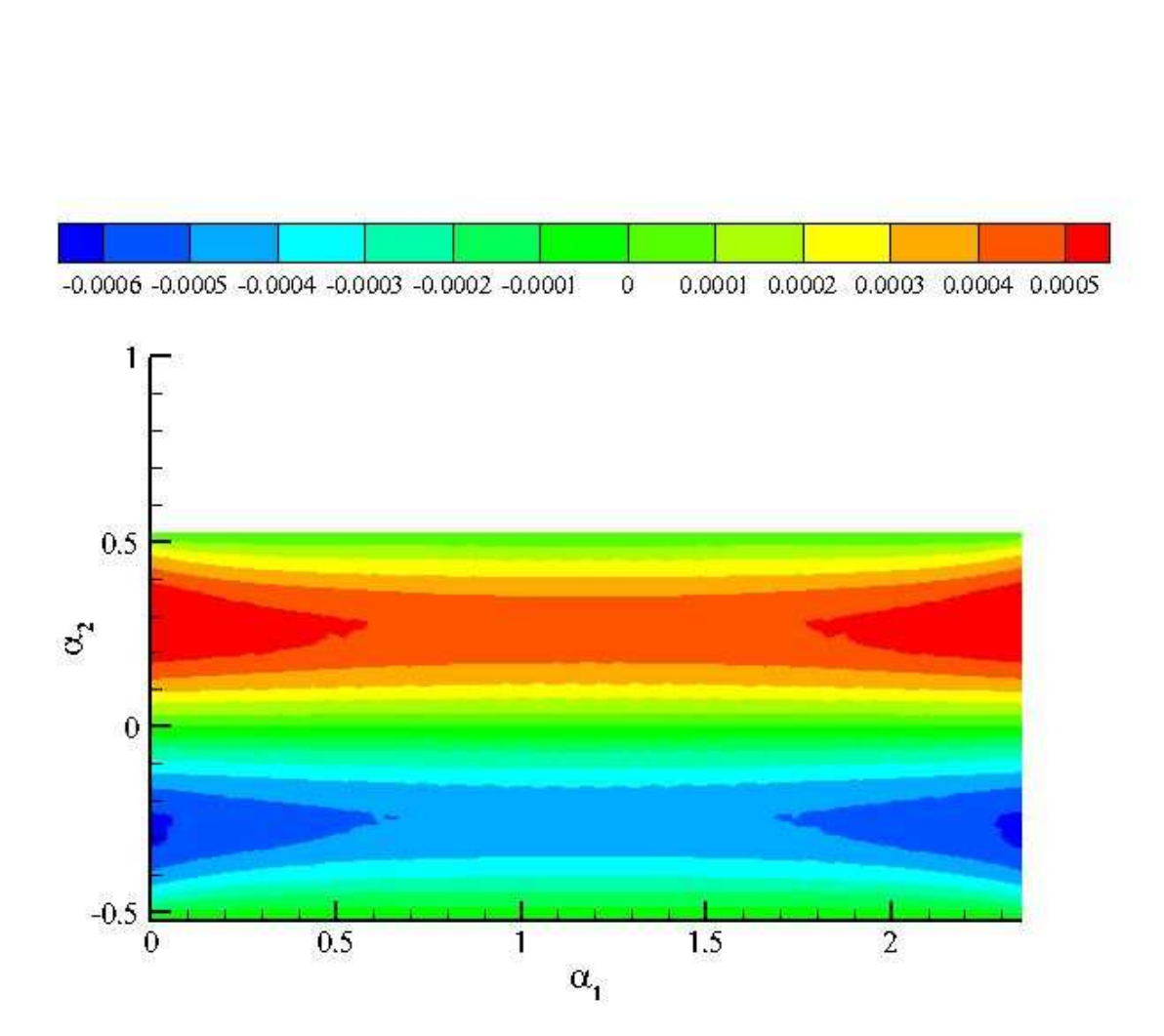}\\
  (c)
\end{minipage}
\caption{The second displacement component in cross section $\alpha_3={\pi}/{72}cm$ at time $t$=0.2s: (a) $u_2^{[0]}$; (b) $u_2^{[1\xi]}$; (c) $u_2^{[2\xi]}$.}
\end{figure}
\begin{figure}[!htb]
\centering
\begin{minipage}[c]{0.29\textwidth}
  \centering
  \includegraphics[width=37mm]{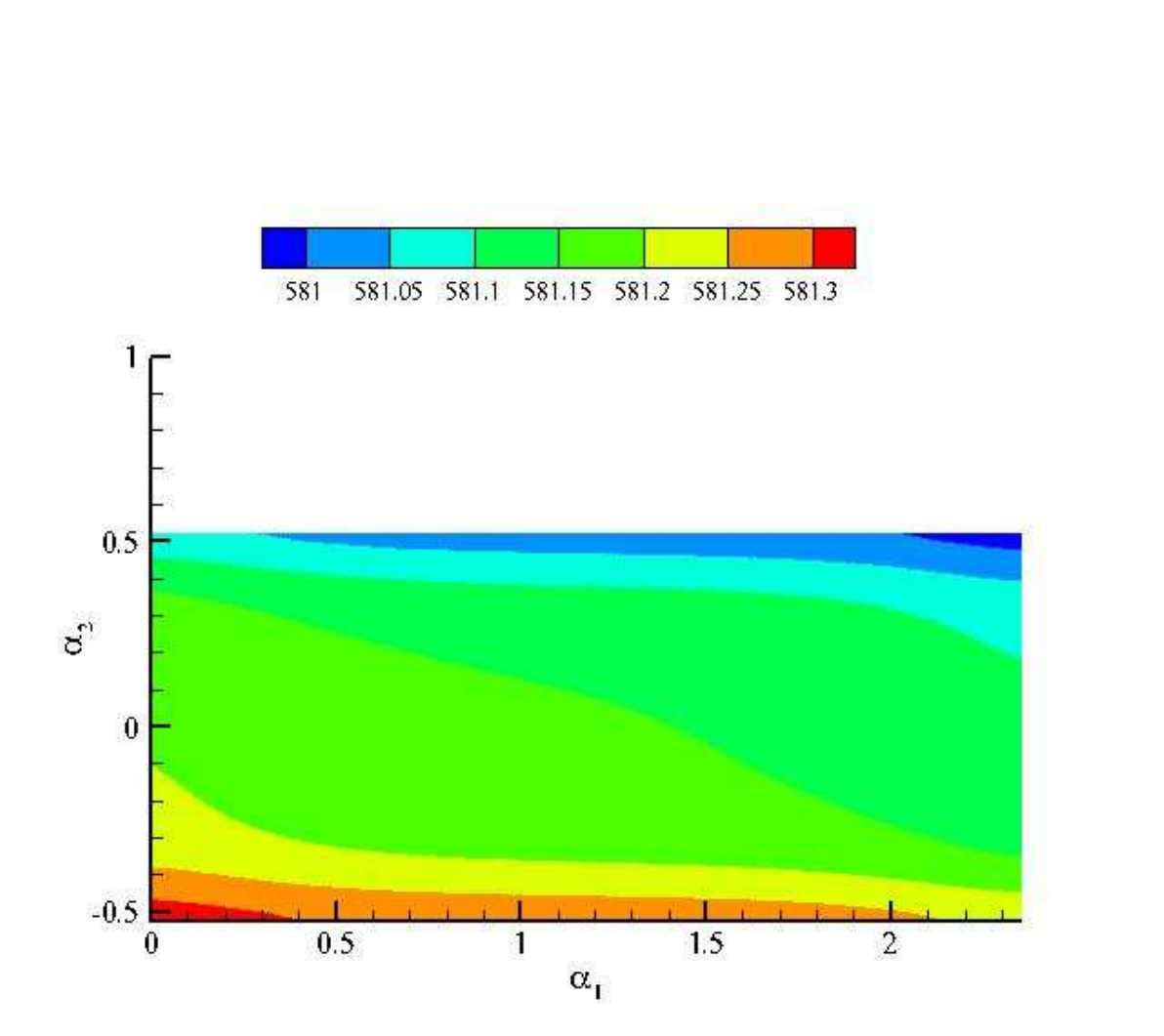}\\
  (a)
\end{minipage}
\begin{minipage}[c]{0.29\textwidth}
  \centering
  \includegraphics[width=37mm]{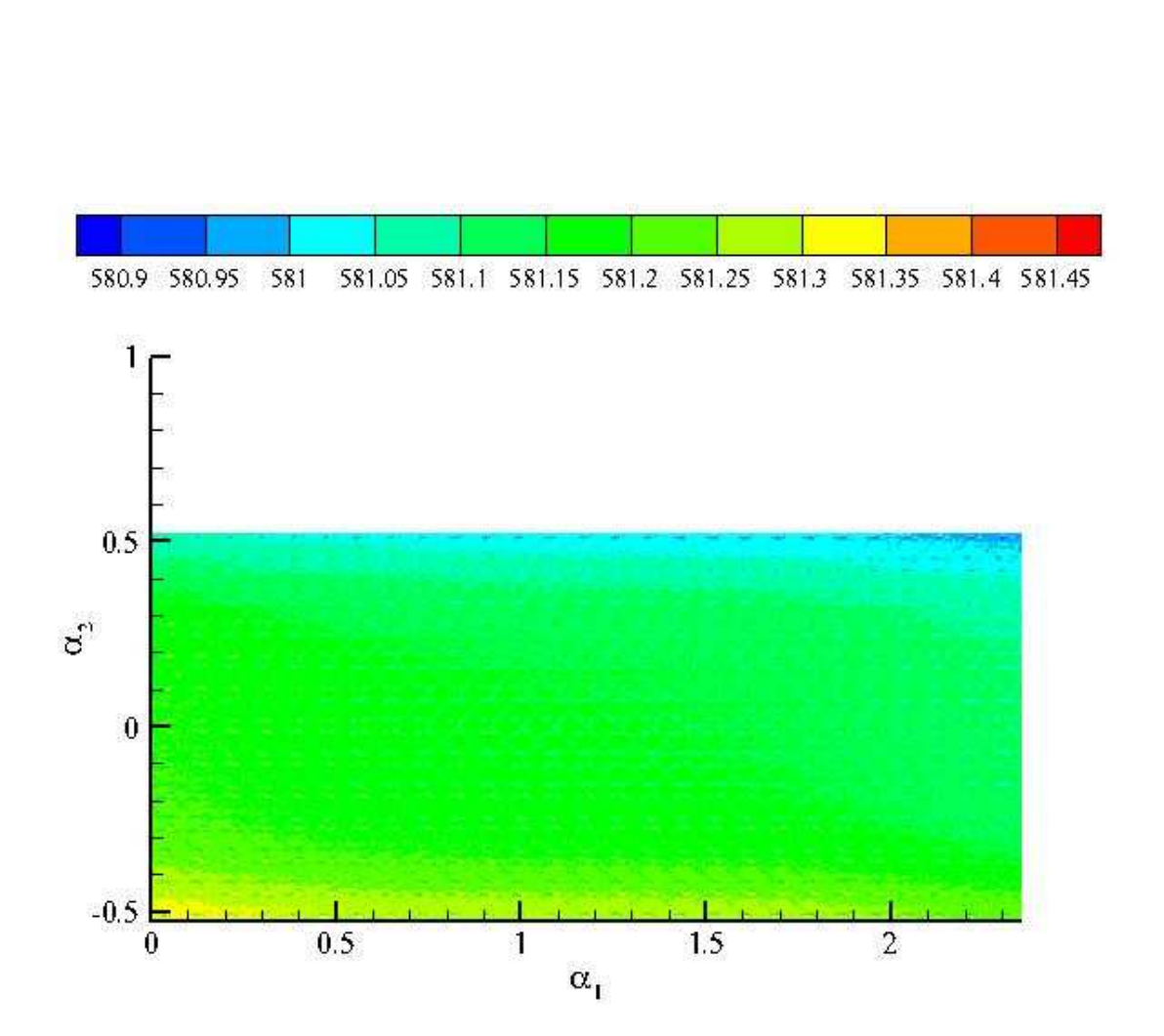}\\
  (b)
\end{minipage}
\begin{minipage}[c]{0.29\textwidth}
  \centering
  \includegraphics[width=37mm]{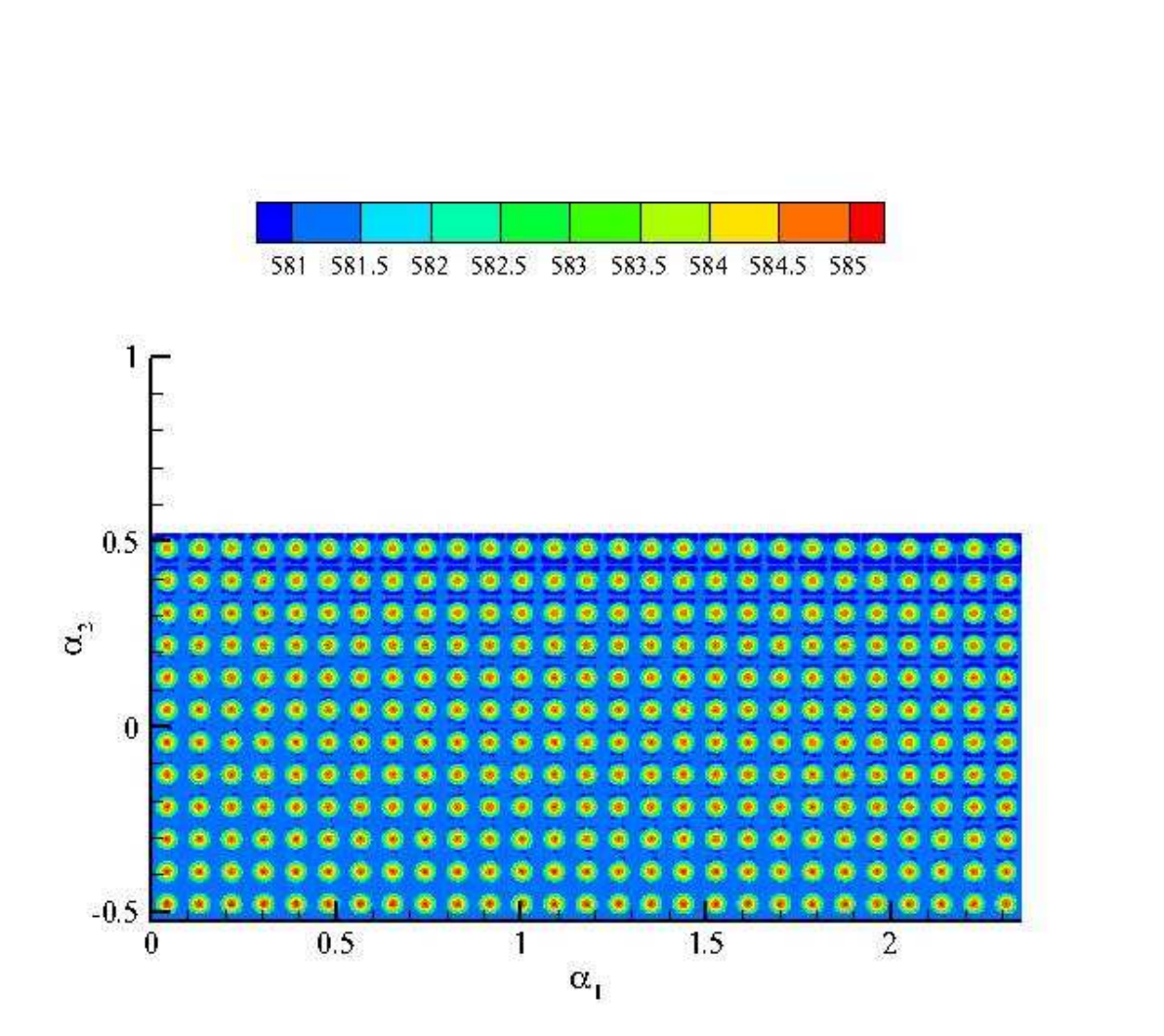}\\
  (c)
\end{minipage}
\caption{The temperature field in cross section $\alpha_3={\pi}/{72}cm$ at time $t$=1.0s: (a) $T^{[0]}$; (b) $T^{[1\xi]}$; (c) $T^{[2\xi]}$.}
\end{figure}
\begin{figure}[!htb]
\centering
\begin{minipage}[c]{0.29\textwidth}
  \centering
  \includegraphics[width=37mm]{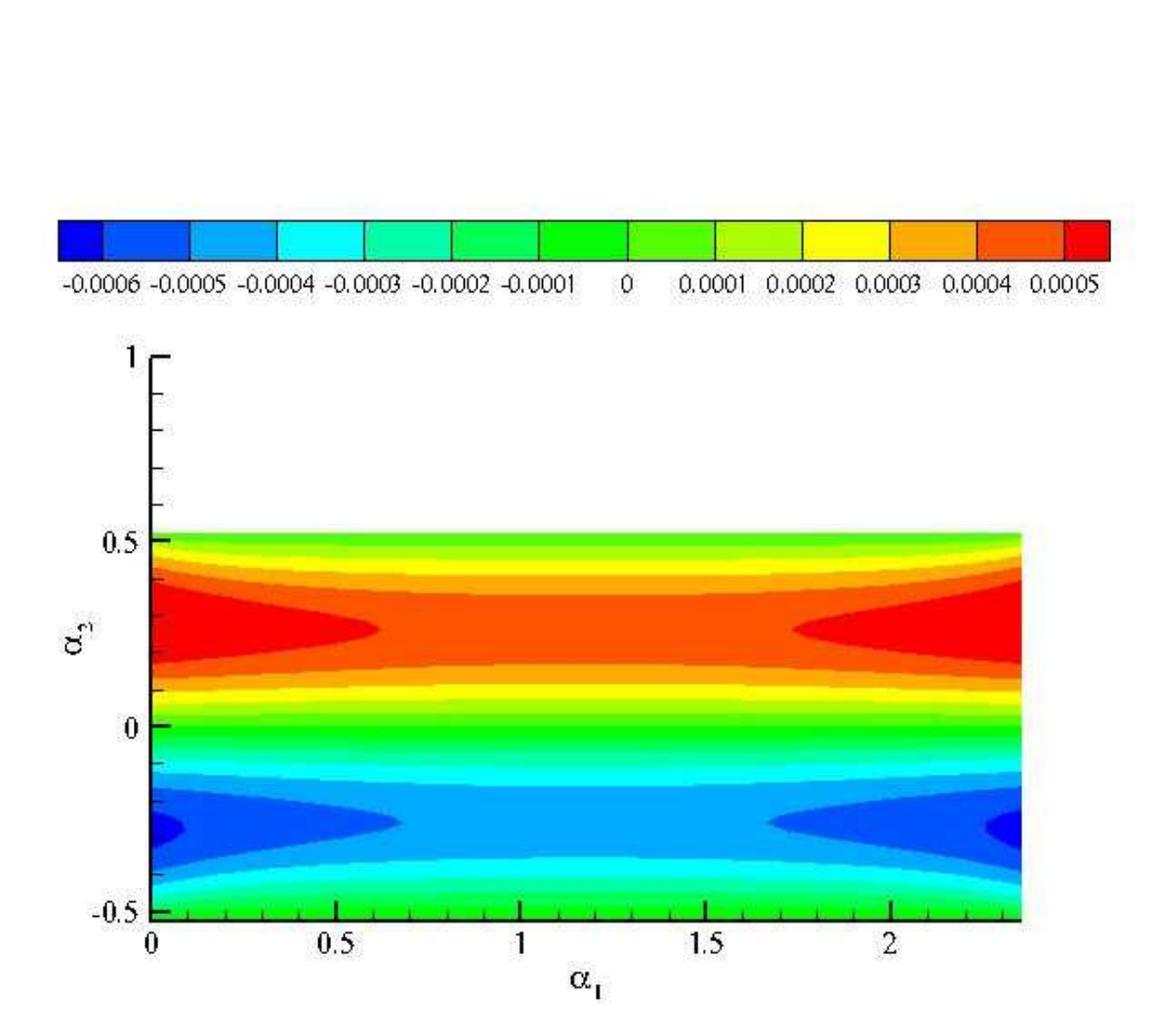}\\
  (a)
\end{minipage}
\begin{minipage}[c]{0.29\textwidth}
  \centering
  \includegraphics[width=37mm]{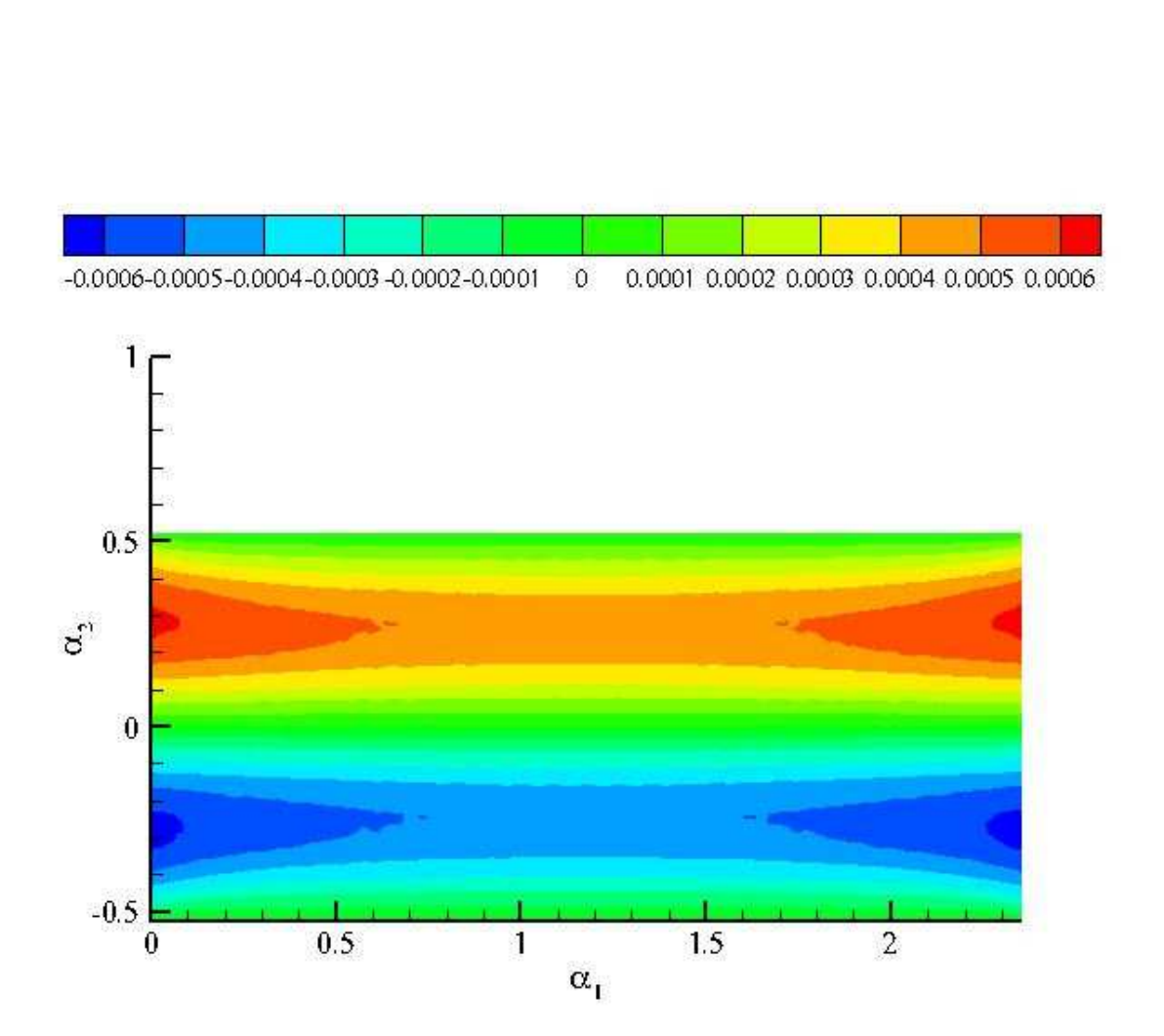}\\
  (b)
\end{minipage}
\begin{minipage}[c]{0.29\textwidth}
  \centering
  \includegraphics[width=37mm]{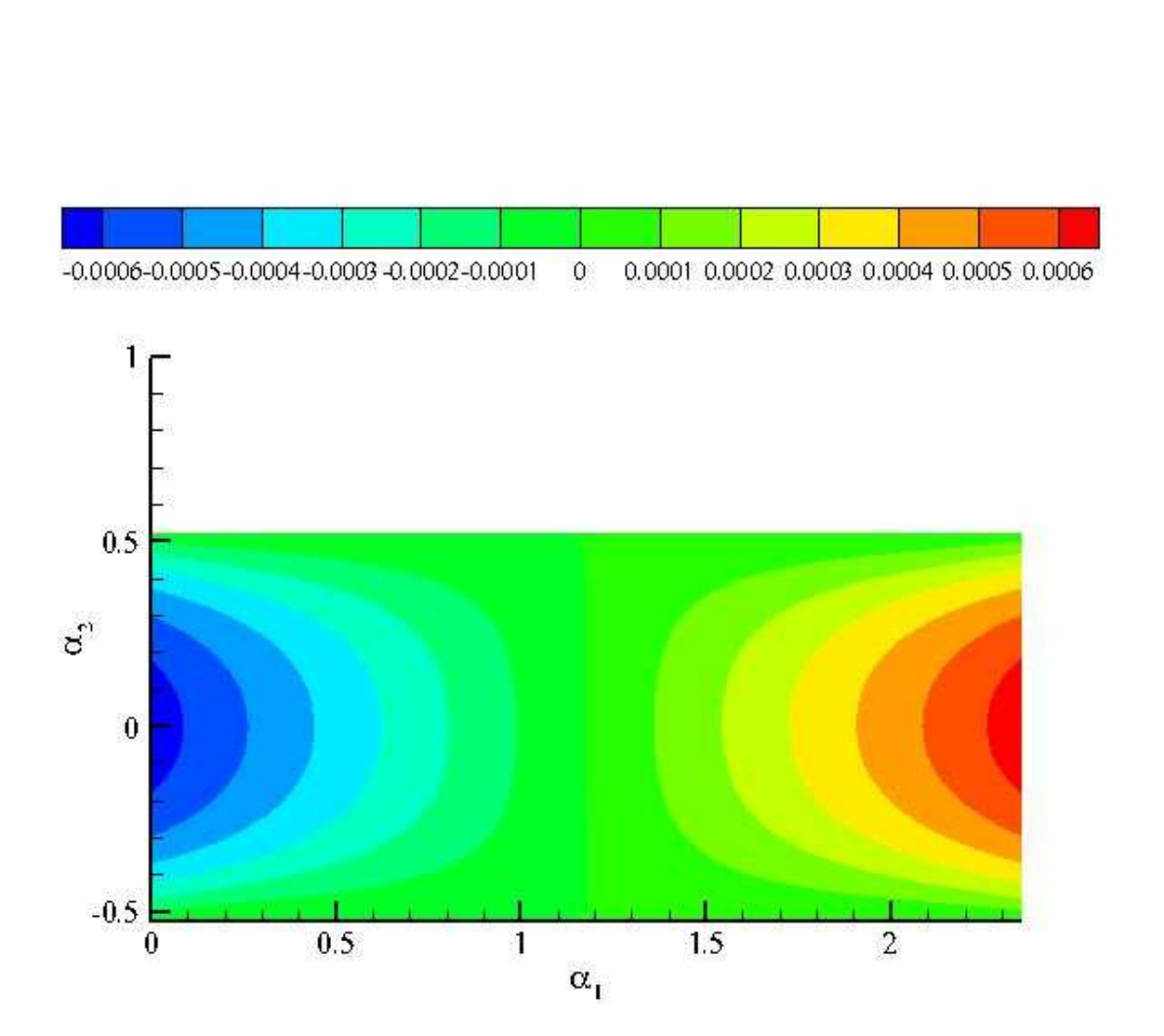}\\
  (c)
\end{minipage}
\caption{The second displacement component in cross section $\alpha_3={\pi}/{72}cm$ at time $t$=1.0s: (a) $u_2^{[0]}$; (b) $u_2^{[1\xi]}$; (c) $u_2^{[2\xi]}$.}
\end{figure}

As illustrated in Table 4, we can clearly see that the computational cost of HOMS approach is far fewer that of the full-scale DNS. By the proposed HOMS approach, a large-scale multi-scale problem is decomposed into a series of small-scale single-scale problems. According to the computational results reported in Figs.\hspace{1mm}13-16, one can discover that merely HOMS solutions precisely catch the steeply microscopic oscillation of inhomogeneous cylindrical shell, and the computational accuracy of homogenized and LOMS solutions is not adequate especially for the temperature field. In addition, it can be summarized that the proposed space-time multi-scale algorithm is stable after longstanding numerical simulation by analyzing the numerical solutions at $t=0.2$s and $t=1.0$s. Accordingly, only HOMS solutions can be employed to high-accuracy computation for heat flux and strains of heterogeneous cylindrical shell, and the HOMS approach should be adopted to compute the practical nonlinear thermo-mechanical problem of inhomogeneous cylindrical shell.
\subsection{Multi-scale nonlinear simulation of heterogeneous doubly-curved shell}
This example considers the nonlinear thermo-mechanical coupling behaviors of a heterogeneous doubly-curved shell with $\displaystyle\xi=\pi/108$ as displayed in Fig.\hspace{1mm}17(a) and Fig.\hspace{1mm}17(b). The Lam$\rm{\acute{e}}$ coefficients of this doubly-curved shell are $H_1=R_1+\alpha_3$, $H_2=R_2+\alpha_3$ and $H_3=1$, where $R_1$ and $R_2$ are the curvature radius of shell's middle plane in $\alpha_1$ and $\alpha_2$ directions with $R_1=R_2=\pi$, respectively. The whole domain $\displaystyle\Omega=(\alpha_1,\alpha_2,\alpha_3)=[-\pi/9,\pi/9]\times[-\pi/9,\pi/9]\times[-\pi/54,\pi/54]cm$  and PUC $\Theta$ are illustrated in Fig.\hspace{1mm}17(c) and Fig.\hspace{1mm}17(d).
\begin{figure}[!htb]
\centering
\begin{minipage}[c]{0.46\textwidth}
  \centering
  \includegraphics[width=0.6\linewidth,totalheight=1.3in]{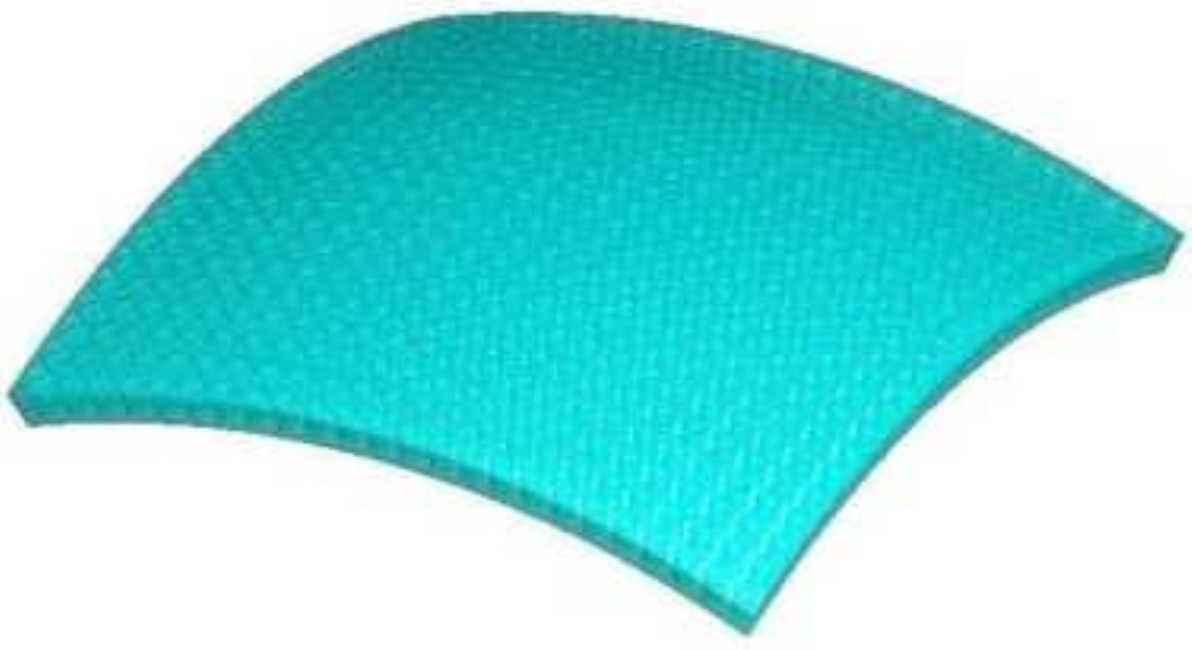} \\
  (a)
\end{minipage}
\begin{minipage}[c]{0.46\textwidth}
  \centering
  \includegraphics[width=0.6\linewidth,totalheight=1.3in]{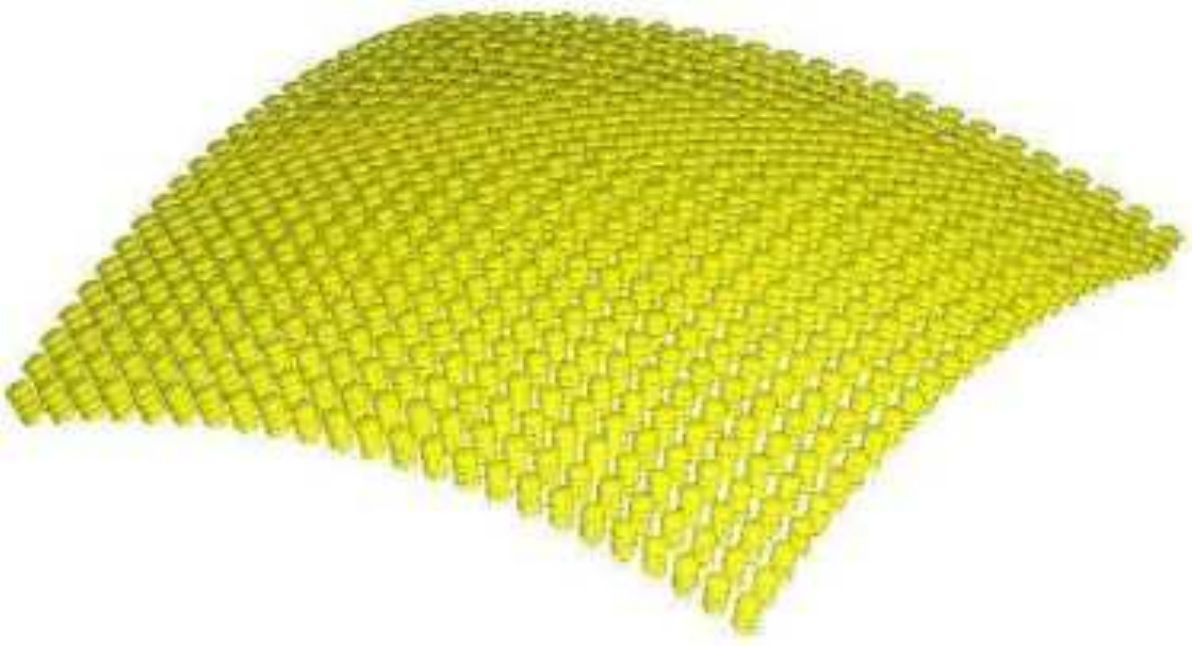} \\
  (b)
\end{minipage}
\begin{minipage}[c]{0.46\textwidth}
  \centering
  \includegraphics[width=0.6\linewidth,totalheight=1.3in]{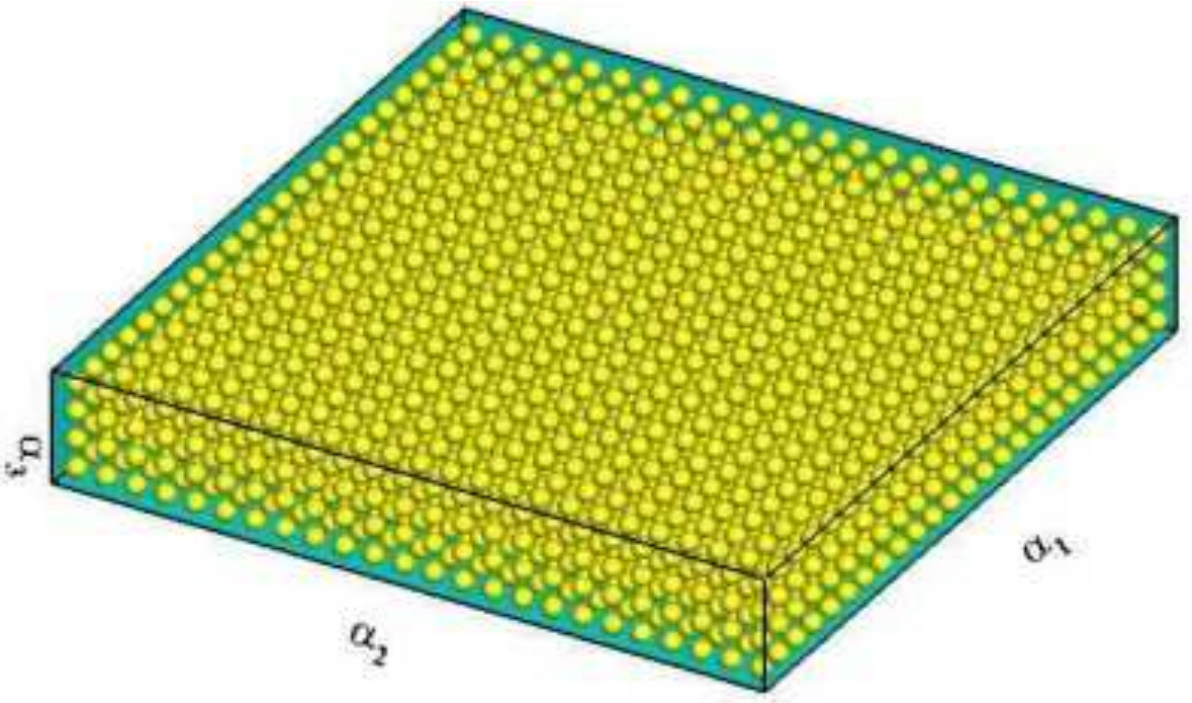} \\
  (c)
\end{minipage}
\begin{minipage}[c]{0.46\textwidth}
  \centering
  \includegraphics[width=0.6\linewidth,totalheight=1.3in]{EX4Nn-new-eps-converted-to.pdf} \\
  (d)
\end{minipage}
\caption{(a) The macrostructure of heterogeneous doubly-curved shell; (b) inclusion distribution in heterogeneous doubly-curved shell; (c) the whole domain $\Omega$; (d) the unit cell $\Theta$.}
\end{figure}

Due to the same reason as second example, we still can not obtain the reference solutions $T_{Fe}^\xi(\bm{\alpha},t)$ and $\bm{u}_{Fe}^\xi(\bm{\alpha},t)$ for the investigated heterogeneous doubly-curved shell. Moreover, its material parameters are defined with the same values as first example. The investigated heterogeneous doubly-curved shell $\Omega$ is clamped on its four side surfaces and the initial temperature is kept at $373.15K$ on its bottom surface. Additionally, we set $h = 100000J/(cm^3\bm{\cdot}s)$ and $(f_1,f_2,f_3) = (0, 0,-120000)N/cm^3$. Next, we create the finite element mesh for original multi-scale equations, auxiliary cell equations and corresponding homogenized equations. The specific mesh information are presented in Table 5.
\begin{table}[!htb]{\caption{Summary of computational cost.}}
\centering
\begin{tabular}{cccc}
\hline
 & Multi-scale eqs. & Cell eqs. & Homogenized eqs. \\
\hline
FEM elements & $\approx$2304$\times$75466 & 75466 & 884736\\
FEM nodes    & $\approx$2304$\times$13062 & 13062 & 159953\\
\hline
\end{tabular}
\end{table}

In this example, off-line 17500 times computation is required for auxiliary cell problems totally, where the quantity of first-order and second-order auxiliary cell functions is set as 33 and 250 respectively. Moreover, we distribute 5 uniformly-spaced interpolating points toward the shell thickness direction and 14 macroscopic interpolation temperature in one unit cell. The nonlinear thermo-mechanical coupling responses of the heterogeneous doubly-curved shell are simulated in the time interval $t\in[0,1]$s. Setting the temporal step as $\Delta t = 0.02$s, the multi-scale equations (1.1) and the homogenized equations (2.18) are on-line simulated respectively. Next, the final simulation results of temperature and displacement fields at $t=0.1$s and $t=1.0$s are depicted in Figs.\hspace{1mm}18 and 19, and Figs.\hspace{1mm}20 and 21, respectively.
\begin{figure}[!htb]
\centering
\begin{minipage}[c]{0.29\textwidth}
  \centering
  \includegraphics[width=37mm]{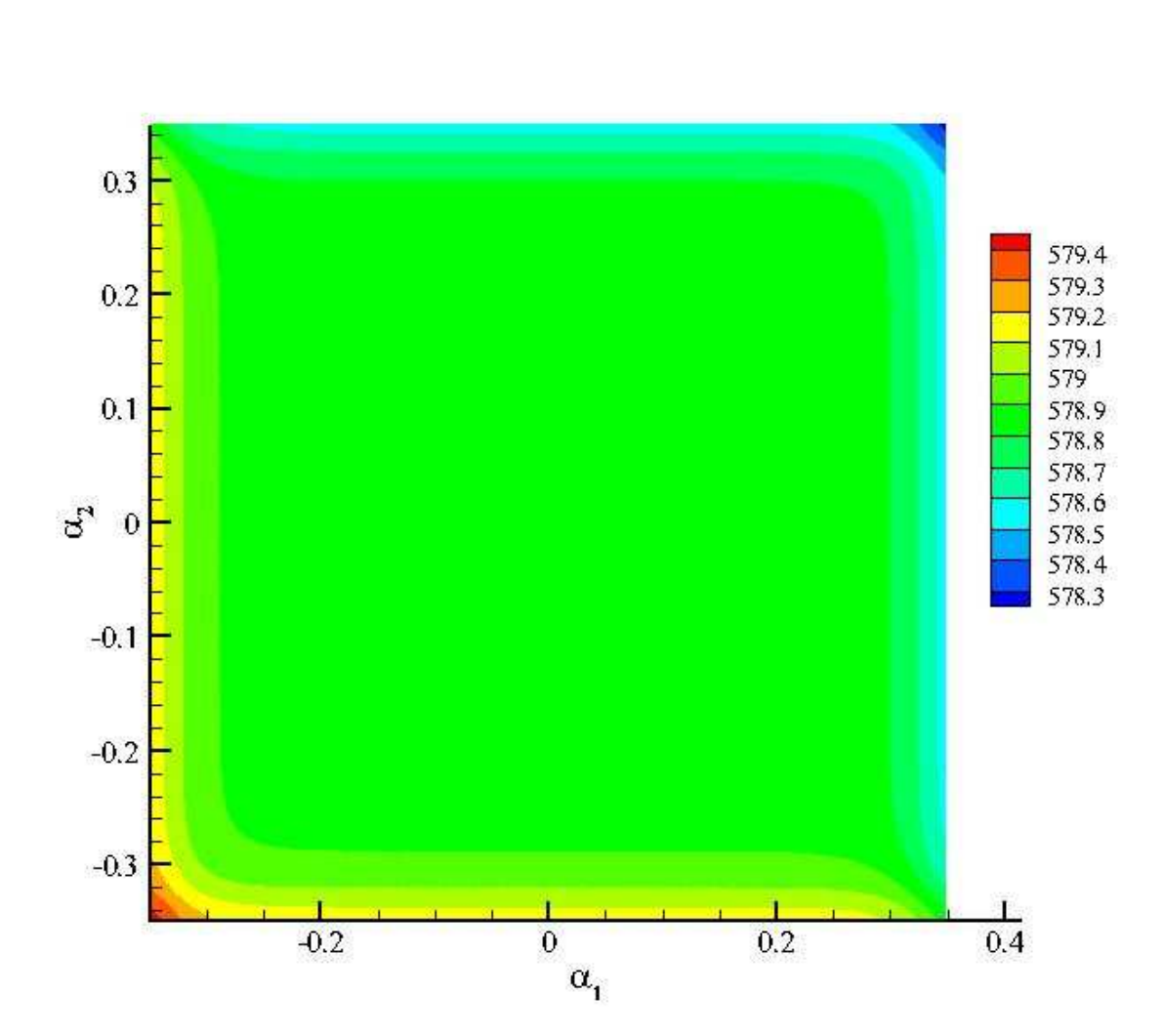}\\
  (a)
\end{minipage}
\begin{minipage}[c]{0.29\textwidth}
  \centering
  \includegraphics[width=37mm]{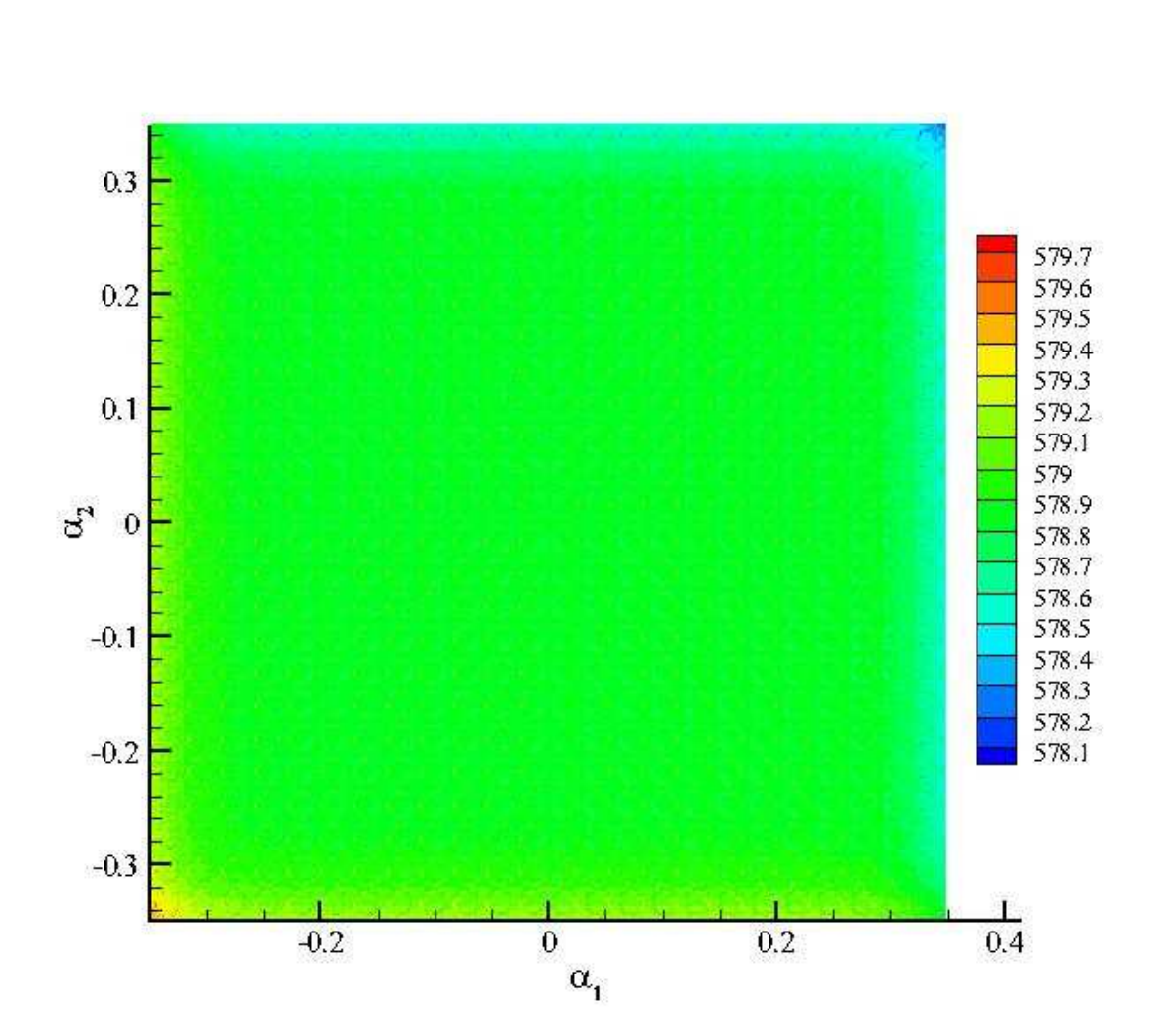}\\
  (b)
\end{minipage}
\begin{minipage}[c]{0.29\textwidth}
  \centering
  \includegraphics[width=37mm]{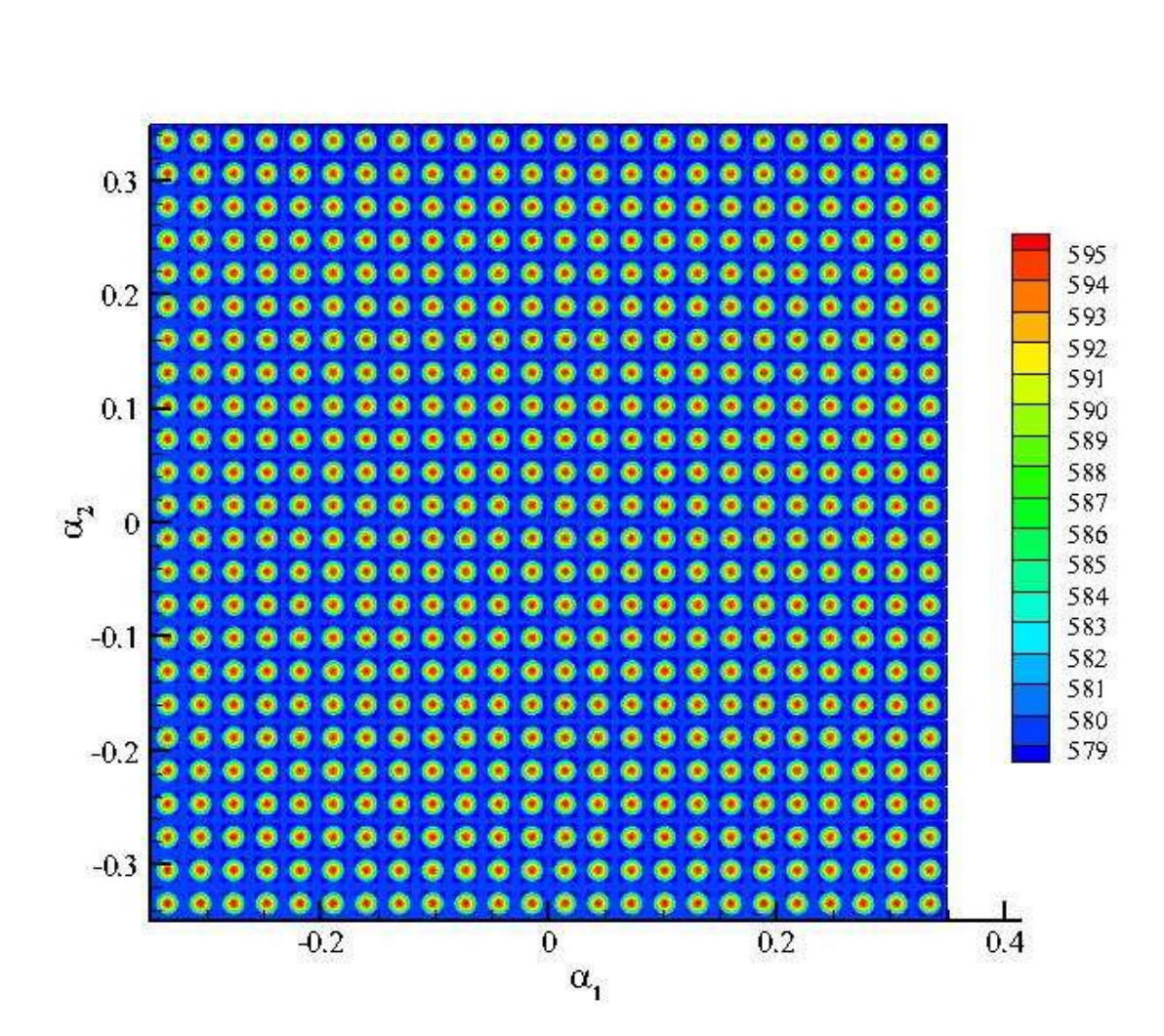}\\
  (c)
\end{minipage}
\caption{The temperature field in cross section $\alpha_3={\pi}/{216}cm$ at time $t$=0.1s: (a) $T^{[0]}$; (b) $T^{[1\xi]}$; (c) $T^{[2\xi]}$.}
\end{figure}
\begin{figure}[!htb]
\centering
\begin{minipage}[c]{0.29\textwidth}
  \centering
  \includegraphics[width=37mm]{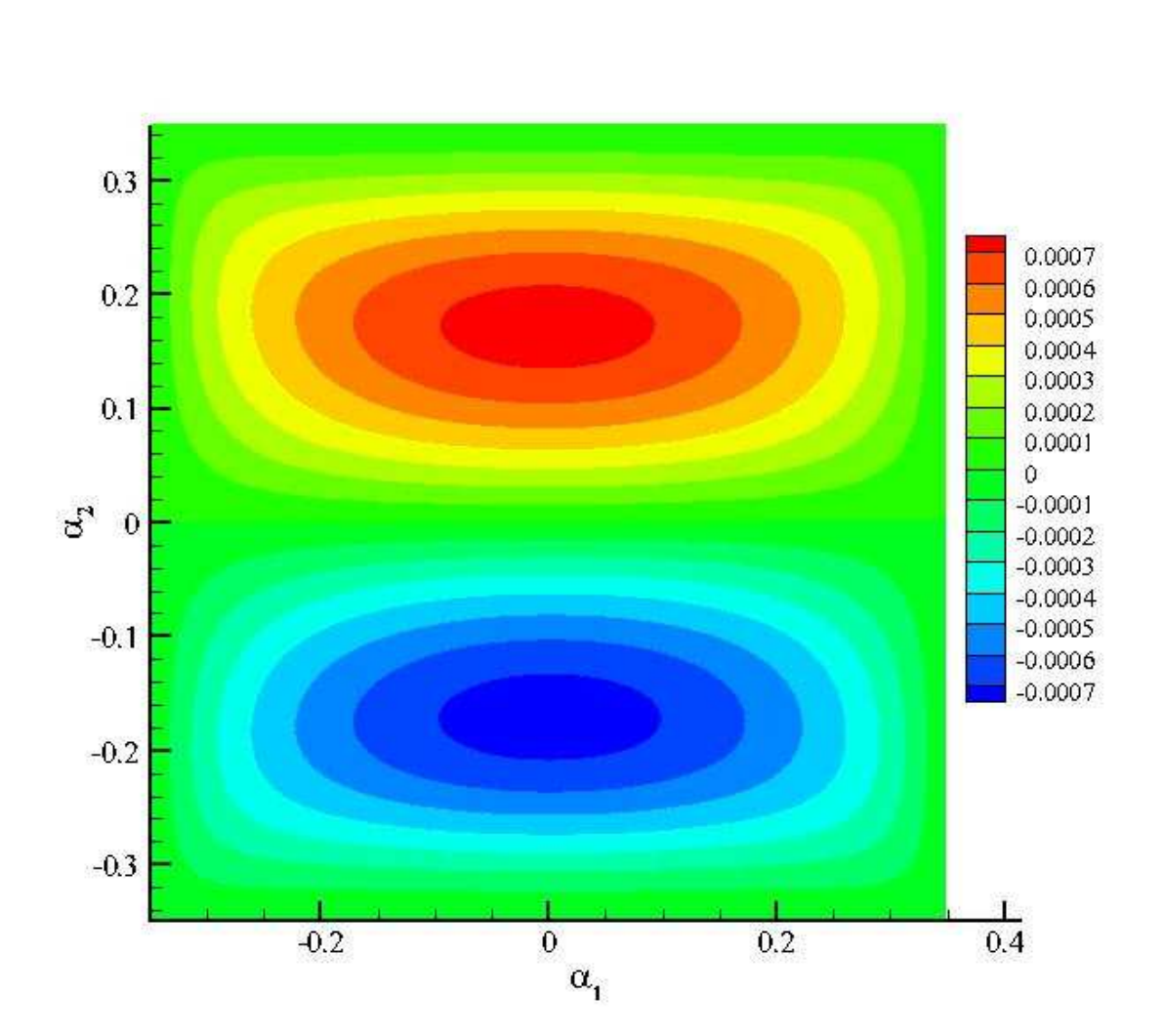}\\
  (a)
\end{minipage}
\begin{minipage}[c]{0.29\textwidth}
  \centering
  \includegraphics[width=37mm]{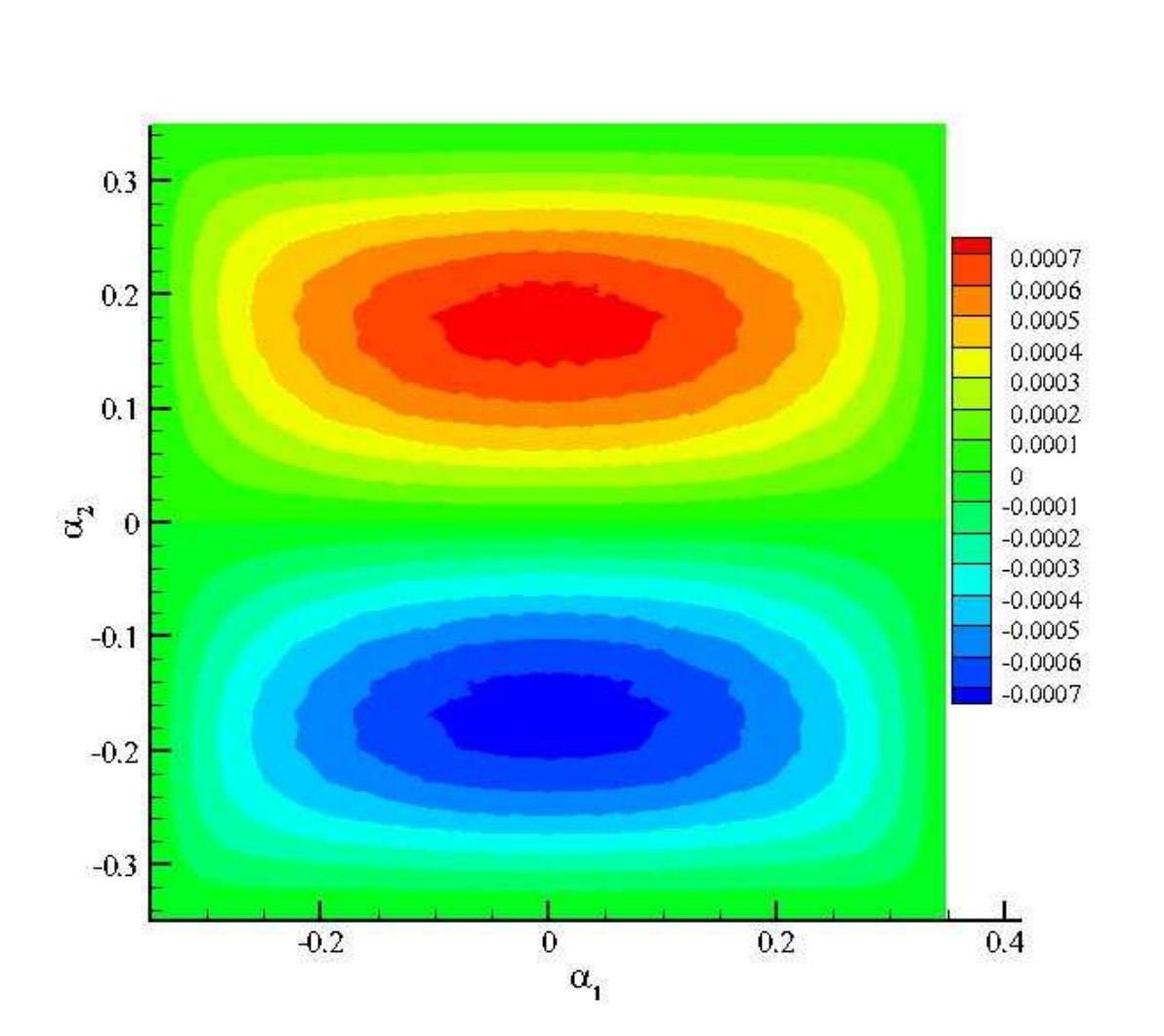}\\
  (b)
\end{minipage}
\begin{minipage}[c]{0.29\textwidth}
  \centering
  \includegraphics[width=37mm]{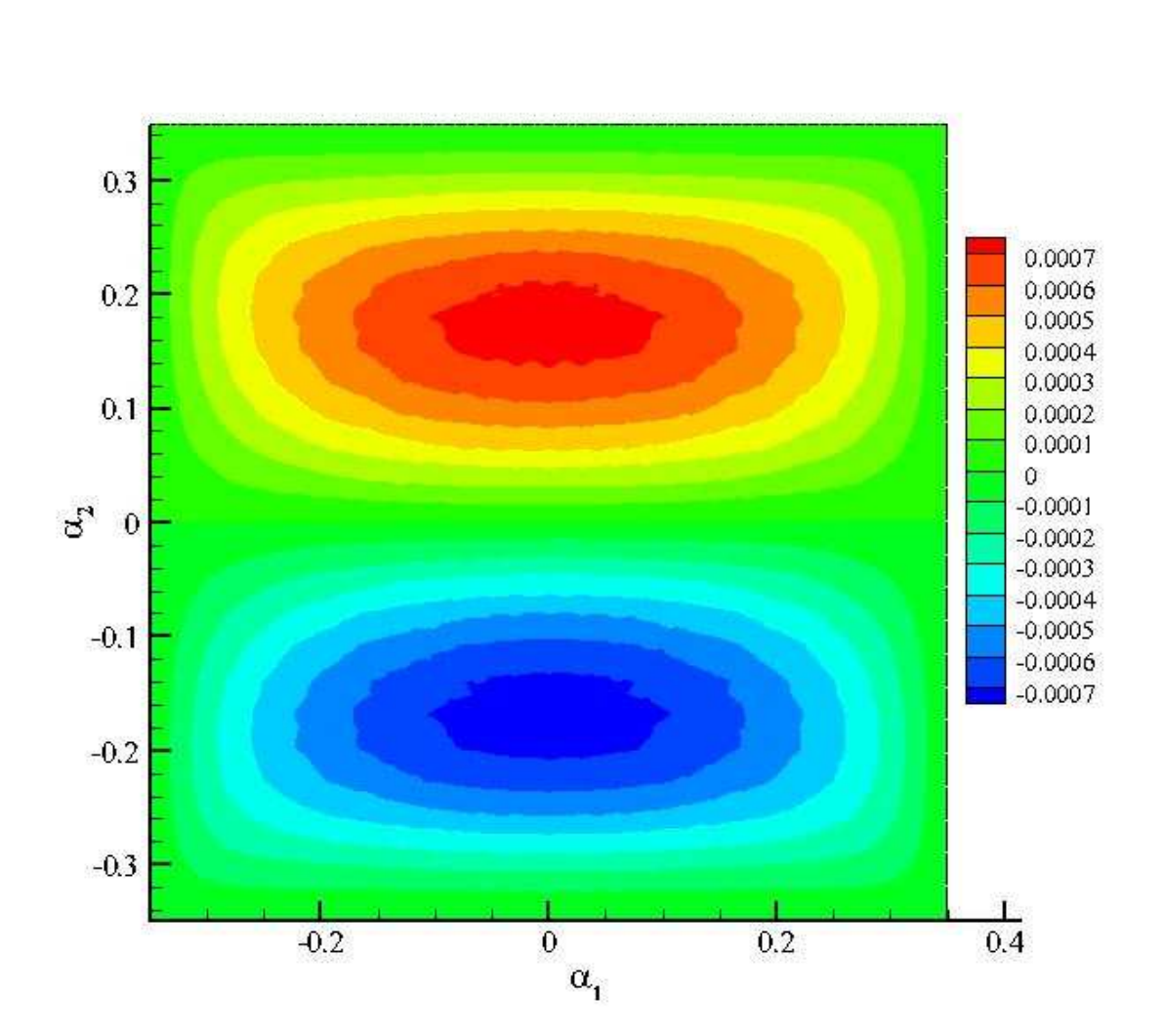}\\
  (c)
\end{minipage}
\caption{The second displacement component in cross section $\alpha_3={\pi}/{216}cm$ at time $t$=0.1s: (a) $u_2^{[0]}$; (b) $u_2^{[1\xi]}$; (c) $u_2^{[2\xi]}$.}
\end{figure}
\begin{figure}[!htb]
\centering
\begin{minipage}[c]{0.29\textwidth}
  \centering
  \includegraphics[width=37mm]{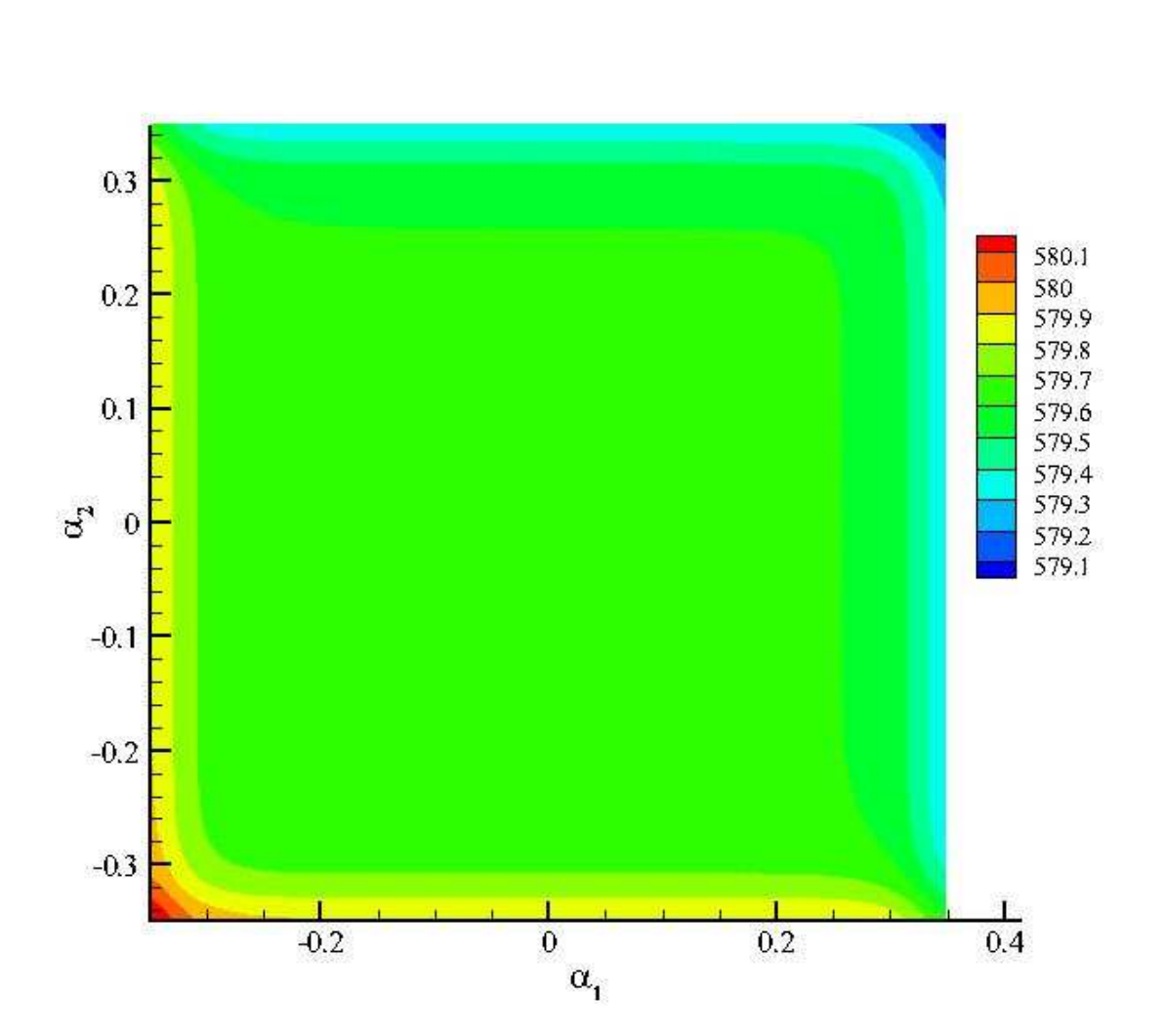}\\
  (a)
\end{minipage}
\begin{minipage}[c]{0.29\textwidth}
  \centering
  \includegraphics[width=37mm]{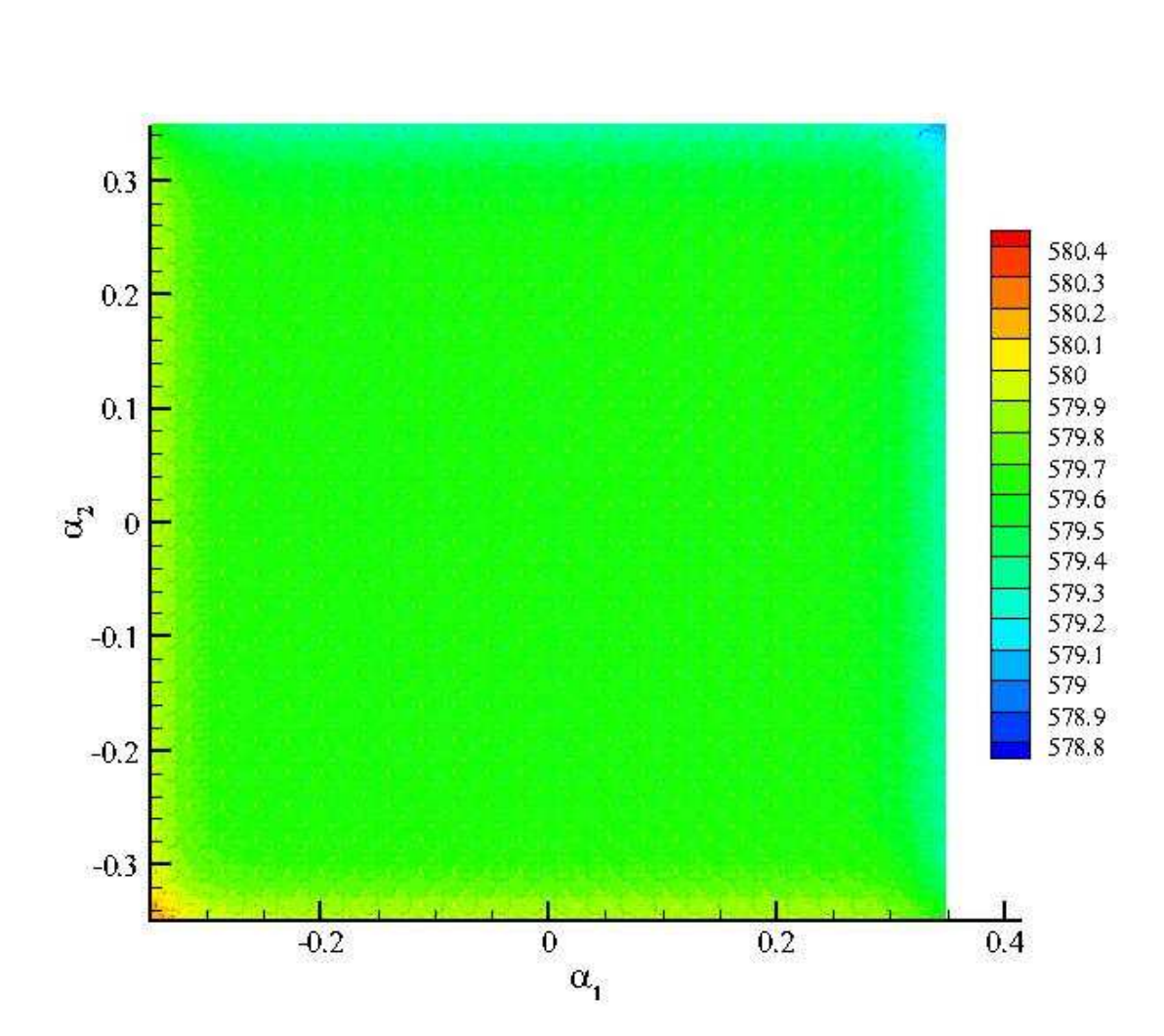}\\
  (b)
\end{minipage}
\begin{minipage}[c]{0.29\textwidth}
  \centering
  \includegraphics[width=37mm]{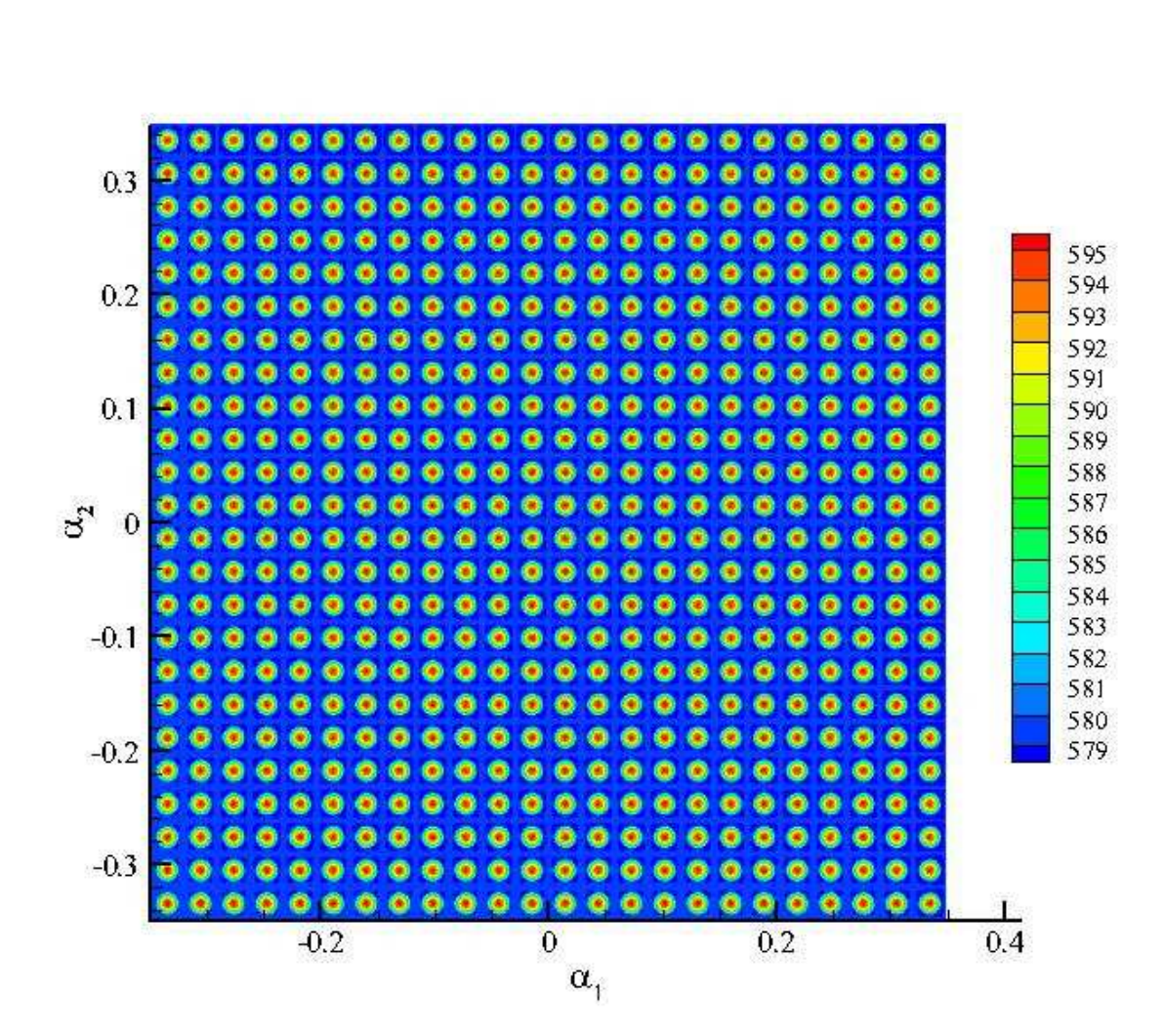}\\
  (c)
\end{minipage}
\caption{The temperature field in cross section $\alpha_3={\pi}/{216}cm$ at time $t$=1.0s: (a) $T^{[0]}$; (b) $T^{[1\xi]}$; (c) $T^{[2\xi]}$.}
\end{figure}
\begin{figure}[!htb]
\centering
\begin{minipage}[c]{0.29\textwidth}
  \centering
  \includegraphics[width=37mm]{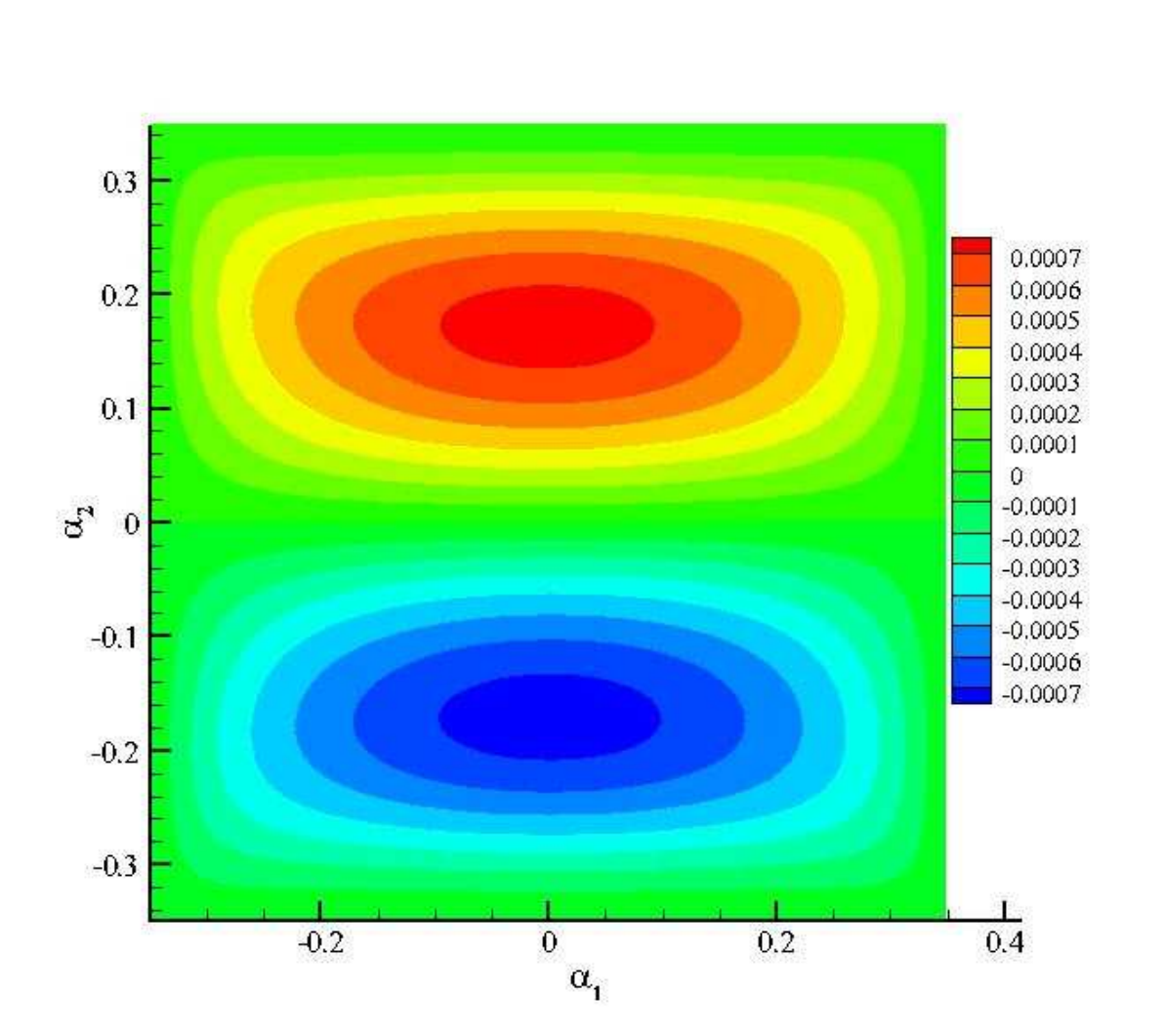}\\
  (a)
\end{minipage}
\begin{minipage}[c]{0.29\textwidth}
  \centering
  \includegraphics[width=37mm]{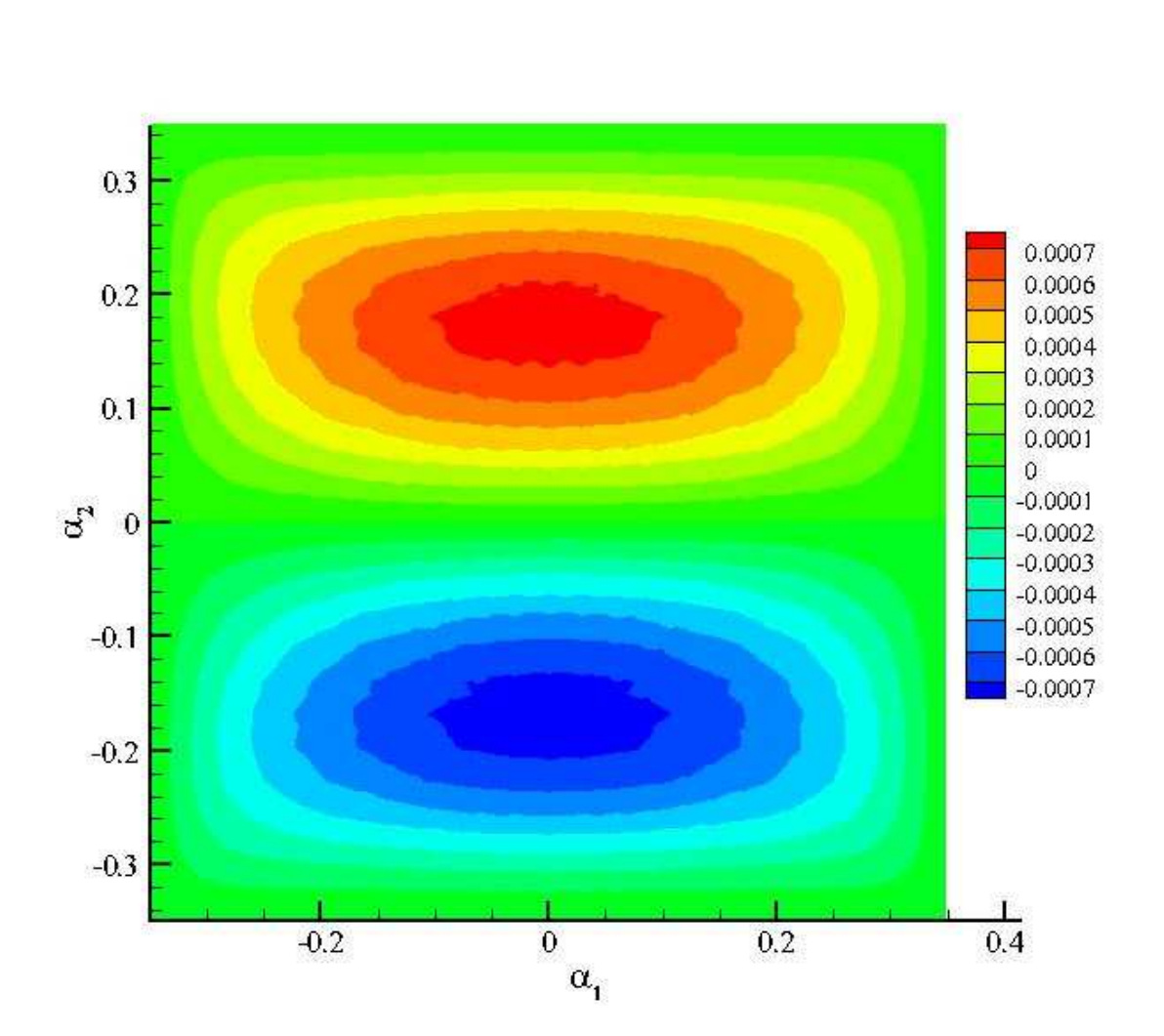}\\
  (b)
\end{minipage}
\begin{minipage}[c]{0.29\textwidth}
  \centering
  \includegraphics[width=37mm]{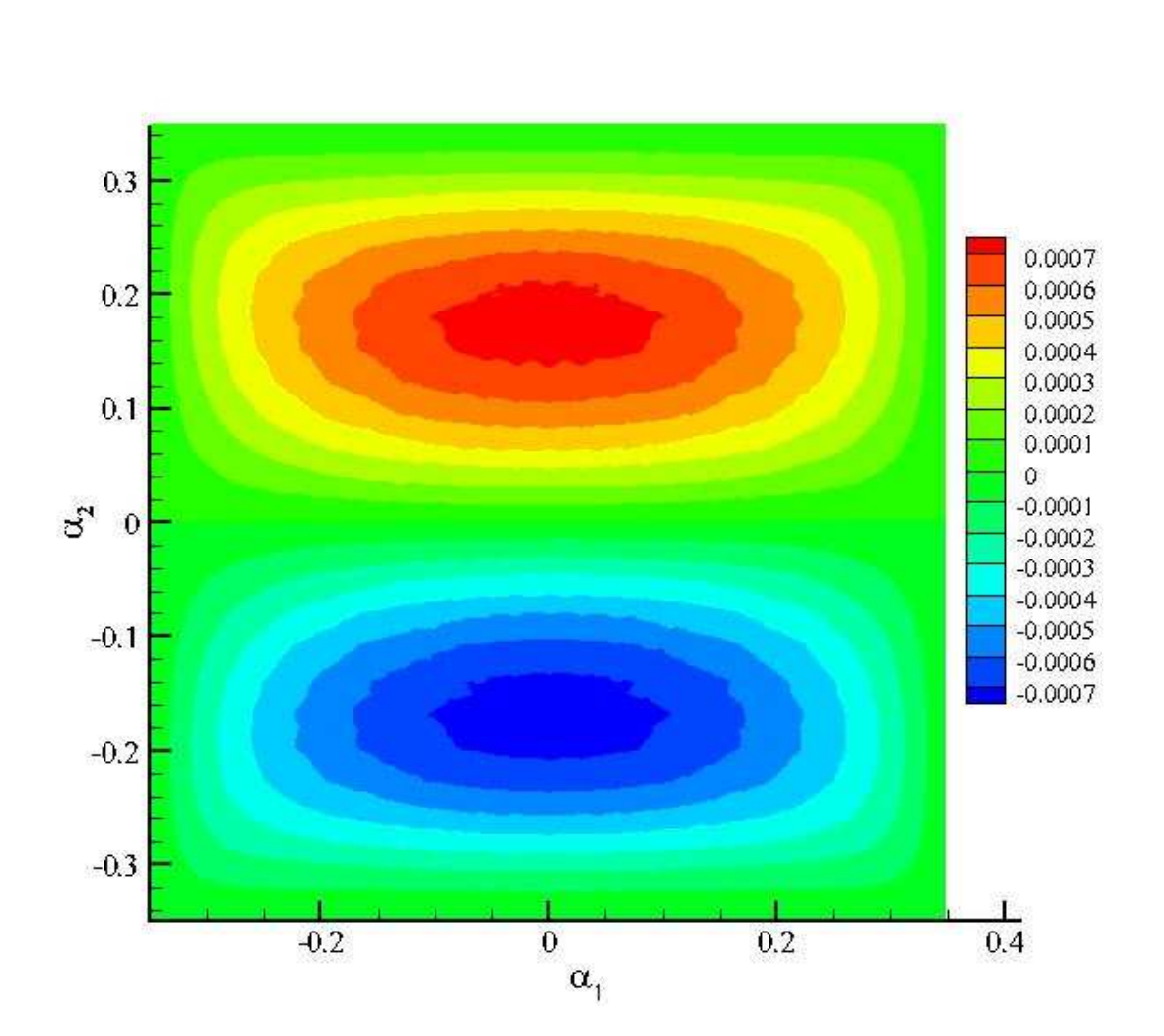}\\
  (c)
\end{minipage}
\caption{The second displacement component in cross section $\alpha_3={\pi}/{216}cm$ at time $t$=1.0s: (a) $u_2^{[0]}$; (b) $u_2^{[1\xi]}$; (c) $u_2^{[2\xi]}$.}
\end{figure}

As illustrated in Table 5, we can conclude that our HOMS approach can remarkably economize computing resources by comparison with high-resolution DNS. From the computational results in Figs.\hspace{1mm}18-21, it can be clearly found that only HOMS solutions have the capacity to accurately simulate the nonlinear thermo-mechanical coupling behaviors of heterogeneous doubly-curved shell, and the computational accuracy of homogenized and LOMS solutions is significantly inferior particular for the temperature field. For heterogeneous doubly-curved shell, the homogenized method can only catch its macroscopic behaviors and the LOMS method can only capture its inadequate microscopic responses. What's more, it can be summarized that the proposed HOMS approach can retain long-time numerical stability because the blow-up phenomenon does not appear in the numerical solutions from $t=0.1$s to $t=1.0$s. Thus, the HOMS solutions would be utilized to catch the microscopic coupling responses of inhomogeneous doubly-curved shell.
\section{Conclusions and prospects}
Extreme thermo-mechanical coupled environments incur complicated nonlinear responses of heterogeneous shells in aerospace, aviation and other engineering fields. However, highly inhomogeneous and temperature-dependent properties of heterogeneous shells are laborious to handle by traditional numerical methods. In the present study, we first propose, to our best knowledge, a HOMS computational methodology for accurately and efficiently computing time-dependent nonlinear thermo-mechanical equations of inhomogeneous shells possessing orthogonal periodic configurations. The principal innovations of this study are threefold: to establish the HOMS computational model, to attain the local and global error estimations for the multi-scale solutions, and to present the corresponding space-time numerical algorithm. Numerical results well demonstrated that the proposed HOMS computational approach is stable and efficient for time-variant nonlinear multi-scale simulation of heterogeneous shells, that can be employed to investigate the nonlinear thermal and mechanical coupling behaviors of inhomogeneous block, plate, cylindrical shell and doubly-curved shell. Moreover, computational results indicate that only HOMS solutions can precisely catch the microscopic oscillatory responses and offer high-accuracy numerical solutions for practical applications, that strongly confirm the theoretical conclusions of this study. From the computation point of view, numerical experiments also expressly manifest the proposed HOMS approach significantly reduces the number of computational freedom degrees and memory storage space, comparatively to direct numerical simulation on large-scale heterogeneous shells.

In the future research, the extension of the proposed HOMS computational approach to nonlinear thermo-mechanical analysis with thermal convection and radiation effects of porous heterogeneous shells should be involved to implement more realistic engineering simulations. Also, we expect the extensions of the presented HOMS computational framework to random thermo-mechanical problems, thermo-elasto-plastic problems and thermo-visco-elastic problems of nonlinear heterogeneous shells.
\appendix
\setcounter{equation}{0}
\renewcommand{\theequation}{A.\arabic{equation}}\\
\section*{Appendix A. Mathematical proof of continuous properties of microscopic cell functions}\\
In order to resolve this issue, employing the similar idea in \cite{R43}, it can be demonstrated that all auxiliary cell functions are continuous with regard to macroscopic parameters. Taking the auxiliary cell function ${M^m}(\bm{H},T^{[0]},{\bm{\beta}})$ for example, we first denote the ${M^m_{\bm{\alpha}_1}}$ and ${M^m_{\bm{\alpha}_2}}$ as the auxiliary cell function ${M^m}(\bm{H},T^{[0]},{\bm{\beta }})$ in points $\bm{\alpha}_1$ and $\bm{\alpha}_2$ of $\Omega$ separately. Afterwards, the variational equations for ${M^m_{\bm{\alpha}_1}}$ and ${M^m_{\bm{\alpha}_2}}$ are established on the basis of auxiliary cell problem (2.17).
\begin{equation}
-\int_{\Theta}{k_{ij}^{[0]}}(\bm{\alpha}_1){{\widetilde \Psi }_j}({M^m_{\bm{\alpha}_1}}){\widetilde \Psi }_i( \upsilon^{h_1})d\Theta=\int_{\Theta} k_{im}^{[0]}(\bm{\alpha}_1){\widetilde \Psi }_i( \upsilon^{h_1})d\Theta,\;\forall\upsilon^{h_1}\in V_{h_1}(\Theta).
\end{equation}
\begin{equation}
-\int_{\Theta}{k_{ij}^{[0]}}(\bm{\alpha}_2){{\widetilde \Psi }_j}({M^m_{\bm{\alpha}_2}}){\widetilde \Psi }_i( \upsilon^{h_1})d\Theta=\int_{\Theta} k_{im}^{[0]}(\bm{\alpha}_2){\widetilde \Psi }_i( \upsilon^{h_1})d\Theta,\;\forall\upsilon^{h_1}\in V_{h_1}(\Theta).
\end{equation}
After that, employing the variational equation of ${M^m_{\bm{\alpha}_1}}$ to subtract the variational equation of ${M^m_{\bm{\alpha}_2}}$, a new variational equation can be derived as below.
\begin{equation}
\begin{aligned}
&\int_{\Theta}{k_{ij}^{[0]}}(\bm{\alpha}_1){{\widetilde \Psi }_j}({M^m_{\bm{\alpha}_1}}-{M^m_{\bm{\alpha}_2}}){\widetilde \Psi }_i( \upsilon^{h_1})d\Theta\\
&=-\int_{\Theta}\big[{k_{ij}^{[0]}}(\bm{\alpha}_1)-{k_{ij}^{[0]}}(\bm{\alpha}_2)\big]{{\widetilde \Psi }_j}({M^m_{\bm{\alpha}_2}}){\widetilde \Psi }_i( \upsilon^{h_1})d\Theta\\
&-\int_{\Theta}\big[k_{im}^{[0]}(\bm{\alpha}_1)-k_{im}^{[0]}(\bm{\alpha}_2)\big]{\widetilde \Psi }_i( \upsilon^{h_1})d\Theta,\;\forall\upsilon^{h_1}\in V_{h_1}(\Theta).
\end{aligned}
\end{equation}
Next, if $|{k_{ij}^{[0]}}(\bm{\alpha}_1)-{k_{ij}^{[0]}}(\bm{\alpha}_2)|\leq C|\bm{\alpha}_1-\bm{\alpha}_2|$ and we replace $\upsilon^{h_1}$ with ${M^m_{\bm{\alpha}_1}}-{M^m_{\bm{\alpha}_2}}$, it follows that
\begin{equation}
\begin{aligned}
&\;\;\;\;\gamma_0\left \|{M^m_{\bm{\alpha}_1}}-{M^m_{\bm{\alpha}_2}}\right\|_{H^1_0(\Omega)}^2\\
&\leq\int_{\Theta}{k_{ij}^{[0]}}(\bm{\alpha}_1){{\widetilde \Psi }_j}({M^m_{\bm{\alpha}_1}}-{M^m_{\bm{\alpha}_2}}){\widetilde \Psi }_i({M^m_{\bm{\alpha}_1}}-{M^m_{\bm{\alpha}_2}})d\Theta\\
&=-\int_{\Theta}\big[{k_{ij}^{[0]}}(\bm{\alpha}_1)-{k_{ij}^{[0]}}(\bm{\alpha}_2)\big]{{\widetilde \Psi }_j}({M^m_{\bm{\alpha}_2}}){\widetilde \Psi }_i({M^m_{\bm{\alpha}_1}}-{M^m_{\bm{\alpha}_2}})d\Theta\\
&-\int_{\Theta}\big[k_{im}^{[0]}(\bm{\alpha}_1)-k_{im}^{[0]}(\bm{\alpha}_2)\big]{\widetilde \Psi }_i( {M^m_{\bm{\alpha}_1}}-{M^m_{\bm{\alpha}_2}})d\Theta\\
&\leq C|\bm{\alpha}_1-\bm{\alpha}_2|\left \|{M^m_{\bm{\alpha}_1}}-{M^m_{\bm{\alpha}_2}}\right\|_{H^1_0(\Omega)}
\end{aligned}
\end{equation}
Hence, when we assume $\bm{\alpha}_1\rightarrow\bm{\alpha}_2$, it can deduce ${M^m_{\bm{\alpha}_1}}\rightarrow{M^m_{\bm{\alpha}_2}}$. Finally, we can testify the continuous properties of microscopic cell functions. In the future, we shall further study how to distribute representative material points for obtaining the optimal numerical accuracy. In addition, we hope to further excavate the computational potential of the proposed space-time multi-scale algorithm by introducing parallel algorithms in off-line computation stage to improve the computational efficiency.
\renewcommand{\appendixname}{Appendix~\Alph{section}}

\bibliographystyle{siamplain}
\bibliography{references}
\end{document}